\documentclass[a4paper,reqno]{amsart}

\usepackage[utf8]{inputenc}
\usepackage[T1]{fontenc}
\usepackage{lmodern}
\usepackage[english]{babel}
\usepackage{microtype}
\usepackage{stackengine}

\usepackage[
  backend=biber,
  style=authoryear,
  natbib=true
]{biblatex}

\usepackage{amsmath,amssymb,amsfonts,amsthm}
\usepackage{mathtools,accents, upgreek,xcolor}
\usepackage{mathrsfs,comment,aliascnt,braket,bm,esint}
\usepackage{graphicx} 

\usepackage[colorlinks=true,
            linkcolor=blue,      
            citecolor=green,     
            urlcolor=blue]{hyperref}

\DeclareFieldFormat{citehyperref}{%
  \DeclareFieldAlias{bibhyperref}{noformat}%
  \bibhyperref{#1}}

\savebibmacro{cite}
\renewbibmacro*{cite}{%
  \printtext[citehyperref]{%
    \restorebibmacro{cite}%
    \usebibmacro{cite}}}

\usepackage{etoolbox,comment}
\makeatletter
\patchcmd{\@maketitle}
  {\ifx\@empty\@dedicatory}
  {\ifx\@empty\@date \else {\vskip3ex \centering\footnotesize\@date\par\vskip1ex}\fi
   \ifx\@empty\@dedicatory}
  {}{}
\patchcmd{\@adminfootnotes}
  {\ifx\@empty\@date\else \@footnotetext{\@setdate}\fi}
  {}{}{}
\makeatother

\newcommand{\dist}{\mathsf{d}}

\newcommand{\vol}{\mathrm{vol}}

\theoremstyle{theorem}
\newtheorem{theorem}{Theorem}

\theoremstyle{definition}
\newtheorem{bump}{Bump}[section]
\newtheorem{definition}{Definition}[section]

\newtheorem{lemma}[bump]{Lemma}

\newtheorem{prop}[bump]{Proposition}
\newtheorem{rem}[bump]{Remark}

\theoremstyle{remark}

\DeclareMathOperator{\Lip}{Lip}

\DeclareMathOperator{\inj}{inj}

\newcommand{\euc}{\text{euc}}

\makeatletter

\makeatletter

\allowdisplaybreaks

\mathtoolsset{showonlyrefs}  
\title{$L^p$ convergence of AMV Laplacians on completions of Riemannian manifolds}

\author{Manuel Dias}
\address{M.~Dias: Department of Mathematics and Data Science, Vrije Universiteit Brussel, Pleinlaan 2, B-1050 Elsene, Belgium.}
\email{manuel.dias@vub.be}

\date{\today}

\begin{document}

\maketitle

\begin{abstract}
We study conditions for the $L^p$ convergence of the AMV Laplacian $\Delta_r$ on completions of Riemannian manifolds $(M,g)$. The condition invloves the Minkowski dimension of the set $\overline{M} \setminus M$ and how much the curvature is blowing up around $S$. We build a concrete family of sets that satisfies these conditions, and provide an example of non convergence when these conditions are not met.
\end{abstract}

\tableofcontents

\section{Introduction}

\subsection{Graph Laplacians}

Graph Laplacians are a basic tool in geometric data analysis. Starting from a finite data set, one builds a weighted graph and studies the spectrum of a graph Laplacian. The resulting eigenvectors are used in spectral clustering, Laplacian eigenmaps, diffusion maps, and related dimensionality-reduction methods; see, for example, \citep{BelkinNiyogi2008Foundation,BelkinNiyogi2006Eigenmaps,CoifmanLafon2006DiffusionMaps,HeinAudibertVonLuxburg2007}. A central mathematical question is therefore to understand when these graph-based operators converge to differential operators on the underlying space.

Let \(X_1,\ldots,X_n\) be sample points in a metric space \((X,d)\), let \(k\colon [0,\infty)\to [0,\infty)\) be a kernel profile, and let \(r>0\) be a bandwidth. A standard choice of weights is
\[
w_{ij}^{(r)}
=
k\left(\frac{d(X_i,X_j)^2}{r^2}\right).
\]
For a function \(f\) on the sample, the corresponding \textit{unnormalized} and \textit{random-walk} graph Laplacians may be written as
\begin{align}
(\mathcal L^{\mathrm{un}}_{r,n}f)(X_i)
=&
\frac{1}{n r^{m+2}}
\sum_{j=1}^n
w_{ij}^{(r)}
\bigl(f(X_i)-f(X_j)\bigr),\\
(\mathcal L^{\mathrm{rw}}_{r,n}f)(X_i)
=&
\frac{1}{r^2}
\frac{
\sum_{j=1}^n
w_{ij}^{(r)}
\bigl(f(X_i)-f(X_j)\bigr)
}{
\sum_{j=1}^n
w_{ij}^{(r)}
},
\end{align}
whenever the denominator in the second expression is nonzero. The factor \(r^{-m}\) in the unnormalized operator reflects the volume scale in dimension \(m\). In the random-walk normalization this factor cancels between numerator and denominator.

The continuum counterpart of the random-walk graph Laplacian is obtained by replacing the empirical measure by a measure \(\nu\). On a metric measure space \((X,d,\nu)\), it is
\[
(L^{\mathrm{rw}}_{r,\nu}f)(x)
=
\frac{1}{r^2}
\frac{
\int_X
k\left(\frac{d(x,y)^2}{r^2}\right)
\bigl(f(x)-f(y)\bigr)\,d\nu(y)
}{
\int_X
k\left(\frac{d(x,y)^2}{r^2}\right)\,d\nu(y)
},
\]
whenever the denominator is positive. If \(k=\mathbf 1_{[0,1]}\), then the kernel selects precisely the ball \(B_r(x)\). Thus
\[
(L^{\mathrm{rw}}_{r,\nu}f)(x)
=
\frac{1}{r^2\nu(B_r(x))}
\int_{B_r(x)}
\bigl(f(x)-f(y)\bigr)\,d\nu(y).
\]
With the sign convention used in this paper, the asymptotic mean value Laplacian is
\[
\Delta^{\mathrm{AMV}}_{\nu,r}f(x)
=
\frac{1}{r^2\nu(B_r(x))}
\int_{B_r(x)}
\bigl(f(y)-f(x)\bigr)\,d\nu(y).
\]
Therefore, for the ball kernel,
\[
L^{\mathrm{rw}}_{r,\nu}f
=
-\Delta^{\mathrm{AMV}}_{\nu,r}f.
\]
This is the precise sense in which the AMV operator is tied to graph Laplacians: it is, up to sign, the population random-walk graph Laplacian associated with the indicator kernel. The present paper studies this population operator. It does not address the separate empirical concentration problem needed to pass from a finite random graph to the population operator.

For smooth manifolds, convergence of graph Laplacians has been studied from several points of view. Pointwise convergence results identify the local differential operator obtained in the small-bandwidth limit. In the usual smooth interior case, graph Laplacians approximate weighted Laplace--Beltrami operators, with the precise limit depending on the kernel, the normalization, and the sampling density; see \citep{BelkinNiyogi2008Foundation,HeinAudibertVonLuxburg2007,Singer2006GraphToManifold}. Spectral convergence is a stronger form of convergence: it concerns eigenvalues and eigenfunctions, and is directly connected to the use of graph Laplacians in clustering and dimensionality reduction. Results in this direction include \citep{BelkinNiyogi2006Eigenmaps,VonLuxburgBelkinBousquet2008,BuragoIvanovKurylev2014,GarciaTrillosGerlachHeinSlepcev2020,CalderGarciaTrillos2022, Dias2026DiscontinuousKernels}. 

Spectral convergence is a stronger form of convergence: it concerns eigenvalues
and eigenfunctions, and is directly connected to the use of graph Laplacians in
clustering and dimensionality reduction. Results in this direction include
\citep{BelkinNiyogi2006,LuxburgBelkinBousquet2008,BuragoIvanovKurylev2014,
GarciaTrillosEtAl2020,CalderGarciaTrillos2022}; for the intrinsic ball-kernel
setting of symmetrized AMV operators, see \citep{DiasTewodrose2026SpectralAMV}.

A complementary problem, lying between local pointwise consistency and full
spectral convergence in terms of global information, is convergence in \(L^p\).
Pointwise asymptotics identify the limiting differential operator at individual
points, whereas \(L^p\)-convergence gives a global function-space consistency
statement for the action of the operators on fixed test functions. This is not strictly weaker than spectral convergence, but can be a stepping stone to proving the spectral convergence; see, for instance, \citep{GarciaTrillosSlepcev2016,
GarciaTrillosEtAl2020,DiasTewodrose2026SpectralAMV}. This \(L^p\)-level
problem is the first motivation of the present paper.

Pointwise asymptotics identify the local limiting operator, but by themselves they do not control the global behavior of the operator. Spectral convergence usually requires additional compactness, variational, or operator-theoretic input. \(L^p\)-convergence provides a global function-space form of convergence and can be viewed as a step toward later spectral results. This is the first motivation of the present paper.

The second motivation is that data spaces are often not smooth manifolds. Singularities arise in many ways: data may be constrained, truncated, glued, partially observed, or identified under symmetries. Intersections of submanifolds, quotient-type singularities, collapsed directions, boundaries, corners, cusps, and other degeneracies all lead to spaces that are smooth only away from a singular set. For graph Laplacians and related averaging operators, the presence of such a singular set can change the asymptotic behavior. Near boundaries and singularities, the local averaging region need not look like a full Euclidean ball, and first-order terms or different scaling regimes may appear. This phenomenon has been studied for graph Laplacians on manifolds with boundary and related singular spaces in \citep{ZhouBelkin2011Boundary,BelkinQueWangZhou2012,PalTewodrose2025Kinks,Pal2025IsolatedSingularity}; see also \citep{PeoplesHarlim2026Boundary,BungertSlepcev2025Intersecting} for recent spectral-convergence results in boundary and singular settings.

\subsection{Asymptotic Mean Value (AMV) Laplacians}
The AMV Laplacian was introduced and studied on metric measure spaces in \citep{MinneTewodrose2020AMV}; related symmetrized AMV operators and their spectral properties were studied in \citep{MinneTewodrose2023SymAMV,DiasTewodrose2026SpectralAMV}. On a smooth \(m\)-dimensional Riemannian manifold, the AMV operator converges at interior points to \(\frac{1}{2(m+2)}\Delta_g f\). The construction is useful here because it is intrinsic to the metric-measure structure and uses only averages over metric balls. In graph-Laplian language, it corresponds exactly to the ball-kernel random-walk normalization described above, up to the sign convention.

\subsection{Our contribution}
The purpose of this paper is to prove \(L^p\)-convergence results for AMV operators on singular metric-measure spaces obtained from incomplete Riemannian manifolds. More precisely, let \((M,g)\) be a smooth \(m\)-dimensional Riemannian manifold and let \((\overline M,d)\) be its metric completion. We write \(S=\overline M\setminus M\) for the singular set and extend the Riemannian volume measure by giving zero mass to \(S\). This gives a metric measure space \((\overline M,d,\mu)\). Theorem~\ref{thm:convergenceOfAMV} gives conditions under which
\[
\Delta^{\mathrm{AMV}}_{\mu,r} f
\longrightarrow
\frac{1}{2(m+2)}\Delta_g f
\]
in \(L^p(\overline M,\mu)\), or weakly in a borderline case, as \(r\downarrow 0\). Equivalently, with the ball kernel, this gives the corresponding \(L^p\)-convergence statement for the population random-walk graph Laplacian \(L^{\mathrm{rw}}_{r,\mu}\), with the opposite sign.

The main theorem gives such convergence under three types of assumptions. First, the singular set is required to be small in the sense of an upper Minkowski-content condition. Second, the curvature is allowed to blow up near \(S\), but only at a controlled rate of the form
\[
K(x)\leq \frac{C}{d(x,S)^\alpha}.
\]
Third, the injectivity radius away from \(S\) is assumed to be controlled from below in terms of the distance to \(S\). Under these hypotheses, and assuming a local doubling condition for the completed metric-measure space, we prove strong \(L^p\)-convergence for Lipschitz functions in the range determined by the size of \(S\), the curvature blow-up exponent, and \(p\). We also prove corresponding \(L^q\)-convergence results for Hölder functions, with additional integrability conditions involving the gradient term.

The proof separates the space into two regions. Away from the singular set,
exponential coordinates are available at the scale \(r\), and the AMV operator is
compared with the usual second-order Taylor expansion. The curvature and
injectivity assumptions give uniform control of the exponential map and volume
density through Lemmas~\ref{lem:exponentialMapDerivative} and~\ref{lem:volumeFormExp},
while Proposition~\ref{prop:spaceIsDoubling} gives the corresponding local
Ahlfors/doubling control for metric balls. Near the singular set, no smooth
coordinate argument is used; instead, the contribution is controlled by Hölder or
Lipschitz regularity and by the size of tubular neighborhoods of \(S\). This is
where the Minkowski-content hypothesis enters.

We also give a more concrete geometric class for which the abstract assumptions can be checked, cf. Definition \ref{def:riemannianManifoldWithSingularitySet}. In this class, the completion is described by finitely many parametrizations that
are smooth away from a compact singular set and whose pullback metrics are
uniformly comparable with the Euclidean metric, as in
\eqref{eq:boundBellowMetric} and \eqref{eq:boundAboveMetric}. Under suitable size assumptions on the inverse images of the singular set, the completion has the expected topology, the resulting metric-measure space is locally doubling, and the required Sobolev and injectivity-radius estimates hold. This provides a family of examples to which the main convergence theorem applies.

Finally, we include an example of non-convergence (cf. Proposition \ref{prop:nonConvergenceResult}). In dimension two, we construct a metric on the punctured sphere with a point singularity and curvature blow-up of order \(d(\cdot,S)^{-2}\). For a Lipschitz function whose Laplace--Beltrami image is still integrable, the AMV operators fail to converge in the expected \(L^p\) range. This illustrates that control of the singular contribution and the Sobolev regularity assumptions cannot simply be ignored.

\subsubsection*{Acknowledgments} The
author is supported by the Research Foundation– Flanders (FWO) via the Odysseus
programme Geometric and analytic properties of metric measure spaces with spectral
curvature constraints, with applications to manifold learning (G0DBZ23N). The author also thanks Susovan Pal for various discussions about the project which helped shape its direction in particular in the construction of the counter-example and writing of the introduction.


\section{Preliminaries}

Let $(M,g)$ be a Riemannian manifold of dimension $m$. The Riemannian metric induces the distance
\begin{align}
\label{eq:distanceByCurves}
    \dist_g(x,y)
:=&
    \inf
    \Bigg\{
        \int_0^1 {g_{\gamma(t)}\left( \gamma' (t),  \gamma' (t)\right) }^\frac{1}{2} dt: \gamma(0) = x, \gamma(1) = y,\\
&\quad\quad\quad
    \gamma:[0,1] \rightarrow M, \text { and } \gamma \text{ is piecewise } C^1
    \Bigg\}.
\end{align}
The admissible class of curves in \eqref{eq:distanceByCurves} can be enlarged by absolutely continuous or Lipschitz curves, without changing the distance.

Given a point $x \in M$ and $r > 0$, define
\begin{align}
    B_r(0_x)
:=&
    \left\{
        v \in T_x M:
        g_x(v,v)^\frac{1}{2}<r
    \right\},
\end{align}
which is the ball in the tangent space, and
\begin{align}
    B_r^M(x)
:=&
    \left\{
        y \in M :
        \dist_g(y,x) < r
    \right\}
\end{align}
which is the ball in the manifold $M$.

Let $(\overline{M},\overline{\dist}_g)$ be the completion of the metric space $(M,\dist_g)$, and define the singular set by
\begin{equation}
    S
:=
    \overline{M}\setminus M.
\end{equation}
Given $x \in \overline{M}$ and $r >0$, define
\begin{equation}
    B_r(x) := \{y \in \overline{M}:\overline{\dist}_g(x,y)<r\},
\end{equation}
the ball in $\overline{M}$. For any point $x \in M$ such that
\begin{equation}
    \overline{\dist}_g(x,S)
>
    r,
\end{equation}
the exponential map $\exp_x: B_r(0_x) \rightarrow M$ is well defined. This allows us to define the injectivity radius at a point by
\begin{equation}
    \inj(x) 
:=
    \sup\left\{
        r > 0: \exp_x : B_r(0_x) \rightarrow B_r(x) \text{ is a diffeomorphism}
    \right\}.
\end{equation}
The metric also induces a measure $\vol_{g}$ on $M$, which extends to a measure on $\overline{M}$ by
\begin{equation}
    \mu(E) := \vol_g(E \setminus S).
\end{equation}
This defines a metric measure space 
\begin{align}
    \overline{\mathfrak{M}}
:=
    \left(
        \overline{M}, \overline{\dist}_g, \mu
    \right).
\end{align}

Let $\nabla$ be the Levi-Civita connection associated with the Riemannian metric $g$. Given a point $x \in M$, let $R_x\left(\cdot, \cdot\right)\cdot : T_xM^3 \rightarrow T_xM$ be the Riemannian curvature tensor. 

\begin{definition}
\label{def:suppremumOfCurvatures}
    We define $K : M \rightarrow \mathbb{R}$ as the pointwise norm of the Riemann tensor; that is,
    \begin{align}
        \label{eq:supremmumOfTheSectionalCurvatures}
            K(x)
        :=&
            \sup_{\{v,w,z \in T_xM : g_x(v,v) = g_x(w,w) = g_x(z,z) = 1\}}
            g_x(R_x(v, w)z, R_x(v, w)z)^\frac{1}{2}\\
        =&
            \sup_{\left\{u,v,w,z \in T_xM : g_x(u,u) = g_x(v,v) = g_x(w,w) = g_x(z,z) = 1\right\}}
            g_x(R(v,w)z,u).
    \end{align}
\end{definition}
We now introduce the notion of Minkowski content.
\begin{definition}\label{def:tubular neighborhood}
    Given a metric space $(X,d)$, a set $A \subset X$, and $r > 0$, define the \textit{tubular neighborhood} of radius $r$ around $A$ by
    \begin{equation}
        T_{r,d}(A)
    :=
    \{
        y \in X:
        d(y, A) < r
    \}.
    \end{equation}
\end{definition}

\begin{definition}
\label{def:minkowskyContent}
    Consider $m > 0$ and a metric measure space $\mathfrak{X} = (X, d, \nu)$. Given $A  \subset X$ and $0 \leq k \leq m$, we define the upper and lower $(m,k)$-dimensional Minkowski contents, denoted respectively by $\mathcal{M}^{(m,k,+)}_{\mathfrak{X}}(A)$ and $\mathcal{M}^{(m,k,-)}_{\mathfrak{X}}(A)$, by
    \begin{align}
        \mathcal{M}^{(m,k,+)}_{\mathfrak{X}}(A)
    &:=
        \limsup_{r \rightarrow 0}
        \frac{\nu(T_{r,d}(A))}{r^{m-k}},\\
        \mathcal{M}^{(m,k,-)}_{\mathfrak{X}}(A)
    &:=
        \liminf_{r \rightarrow 0}
        \frac{\nu(T_{r,d}(A))}{r^{m-k}}.
    \end{align}
    The Minkowski dimensions are defined, respectively, by
    \begin{align}
        \text{dim}^{(m,+)}_{\mathcal{M}, \mathfrak{X}}(A)
    &=
        \inf
        \left\{k\in [0,m] : \mathcal{M}^{(m,k,+)}_{\mathfrak{X}}(A) = 0\right\},\\
        \text{dim}^{(m,-)}_{\mathcal{M}, \mathfrak{X}}(A)
    &=
        \inf
        \left\{k\in [0,m] : \mathcal{M}^{(m,k,-)}_{\mathfrak{X}}(A) = 0\right\}.
    \end{align}
\end{definition}
In this definition, $m$ is understood to be the intrinsic dimension of the metric measure space $\mathfrak{X}$. One can find the definition and results for the Euclidean metric measure space $\mathfrak{E}:=(\mathbb{R}^m, \dist_\euc, \mathcal{L}^m)$ in \cite{AFP}.
We also recall the Sobolev spaces on Riemannian manifolds. We only need Sobolev spaces up to order two. For a reference on Sobolev spaces on manifolds, see \cite{sobolevSpacesManifoldHebey}.
\begin{definition}
\label{def:sobolevSpace}
    Given the Riemannian manifold $(M, g)$, let $\nabla$ be the Levi-Civita connection on $M$.
    Given $f \in C^\infty(M)$ and $x \in M$, define
    \begin{equation}
        \|\nabla f(x)\|_g
    =
        \sup_{\{v \in T_xM: g_x(v,v)=1\}}
        |\nabla_v f(x)|
    =
        g\left(\nabla f\left(x\right),\nabla f\left(x\right)\right)^\frac{1}{2}
    \end{equation}
    and
    \begin{equation}
        \|H f(x)\|_g
    =
        \sup_{\{v,w \in T_xM : g_x(v,v) = g_x(w,w) = 1\}}
        \left|\left(\nabla d f\right)_x(v,w)\right|.
    \end{equation}
    Consider the norms given by
    \begin{align}
        \|f\|_{W^{1,p}(M, g)}
    :=&
        \|f\|_{L^p(M, \vol_g)}
    +
        \left\|\left\|\nabla f\right\|_g\right\|_{L^p(M, \vol_g)},\\
        \|f\|_{W^{2,p}(M,g)}
    :=&
        \|f\|_{L^p(M, \vol_g)}
    +
        \left\|\left\|\nabla f\right\|_g\right\|_{L^p(M, \vol_g)}
    +   
        \left\|\left\|H f\right\|_g\right\|_{L^p(M, \vol_g)}.
    \end{align}
    We define the Sobolev spaces $W^{1,p}(M,g)$ and $W^{2,p}(M,g)$ as the completions of
    \begin{equation}
        \left\{
            f \in C^\infty(M):
            \left\|f\right\|_{W^{1,p}(M,g)} < \infty
        \right\},
        \text{ or }
        \left\{
            f \in C^\infty(M):
            \left\|f\right\|_{W^{2,p}(M,g)} < \infty
        \right\},
    \end{equation}      
    with respect to the corresponding norm.
\end{definition}
Sobolev inequalities form an important class of estimates in Sobolev theory. These are not always satisfied on Riemannian manifolds, and the subject has been well developed on complete manifolds (see, for example, \cite{sobolevSpacesManifoldHebey}). However $(M,g)$ is not necessarily complete, and the theory in this case is not as well studied.
\begin{definition}
    We say that the space $(M,g)$ satisfies property \eqref{eq:SobolevIneq} with exponent $1\leq p < m$ if, for $p^*>0$ satisfying $\frac{1}{p^*} = \frac{1}{p} - \frac{1}{m}$, there exists a constant $C_{\mathcal{S},p} > 0$ such that for all $f \in W^{1,p}(M,g)$ we have
    \begin{equation}
    \label{eq:SobolevIneq}
    \tag{S}
        \|
            f
        \|_{L^{p^*}(M,\vol_g)}
    \leq
        C_{\mathcal{S},p}
        \left(
        \left\|
        \left\|\nabla f\right\|_g
        \right\|_{L^p(M,\vol_g)}
        +
        \left\|
            f
        \right\|_{L^p(M,\vol_g)}
        \right).
    \end{equation}
\end{definition}

We now define the asymptotic mean value (AMV) Laplacian operators.
\begin{definition}
    \label{def:amvLaplacian}
    Given a metric measure space $\mathfrak{X} = (X,d,\nu)$ such that every ball has positive finite measure, we define the AMV Laplacian operators $\Delta_r$ for $f \in L^1_{\text{loc}}\left(X,\nu\right)$ by
    \begin{equation}
        \Delta_r f(x)
    :=
        \frac{1}{r^2\nu(B_r^d(x))}
        \int_{B_r^d(x)}
            \left(
                f(y)
            -
                f(x)
            \right)
            d\nu(y),
    \end{equation}
    where $B_r^d(x):=\{y \in X:d(x,y)<r\}$.
\end{definition}
A metric-measure property that will be useful for the AMV Laplacian is local doubling.

\begin{definition}
    Given a metric measure space $\mathcal{X}=(X,d,\nu)$, we say that the space is locally doubling if there exist $C,r_0>0$ such that for all $r \leq r_0$ and $x \in X$ we have
    \begin{equation}
        \nu(B_{2r}(x))
    \leq
        C\nu(B_r(x)).
    \end{equation}
\end{definition}

The objective of this paper is to identify conditions under which, for a function $f \in L^p(\overline{M},\mu)$, the operators $\Delta_r f$ converge in $L^p(\overline{M},\mu)$ as $r \rightarrow 0$ on the metric measure space $\overline{\mathfrak{M}}=(\overline{M},\overline{\dist}_g,\mu)$.

\begin{theorem}
\label{thm:convergenceOfAMV}
    Let $(M,g)$ be a Riemannian manifold such that the completion $(\overline{M},\overline{\dist}_g)$ of $(M, \dist_g)$ is compact. Consider the set $S = \overline{M} \setminus M$ and suppose that, for some $0 \leq k < m-1$, we have $\mathcal{M}^{(m,k,+)}_{\overline{\mathfrak{M}}}(S) < \infty$. Let $m = \dim(M)$, $C_m = \frac{1}{2(m+2)}$, $1 < p < \infty $, $0 \leq \alpha \leq m$, and $C_K > 0$. Suppose that
    \begin{equation}
    \label{eq:curvatureCondition}
    \tag{K}
        K(x) 
    \leq 
        \frac{C_K}{\overline{\dist}_g(x,S)^\alpha}.
    \end{equation}
    Furthermore, suppose there exist $r_0>0$ and $C_{\alpha, \inj} > 2$ such that
    \begin{equation}
    \label{eq:InjectivityCondition}
    \tag{I}
        \inj(x)
    \geq
        \min\left\{
            r_0,
            \left(
                \frac{\overline{\dist}_g(x,S)}{C_{\alpha, \inj}}
            \right)^{\max\{1,\frac{\alpha}{2}\}}
        \right\},
    \end{equation}
    and that $\left(\overline{M},\overline{\dist}_g, \mu\right)$ is a locally doubling metric measure space.
    
    Let $f:\overline{M} \rightarrow \mathbb{R}$ satisfy $f \in C^\infty(M)$ and $f \in W^{2,p}(M,g)$.
    \begin{enumerate}
        \item If $f$ is Lipschitz on $\overline{M}$, then 
        \begin{equation}
        \label{eq:conditionForConvergence}
            m-k-
            p
            \max
            \left\{1,\frac{\alpha}{2}
            \right\} > 0
            \quad
            \Rightarrow
            \quad
            \Delta_r f
            \rightarrow_{L^p(\overline{M},\mu)} C_m \Delta_g f.
        \end{equation}
        as $r \rightarrow 0$.
        If $\alpha < 2$ and $p>1$ then
        \begin{equation}
            m-k-p = 0 
        \quad 
        \Rightarrow
        \quad
        \Delta_r f
        \rightharpoonup_{L^p(\overline{M},\mu)}
        C_m \Delta_g f
        \end{equation}
        as $r \rightarrow 0$.
    \item If $f$ is H\"older continuous with exponent $0<\beta<1$ and $1<q<p<m$ satisfies
    \begin{equation}
        \label{eq:holderResultconvergencestrong}
        m-k-q\max\left\{2-\beta, \frac{\alpha(2-\beta)}{2}\right\}>0,
    \quad
        m-k-\frac{\alpha pq}{2(p-q)}> 0
    \end{equation}
    then
    \begin{equation}
    \Delta_r f \rightarrow_{L^q(\overline{M}, \mu)} C_m \Delta_g f
    \end{equation}
    as $r \rightarrow 0$. If $m-k-\frac{\alpha pq}{2(p-q)} > 0$ and $p>1$ then
    \begin{equation}
    \label{eq:holderconvergenceweak}
        m-k-q\max\left\{2-\beta, \frac{\alpha(2-\beta)}{2}\right\}
        =
        0
        \Rightarrow
        \Delta_r f
        \rightharpoonup_{L^q(\overline{M},\mu)}
        C_m \Delta_g f.
    \end{equation}
    as $r \rightarrow 0$. Furthermore, if $(M,g)$ satisfies the Sobolev inequality \eqref{eq:SobolevIneq} with exponent $p$, then in the condition involving $\frac{\alpha pq}{2(p-q)}$, one may replace $p$ by $p^* = \frac{pm}{m-p}$, provided $1<q < p$.
    \end{enumerate}
\end{theorem}

Since conditions \eqref{eq:InjectivityCondition} and \eqref{eq:SobolevIneq} may not be easy to verify, we introduce a family of manifolds that satisfies them. This relatively general condition allows us to identify more easily examples satisfying the hypotheses of Theorem~\ref{thm:convergenceOfAMV}.
\begin{definition}
\label{def:riemannianManifoldWithSingularitySet}
    We say that the tuple $(N, S, g, \{\psi_i\}_{i=1}^l)$ is a compact Riemannian manifold of dimension $m$ with $C^\infty$ Riemannian metric $g$ on $N \setminus S$ and a compact singularity set $S \subset N$ with a finite number of parametrizations $\{\psi_i:\mathbb{B}_2^m(0) \rightarrow N\}_{i=1}^l$ if
    \begin{enumerate}
        \item Each map $\psi_i : \mathbb{B}^m_2(0) \rightarrow N$ is a homeomorphism onto the image. Furthermore, the parametrizations $\psi_i:\mathbb{B}_2^m(0) \setminus \psi_i^{-1}(S) \rightarrow N \setminus S$ form a $C^\infty$ atlas of $N \setminus S$, that is they are homeomorphisms onto the image and for $i,j \in \{1,...,l\}$ the map $\psi_j^{-1}\circ \psi_i:\psi_i^{-1}\left(\psi_j\left(\mathbb{B}_2^m(0)\right)\setminus S\right) \rightarrow \mathbb{R}^m$ is $C^\infty$, and $\psi_i\left(\mathbb{B}^m_2(0)\right)$ is open in $N$. Furthermore, $N$ is a Hausdorff topological space.
        \item $N$ is covered by
            \begin{equation}
                \label{eq:NcoveredFiniteParametrizations}
                N = \cup_{i=1}^l\psi_i\left(\mathbb{B}_1^m\left(0\right)\right).
            \end{equation}
        \item $(N\setminus S, g)$ is a Riemannian manifold  with metric $g$ such that there exist constants $c_{LB},C_{UB}>0$ such that for all $i \in \{1,...,l\}$ and $x \in \psi_i^{-1}(N\setminus S)$
        \begin{align}
            \tag{LB}
            \label{eq:boundBellowMetric}
                (\psi_i^*g)_x(v, v)
            &\geq
                c_{LB}\|v\|^2,\\
            \tag{UB}
            \label{eq:boundAboveMetric}
                (\psi_i^*g)_x(v,v)
            &\leq
                C_{UB}\|v\|^2.
            \end{align}
    \end{enumerate}
\end{definition}
Let $\mathfrak{E} = (\mathbb{R}^m, \dist_\euc, \mathcal{L}^m)$ be the Euclidean metric measure space.
\begin{theorem}
\label{thm:riemannianManifoldWithSingularityThm}
        Consider a tuple $(N,S,g, \{\psi_i\}_{i=1}^l)$ that \textcolor[HTML]{5C4033}{satisfies Definition~\ref{def:riemannianManifoldWithSingularitySet}} of dimension $m$. Consider $0\leq k < m-1$ and suppose that for all $i \in \{1,...,l\}$ we have $\mathcal{M}^{(m,k,+)}_{\mathfrak{E}}(\psi_i^{-1}(S)) < \infty$. Then the completion of $(N\setminus S, \dist_g)$ is homeomorphic to $N$ with the topology given by the parametrizations $\left\{\psi_i\right\}_{i=1}^l$, and $(N\setminus S, g)$ satisfies \eqref{eq:SobolevIneq} for exponents $p$ satisfying $m-k-p>0$. Moreover, given the measure $\nu(E) := \vol_g(E\setminus S)$ and the metric measure space $\mathfrak{N}:=(N,\overline{\dist}_g, \nu)$ we have that $\mathcal{M}^{(m,k,+)}_{\mathfrak{N}}(S) < \infty$ and that the space is a locally doubling metric measure space.

Furthermore, given $1 \leq \alpha \leq m$ and $C_K>0$ if the function $K:N \setminus S \rightarrow \mathbb{R}$ from Definition~\ref{def:suppremumOfCurvatures} satisfies equation \eqref{eq:curvatureCondition} then the space $(N,\overline{\dist}_g)$ satisfies condition \eqref{eq:InjectivityCondition} for some $C_{\dist_g,\alpha} > 2$.
\end{theorem}

Beyond the parametrized class in Definition~\ref{def:riemannianManifoldWithSingularitySet}, Proposition~\ref{prop:spaceIsDoubling} shows that if a space satisfies \eqref{eq:curvatureCondition} with $\alpha=2$ and \eqref{eq:InjectivityCondition}, and is locally Ahlfors regular on $S$, then $(\overline{M},\overline{\dist}_g, \mu)$ is locally Ahlfors and hence locally doubling.

\section{Geometric preliminaries}

In what follows, whenever $\alpha=0$, $\min\left\{1,\frac{2}{\alpha}\right\} = 1$.

\begin{lemma}
\label{lem:hessianAndInjecRadius}
    Suppose that $(M,g)$ is a Riemannian manifold that satisfies \eqref{eq:InjectivityCondition} with constants $C_{\alpha,\inj} > 2$ and $r_0>0$. Let $f \in C^\infty(M)$. Given $0 < r \leq r_0$, if $\overline{\dist}_g(x,S) \geq C_{\alpha, \inj} r^{\min\left\{1,\frac{2}{\alpha}\right\}}$ then $\exp_x : B_r(0_x) \rightarrow B_r(x)$ is a diffeomorphism and, for every $w \in \partial B_1(0_x)$ and $0<t<r$, we have
    \begin{equation}
    \left|
        \frac{d^2}{dt^2}
        f\circ \exp_x(tw)
    \right|
    \leq
    \|
        Hf(\exp_x(tw))
    \|_g.
    \end{equation}
\end{lemma}
    \begin{proof}
        For $0<r\leq r_0$, condition \eqref{eq:InjectivityCondition}
        \begin{equation}
            \inj(x)
        \geq
            \min\left\{r_0, \left(\frac{\overline{\dist}_g(x,S)}{C_{\alpha,\inj}}
            \right)^{\max\left\{1,\frac{\alpha}{2}\right\}}
            \right\}
        \geq
            \min\left\{r_0, r
            \right\}
            =
            r,
        \end{equation}
        showing that $\exp_x:B_r(0_x) \rightarrow B_r(x)$ is a diffeomorphism.

        Now consider the radial vector field in $T_xM \setminus \{0\}$ given by
        \begin{equation}
            \Tilde{X}_w:= \frac{w}{g_x(w,w)^\frac{1}{2}},
        \end{equation}
        and we may therefore define the vector field in $B_r(x) \setminus \{x\}$ given by
        \begin{equation}
            X_z
        :=
            (d\exp_x)_{\exp_x^{-1}(z)}
            \Tilde{X}_{\exp_x^{-1}(z)}.
        \end{equation}
        This vector field satisfies
        \begin{equation}
            \nabla_{X}X = 0.
        \end{equation}
        Thus, by the definition of the covariant derivative of a one-form, for $w \in \partial B_1(0_x)$ we obtain
        \begin{align}
        \left|
        \frac{d^2}{dt^2}
        f\circ \exp_x(tw)
        \right|
    =&
        \left|
        \left.\frac{d}{ds}\right|_{s=t}
        \left(
            \nabla_{X}
            f
        \right)({\exp_x(sw)})
        \right|\\
    =&
        \left|
        X
        \cdot
        \left(
            \nabla_{X}
            f
        \right)({\exp_x(sw)})
        \right|
        \\
    =&
        \left|
            X
            \cdot
            (
            X\cdot
            f
            )_{\exp_x(tw)}
            +
            \nabla_{\nabla_{X}X}
            f(\exp_x(tw))
        \right|\\
    \leq&
        \|
            Hf(\exp_x(tw))
        \|_g
        g_{\exp_{x}(tw)}(X,X)
    =
        \|
            Hf(\exp_x(tw))
        \|_g.
    \end{align}
\end{proof}

\begin{theorem}[Rauch comparison theorem \cite{Chavel_2006}]
\label{thm:RauchComparisonTheorem}
    Given a real number $k \in \mathbb{R}$, define
    \begin{equation}
        s_k(r)
    :=
        \begin{cases}
          \frac{\sin(\sqrt{k}r)}{\sqrt{k}}  & k>0 \\
          r & k = 0\\
          \frac{\sinh(\sqrt{-k}r)}{\sqrt{-k}} & k<0
        \end{cases}
    \end{equation}
    Let $T>0$ and consider an arc-length-parametrized geodesic $\gamma:[0,T]\rightarrow M$ such that $g_{\gamma(0)}\left(
        \gamma'(0), 
        \gamma'(0)
    \right) = 1$. Suppose that there exists $k > 0$ such that  
    \begin{equation}
        \sup_{u \in T_{\gamma(t)}M, g_{\gamma(t)}(u,u) = 1, g_{\gamma(t)}(u,\gamma'(t))  = 0}|k_{\gamma(t)}(u,\gamma'(t))|\leq k,
    \end{equation}
    where $k_x(u,w) = \frac{g_x\left(R_x\left(u,w\right)w,u\right)}{\left\|u\right\|_g^2\left\|w\right\|_g^2 - g_x\left(u,w\right)^2}$ is the sectional curvature of the plane spanned by $u$ and $w$. If $J$ is a Jacobi field along $\gamma$ such that 
    \begin{align}
    g_{\gamma(0)}(J_{\gamma(0)}, J_{\gamma(0)}) 
    &= 0\\
    g_{\gamma(0)}\left(\nabla_{\gamma'(0)}J_{\gamma(0)}, \nabla_{\gamma'(0)}J_{\gamma(0)}\right) &= 1
    \end{align}
    then we have
    \begin{equation}
    s_{k}(t)
    \leq
        g_{\gamma(t)}\left(
            J_{\gamma(t)},
            J_{\gamma(t)}
        \right)^\frac{1}{2}
    \leq
    s_{-k}(t).
    \end{equation}
    for all $t < \min \{T, \frac{\pi}{\sqrt{k}}\}$.
\end{theorem}
Given $w \in \partial B_1(0_x)$, let $\{v_{1,w},...,v_{m-1,w}\} \subset T_xM$ be an orthonormal set of vectors orthogonal to $w$. Consider the Jacobi fields along the geodesic $\gamma(t) = \exp_x(tw)$ given by
\begin{equation}
\label{eq:JacobiFieldsAlongGeo}
    J_{i,w} := (d\exp_x)_{tw}tv_{i,w}.
\end{equation}
Notice that if $\overline{\dist}_g(x,S) > r$ then the exponential map $\exp_x:B_r(0_x) \rightarrow M$ is well defined.

\begin{lemma}
    \label{lem:exponentialMapDerivative}
    Suppose the manifold $(M,g)$ satisfies \eqref{eq:curvatureCondition}. There exist $C_{\alpha,\dist_g}>2$ and $0<r_0\leq1$ with the following property. For all $0<r\leq r_0$ and $x\in M$ satisfying $\overline{\dist}_g(x,S)\geq C_{\alpha,\dist_g}r^{\min\left\{1,\frac{2}{\alpha}\right\}}$, and for each $w\in\partial B_1(0_x)$, let $\gamma(t)=\exp_x(tw)$. Then, for all $u \in \partial B_1(0_x)$ and $0 \leq t \leq r$
    \begin{align}
    \left(
    \frac{\sin\left(
        \sqrt{C_{K,r}(x)}t
        \right)}{\sqrt{C_{K,r}(x)}t}
    \right)^2
    \leq
        g_{\gamma(t)}
        \left(
            (d\exp_x)_{tw}
            u,
            (d\exp_x)_{tw}
            u
        \right)
    \leq
    \left(
        \frac{\sinh\left(\sqrt{C_{K,r}(x)}t\right)}{ \sqrt{C_{K,r}(x)} t}
    \right)^2,
    \end{align}
    where $C_{K,r}(x) := \max\left\{\sup_{y \in B_r(x)} K(y), 1
    \right\}$. Also in particular for all $u \in \partial B_1(0_x)$ and $0 \leq t \leq r$ 
    \begin{align}
    \label{eq:ineqsForExponentialJacobian}
        \frac{1}{2}
    \leq
        g_{\gamma(t)}
        \left(
            (d\exp_x)_{tw}
            u,
            (d\exp_x)_{tw}
            u
        \right)
    \leq
        2.
    \end{align}
\end{lemma}
\begin{proof}
    Because $\overline{\dist}_g(x,S) \geq C_{\alpha,\dist_g} r^{\min\left\{1,\frac{2}{\alpha}\right\}}$, we have $\sup_{y \in B_r(x)} K(y) \leq \frac{C_K}{r^2(C_{\alpha, \dist_g} -1)^\alpha}$. By the triangle inequality and the assumption $r\leq1$,
    \begin{align}
        C_{K,r}(x)
    \leq
    \max\left\{
        \sup_{y \in B_r(x)}
        \frac{C_K}{\overline{\dist}_g(y,S)^\alpha},1
    \right\}
    &\leq
    \max\left\{
    \frac{C_K}{\left( C_{\alpha,\dist_g}r^{\min\left\{1,\frac{2}{\alpha}\right\}} - r\right)^\alpha},1
    \right\}\\
    &\leq
    \max\left\{
        \frac{C_K}{\left( C_{\alpha,\dist_g}r^{{\min\left\{1,\frac{2}{\alpha}\right\}}} - r^{{\min\left\{1,\frac{2}{\alpha}\right\}}}\right)^\alpha},1
        \right\}\\
    &\leq
        \max\left\{
        \frac{C_K}{r^{\min\left\{\alpha,2\right\}}(C_{\alpha,\dist_g} - 1)^\alpha},1
        \right\}\\
    &\leq
        \max\left\{
        \frac{C_K}{r^2(C_{\alpha,\dist_g} - 1)^\alpha},1
        \right\}.
    \end{align}
    Consider $w \in \partial B_1(0_x) \subset T_xM$ and the Jacobi fields $J_{i,w}$ defined by \eqref{eq:JacobiFieldsAlongGeo} along the geodesic $\gamma(t) = \exp_x(tw)$. Consider now $v = \sum_{i=1}^{m-1}b_i v_{i,w}$ with $b_i \in \mathbb{R}$ and by orthogonality of these vectors $g_x(v,v) = \sum_{i=1}^{m-1}b_i^2$. Notice that $J_{v,w} = \sum_{i=1}^{m-1}b_i J_{i,w}$ is a solution of the Jacobi equation. Furthermore, since $\left\|\gamma'(t)\right\|_g = 1$, for all $0< t \leq r$, $u \in T_{\gamma(t)}M$ such that $\|u\|_g = 1$ and $g_{\gamma(t)}\left(u,\gamma'(t)\right) = 0$, we have $\left|
    k_{\gamma(t)}\left(u,\gamma'(t)\right)
    \right|
    \leq
    \sup_{y \in B_r(x)} K(y)
    \leq C_{K,r}(x)$. Thus we can apply Theorem \ref{thm:RauchComparisonTheorem} to conclude that for all $t \leq \min\left\{r, \frac{\pi}{\sqrt{C_{K,r}(x)}}\right\}$ 
    \begin{align}
        \left(
            \sum_{i=1}^{m-1}b_i^2
        \right)
        \left(
        \frac{\sin\left(
            \sqrt{C_{K,r}(x)}t
            \right)}{\sqrt{C_{K,r}(x)}}
        \right)^2
        \leq
            g_{\gamma(t)}\left(J_{v,w},J_{v,w}\right)
        \leq
        \left(
            \sum_{i=1}^{m-1}b_i^2
        \right)
        \left(
            \frac{\sinh\left(\sqrt{C_{K,r}(x)}t\right)}{ \sqrt{C_{K,r}(x)}}
        \right)^2.
    \end{align}
    Since $\frac{\pi}{\sqrt{C_{K,r}(x)}} \geq
    \min\left\{
    \frac{\pi (C_{\alpha,\dist_g}-1)^{\frac{\alpha}{2}}r}{\sqrt{C_K}}, \pi
    \right\}$, if $\alpha > 0$ we can take $C_{\alpha,\dist_g} > 2$ large enough so that $\frac{\pi (C_{\alpha,\dist_g}-1)^\frac{\alpha}{2}r}{\sqrt{C_K}}>r$ and so the inequalities are satisfied for all $0 \leq t \leq r$. If $\alpha = 0$, we can find $r_0 > 0$ small enough so that $\frac{\pi}{\sqrt{C_{K,r}(x)}} \geq \min\left\{\frac{\pi}{\sqrt{C_K}}, \pi\right\} > r_0$, and so the inequality will also be satisfied for $0\leq t \leq r \leq r_0$.
    
    Also notice that $J_{v,w}$ is a Jacobi field along $\gamma(t)$, and so we have
    \begin{equation}
        g_{\gamma(t)}\left(
            (d\exp_x)_{tw}w, J_{v,w}
        \right)
    =
        g_{\gamma(t)}\left(
            \gamma'(t), J_{v,w}
        \right)
    =
        0,
    \end{equation}  
    (This can be seen by noticing $h(t) = g_{\gamma(t)}\left(
            (d\exp_x)_{tw}w, J_{v,w}
        \right)$ satisfies a second order linear ODE with initial conditions equal to zero). Given a vector $u \in \partial B_1(0_x)$, by the orthonormality of the set $\{w,v_{1,w},...,v_{m-1,w}\}$ we can write it as $u = b_0 w + \sum_{i=1}^{m-1}b_i v_{i,w}$ where $\sum_{i=0}^{m-1}b_i^2 = 1$, thus for all $u \in \partial B_1(0_x)$ and $0 \leq t\leq r$ we have
    \begin{align}
        &g_{\gamma(t)}
        \left(
            (d\exp_x)_{tw}
            u,
            (d\exp_x)_{tw}
            u
        \right)\\
    =&
        g_{\gamma(t)}
        \left(
            b_0(d\exp_x)_{tw}w
            +
            \sum_{i=1}^{m-1}
                b_i
                (d\exp_x)_{tw}v_{i,w},
            b_0(d\exp_x)_{tw}w
            +
            \sum_{i=1}^{m-1}
                b_i
                (d\exp_x)_{tw}v_{i,w}
        \right)\\
    =&
        b_0^2
        +
        g_{\gamma(t)}
        \left(
            \sum_{i=1}^{m-1}
                b_i
                J_{i,w}/t,
            \sum_{i=1}^{m-1}
                b_i
                J_{i,w}/t
        \right)\\
    \geq&
        b_0^2
        +
        \left(\sum_{i=1}^{m-1}
        b_i^2
        \right)
        \left(
        \frac{\sin\left(
        \sqrt{C_{K,r}(x)}t
        \right)}{\sqrt{C_{K,r}(x)}t}
        \right)^2\\
    \geq&
        \left(\sum_{i=0}^{m-1}
        b_i^2
        \right)
        \left(
        \frac{\sin\left(
        \sqrt{C_{K,r}(x)}t
        \right)}{\sqrt{C_{K,r}(x)}t}
        \right)^2
    =
        \left(
        \frac{\sin\left(
        \sqrt{C_{K,r}(x)}t
        \right)}{\sqrt{C_{K,r}(x)}t}
        \right)^2,
        \label{eq:lowerBoundOnExp}
    \end{align}
    where we extend $\frac{\sin\left(
        \sqrt{C_{K,r}(x)}t
        \right)}{\sqrt{C_{K,r}(x)}t}$ to $1$ at $t = 0$. Similarly we can prove that
    \begin{equation}
        g_{\gamma(t)}
        \left(
            (d\exp_x)_{tw}
            u,
            (d\exp_x)_{tw}
            u
        \right)
    \leq
        \left(
        \frac{\sinh\left(
            \sqrt{C_{K,r}(x)}t
            \right)}{\sqrt{C_{K,r}(x)}t}
        \right)^2,
    \end{equation}
    where we extend $\frac{\sinh\left(
            \sqrt{C_{K,r}(x)}t
            \right)}{\sqrt{C_{K,r}(x)}t}$ to $1$ at $t = 0$. Notice now that there exists $\delta > 0$ such that for all $0 <y \leq \delta $ we have $\left(\frac{\sin(y)}{y}\right)^2 \geq \frac{1}{2}$ and $\left(\frac{\sinh(y)}{y}\right)^2 \leq 2$. Since $0 < t < r < r_0 < 1$, we have $\sqrt{C_{K,r}(x)}t \leq \max\left\{\frac{\sqrt{C_K}}{r\left( C_{\alpha,\dist_g} - 1\right)^{\alpha/2}} t,t
    \right\}
    \leq 
    \max\left\{\frac{\sqrt{C_K}}{\left( C_{\alpha,\dist_g} - 1\right)^{\alpha/2}},r_0\right\}$, if $\alpha > 0$ we can choose $C_{\alpha, \dist_g} > 2$ large enough and $r_0$ small enough so that $\sqrt{C_{K,r}(x)}t < \delta$, thus we conclude from the previous inequalities
    \begin{equation}
        \frac{1}{2}
        \leq
        g_{\gamma(t)}
        \left(
            (d\exp_x)_{tw}
            u,
            (d\exp_x)_{tw}
            u
        \right)
        \leq
        2.
    \end{equation}
    Similarly if $\alpha = 0$ we can choose $r_0>0$ small enough so that $\sqrt{C_{K,r}(x)}t\leq \max\left\{\sqrt{C_K}t ,t\right\}\leq \max\left\{\sqrt{C_K} r_0, r_0 \right\}< \delta$, the same inequality can be obtained.
\end{proof}

In exponential coordinates associated with $\exp_x$, define the volume density by
\begin{equation}
\label{eq:defintionOfVolumeFormExpCoord}
    \left(
    \left(\exp_x\right)^*
    \left(
        d\vol_g
    \right)
    \right)_{y}
=
    J_x(y)
    dy_1\wedge\cdots \wedge dy_m.
\end{equation}

\begin{lemma}
\label{lem:volumeFormExp}
    Suppose the manifold $(M,g)$ satisfies \eqref{eq:curvatureCondition} and the closure $\left(\overline{M},\overline{\dist}_g\right)$ is compact. 
    Let $x \in M$, and let $J_x(y)$ be defined by \eqref{eq:defintionOfVolumeFormExpCoord}. Then there exist $C_{\alpha, \dist_g} > 2$, $0< r_0 \leq 1$ and $C_J > 0$ such that for all $0< r \leq r_0$ and $x \in M$ with $\overline{\dist}_g(x,S) \geq C_{\alpha, \dist_g}r^{\min\left\{1,\frac{2}{\alpha}\right\}}$ and $y \in B_r(0_x)$ we have
    \begin{equation}
    \label{eq:ineqForVolumeForm}
        \left|
            1
        -
            J_x(y)
        \right|
        \leq
        \frac{C_J}
        {\overline{\dist}_g(x,S)^\alpha}g_x(y,y).
    \end{equation}
    and
    \begin{equation}
    \label{eq:ineqForVolumeForm2}
        \frac{1}{2}
        \leq
        J_x(y) \leq \frac{3}{2}.
    \end{equation}
    Furthermore, if condition \eqref{eq:InjectivityCondition} is satisfied with constants $r_0^{\text{inj}}>0$ and $C_{\alpha,\inj}>2$, then, after possibly decreasing $r_0$, for $r \leq r_0$ and $\overline{\dist}_g(x,S) \geq C_{\alpha, \dist_g}r^{\min\left\{1,\frac{2}{\alpha}\right\}}$ we have
    \begin{equation}
    \label{eq:ineqVolOfBalls}
        \mathcal{L}^m\left(\mathbb{B}^m_1(0)\right)\frac{r^m}{2}
    \leq
        \mu(B_r(x))
    \leq
        \mathcal{L}^m\left(\mathbb{B}^m_1(0)\right)\frac{3 r^m}{2}.
    \end{equation}
\end{lemma}
\begin{proof}
    Given $y \in B_r(0_x)$, let $y = tw$ where $w \in \partial B_1(0_x)$. Then $\{w,v_{1,w},...,v_{m-1,w}\}$ forms an orthonormal basis of $T_xM$. With respect to this basis, consider the matrix
    \begin{align}
        \left[
            A(y)
        \right]_{i,j}
        :=
        \begin{cases}
          g_{\exp_x(tw)}
          \left(
               \left(
                d\exp_x
                \right)_{tw} w,
                \left(
                d\exp_x
                \right)_{tw} w
          \right) = 1  & i = j = 0 \\
          g_{\exp_x(tw)}
          \left(
               \left(
                d\exp_x
                \right)_{tw} v_{i,w},
                \left(
                d\exp_x
                \right)_{tw} v_{j,w}
          \right) & i \neq 0 \text{ and } j \neq 0\\
          0
          =
          g_{\exp_x(tw)}
          \left(
               \left(
                d\exp_x
                \right)_{tw} w,
                \left(
                d\exp_x
                \right)_{tw} v_{j,w}
          \right) & \text{ otherwise}
        \end{cases}
    \end{align}
we have that $J_x(tw) = \sqrt{\det A(tw)}$. By Lemma \ref{lem:exponentialMapDerivative} there exists $C_{\alpha, \dist_g}>0$ such that for all $0 \leq t \leq r$
\begin{align}
    \left(
    \frac{\sin\left(
        \sqrt{C_{K,r}(x)}t
        \right)}{\sqrt{C_{K,r}(x)}t}
    \right)^2
    \leq
        g_{\exp_x\left(tw\right)}
        \left(
            (d\exp_x)_{tw}
            u,
            (d\exp_x)_{tw}
            u
        \right)
    \leq
    \left(
        \frac{\sinh\left(\sqrt{C_{K,r}(x)}t\right)}{ \sqrt{C_{K,r}(x)} t}
    \right)^2,
\end{align}
where $C_{K,r}(x) = \max\left\{\sup_{y \in B_r(x)}K(y),1\right\}$. Thus we must have
\begin{equation}
\label{eq:volumeDesnityIneq1}
    \left(
        \frac{\sin\left(\sqrt{C_{K,r}(x)}t\right)}{ \sqrt{C_{K,r}(x)} t}
    \right)^{m-1}
    -
    1
\leq
    J_x(tw)
    -
    1
\leq
    \left(
        \frac{\sinh\left(\sqrt{C_{K,r}(x)}t\right)}{ \sqrt{C_{K,r}(x)} t}
    \right)^{m-1}
    -
    1.
\end{equation}
Consider now the functions given by $f(x) = \left(\frac{\sin(x)}{x}\right)^{m-1} - 1$ and $h(x) = \left(\frac{\sinh(x)}{x}\right)^{m-1} - 1$ extended continuously to zero by assigning them the value zero. Their first derivatives vanish at $x=0$, thus there exists $C>0$ and $\delta > 0$ such that for all $|x| \leq \delta$ we have
    \begin{equation}
    \label{eq:ineqAuxilliaryFunctions}
        \left|f(x)\right| \leq Cx^2,
        \quad
        \left|
            h(x)
        \right| \leq Cx^2.
    \end{equation}
Furthermore $\sqrt{C_{K,r}(x)}t \leq \max\left\{\frac{\sqrt{C_K}}{r\left( C_{\alpha,\dist_g} - 1\right)^{\alpha/2}} t,t
    \right\}
    \leq 
    \max\left\{\frac{\sqrt{C_K}}{\left( C_{\alpha,\dist_g} - 1\right)^{\alpha/2}},r_0\right\}$, and so by choosing $C_{\alpha ,\dist_g} > 0$ large enough and $r_0$ small enough we have $\sqrt{C_{K,r}(x)}t \in [-\delta,\delta]$. Thus combining equations \eqref{eq:volumeDesnityIneq1} and \eqref{eq:ineqAuxilliaryFunctions}, we get
\begin{equation}
\label{eq:boundOnJacobiField}
    -C
    C_{K,r}(x)t^2
    \leq
    J_x(tw)-1
    \leq
    C
    C_{K,r}(x)t^2.
\end{equation}

Notice that since $0 < r \leq 1$, if $\alpha > 0$, we have
\begin{equation}
    r 
\leq 
    r^{\min\left\{1,\frac{2}{\alpha}\right\}}
\leq
    \frac{\overline{\dist}_g(x,S)}{C_{\alpha, \dist_g}}.
\end{equation}   
Thus
\begin{align}
    \left|
        J_x(tw)-1
    \right|
\leq&
    CC_{K,r}(x)t^2\\
\leq&
\max\left\{
    \sup_{y \in B_r(x)}
    \frac{CC_K}
    {\overline{\dist}_g(y,S)^\alpha}t^2,
    Ct^2
\right\}\\
\leq&
\max\left\{
    \frac{CC_K}
    {\left(
        \overline{\dist}_g(x,S)
        -
        r
    \right)^\alpha}t^2,
    Ct^2
\right\}\\
\leq&
\max\left\{
    \frac{CC_K}
    {\left(
        \overline{\dist}_g(x,S)
        -
        \frac{\overline{\dist}_g(x,S)}{C_{\alpha,\dist_g}}
    \right)^\alpha}
    t^2
    ,
    Ct^2
\right\}\\
\leq&
\max\left\{
    \frac{CC_K}
    {\left(
        1
        -
        \frac{1}{C_{\alpha,\dist_g}}
    \right)^\alpha}
    \frac{1}{\overline{\dist}_g(x,S)^\alpha}
    t^2,
    \frac{\text{diam}\left(\overline{M}\right)^\alpha C}
    {\overline{\dist}_g(x,S)^\alpha}t^2
\right\},
\end{align} 
so taking $C_J = \max\left\{\frac{CC_K}
    {\left(
        1
        -
        \frac{1}{C_{\alpha,\dist_g}}
    \right)^\alpha},\text{diam}\left(\overline{M}\right)^\alpha C\right\}$ we conclude equation \eqref{eq:ineqForVolumeForm}. Furthermore, since $\overline{\dist}_g(x,S) \geq C_{\alpha,\dist_g}r^{\min\left\{1,\frac{2}{\alpha}\right\}}$ and $0<t<r<r_0<1$ we have for $\alpha > 2$
    \begin{equation}
        \left|
            J_x(tw)-1
        \right|
    \leq
    \max\left\{
        \frac{CC_K}
        {\left(1-\frac{1}{C_{\alpha,\dist_g}}\right)^\alpha}
        \frac{t^2}{C_{\alpha,\dist_g}^\alpha r^2}
        ,
        Ct^2
    \right\}
    \leq
        \max\left\{
        \frac{CC_K}
        {\left(1-\frac{1}{C_{\alpha,\dist_g}}\right)^\alpha}
        \frac{1}{C_{\alpha,\dist_g}^\alpha}
        ,
        Cr_0
        \right\}.
    \end{equation}
    Thus by taking $C_{\alpha,\dist_g}$ large enough so that $\frac{CC_K}
        {\left(1-\frac{1}{C_{\alpha,\dist_g}}\right)^\alpha}
        \frac{1}{C_{\alpha,\dist_g}^\alpha} < \frac{1}{2}$ and $r_0$ small enough so that $Cr_0^2<\frac{1}{2}$, we conclude \eqref{eq:ineqForVolumeForm2}. For $0<\alpha \leq 2$ we have
    \begin{equation}
        \left|
            J_x(tw)-1
        \right|
    \leq
    \max\left\{
        \frac{CC_K}
        {\left(1-\frac{1}{C_{\alpha,\dist_g}}\right)^\alpha}
        \frac{t^2}{C_{\alpha,\dist_g}^\alpha r^{\alpha}}
        ,
        Ct^2
    \right\}
    \leq
        \max\left\{
        \frac{CC_K}
        {\left(1-\frac{1}{C_{\alpha,\dist_g}}\right)^\alpha}
        \frac{1}{C_{\alpha,\dist_g}^\alpha}
        ,
        Cr_0^2
        \right\},
    \end{equation}
    and taking $C_{\alpha,\dist_g}$ large enough and $r_0$ small enough as before concludes \eqref{eq:ineqForVolumeForm2}.
    
    In the case $\alpha = 0$, we can take $r_0>0$ small enough so that $C_{K,r}(x)r_0^2\leq \max\left\{CC_K,C\right\}r_0^2 \leq \min\left\{\delta^2,\frac{1}{2}\right\}$. Then by estimates similar to those in \eqref{eq:boundOnJacobiField}, one obtains
    \begin{equation}
        \left|
        J_x(tw)-1
        \right|
        \leq
        \max
        \left\{
        CC_Kt^2
        ,
        Ct^2
        \right\}
        \leq
        \max
        \left\{
        CC_Kr_0^2
        ,
        Cr_0^2
        \right\}
        <
        \frac{1}{2}.
    \end{equation} 
    Thus we conclude both \eqref{eq:ineqForVolumeForm} and \eqref{eq:ineqForVolumeForm2}.

    Furthermore by choosing $C_{\alpha,\dist_g} \geq C_{\alpha, \inj}$ and $r_0 < r_0^{\text{inj}}$, by Lemma \ref{lem:hessianAndInjecRadius}, the map $\exp_x : B_r(0_x) \rightarrow B_r(x)$ is a diffeomorphism, thus by a change of variables we obtain
    \begin{align}
        \mu(B_r(x))
    =
        \int_{B_r(x)}
        1
        d\vol_g(y)
    =
        \int_{B_r(0_x)}
        J_x(y)
        d\mathcal{L}^m(y)
    \leq
        \mathcal{L}^m\left(\mathbb{B}^m_1(0)\right)\frac{3}{2}r^m,
    \end{align}
    and
    \begin{equation}
        \mu(B_r(x))
    =
        \int_{B_r(x)}
        1
        d\vol_g(y)
    =
        \int_{B_r(0_x)}
        J_x(y)
        d\mathcal{L}^m(y)
    \geq
        \mathcal{L}^m\left(\mathbb{B}^m_1(0)\right)\frac{r^m}{2},
    \end{equation}
    This completes the proof.
\end{proof}

The next proposition can be seen as a way of obtaining an Ahlfors condition on the completion $(\overline{M},\overline{d},\mu)$ from an Ahlfors condition just for the points in $S = \overline{M}\setminus M$.
\begin{prop}
\label{prop:spaceIsDoubling}
    Suppose that $(M,g)$ is a manifold such that the completion $(\overline{M},\overline{\dist}_g)$ is compact, that the curvature satisfies
    \begin{equation}
        K(x)
    \leq
        \frac
        {C_K}
        {\overline{\dist}_g(x,S)^2},
    \end{equation}
    satisfies the injectivity condition \eqref{eq:InjectivityCondition} with $\alpha = 2$, and there exist $\Tilde{c}_B, C_B, r_0 > 0$ such that for all $0< r < r_0$ and $x \in S = \overline{M} \setminus M$
    \begin{equation}
    \label{eq:AhlforsOnSingularity}
        \Tilde{c}_B r^m
    \leq
        \mu(B_r(x)) 
    \leq
        C_B r^m
    \end{equation}
    then the metric measure space $(\overline{M}, \overline{\dist}_g, \mu)$ is locally Alhfors, that is there exists $\overline{r}_0,\Tilde{c}_\vol, C_\vol > 0$ such that for all $0 \leq r < \overline{r}_0$ and $x \in \overline{M}$
    \begin{equation}
    \label{eq:ahlforsOnAllManifold}
        \Tilde{c}_\vol{r^m}
    \leq
        \mu(B_r(x))
    \leq
        C_\vol r^m.
    \end{equation}
    In particular it is a locally doubling space.
\end{prop}
\begin{proof}
    First, we prove the upper bound. By Lemma \ref{lem:volumeFormExp} with $\alpha = 2$, there exists $C_{\alpha, \dist_g} > 0, \tilde{r}_0$ such that if $\overline{\dist}_g(x,S) \geq C_{\alpha, \dist_g} r$ for $0 \leq r \leq \Tilde{r}_0<1$ then
    \begin{equation}
        \mu(B_r(x))
    \leq
        \mathcal{L}^m\left(\mathbb{B}^m_1(0)\right)
        \frac{3}{2}r^m.
    \end{equation}
    On the other hand if $x \in \overline{M}$ is such that $\overline{\dist}_g(x,S) < C_{\alpha, \dist_g} r$, then there exists $\overline{y} \in S$ such that
    \begin{equation}
        \overline{\dist}_g(x,\overline{y})
    \leq
        C_{\alpha, \dist_g} r,
    \end{equation}
    and so by triangle inequality
    \begin{equation}
        B_r(x)
    \subset
        B_{(C_{\alpha,\dist_g}+1)r}(\overline{y}).
    \end{equation}
    Thus for $r \leq \frac{r_0}{(C_{\alpha, \dist_g} + 1)}$ by condition \eqref{eq:AhlforsOnSingularity} we have
    \begin{equation}
        \mu(B_r(x))
    \leq
        \mu(B_{(C_{\alpha,\dist_g}+1)r}(\overline{y}))
    \leq
        C_B\left(C_{\alpha, \dist_g} + 1\right)^m r^m.
    \end{equation}
    Thus we have the upper bound of \eqref{eq:ahlforsOnAllManifold} by taking $C_\vol := \max\left\{\frac{3\mathcal{L}^m\left(\mathbb{B}^m_1(0)\right)}{2}, C_B\left(C_{\alpha, \dist_g} + 1\right)^m\right\}$ and $\overline{r}_0 \leq \min\left\{\Tilde{r}_0, \frac{r_0}{(C_{\alpha, \dist_g} + 1)}\right\}$.

    We now prove the lower bound. Given $r\leq1$, consider $\Tilde{r} = \frac{r}{C_{\alpha,\dist_g}2}$. Then by Lemma \ref{lem:volumeFormExp}, for $\overline{\dist}_g(x,S) \geq C_{\alpha, \dist_g}\Tilde{r} = \frac{r}{2}$ we have
    \begin{equation}
        \frac{\mathcal{L}^m\left(\mathbb{B}^m_1(0)\right)}{2}
        \frac{1}{(C_{\alpha,\dist_g}2)^m}
        r^m
    =
        \frac{\mathcal{L}^m\left(\mathbb{B}^m_1(0)\right)}{2}\left(\Tilde{r}\right)^m
    \leq
        \mu(B_{\Tilde{r}}(x))
    \leq
        \mu(B_r(x)),
    \end{equation}
    since $C_{\alpha,\dist_g} > 2$. On the other hand if $\overline{\dist}_g(x,S) \leq \frac{r}{2}$, then there exists $\overline{y} \in S$ such that
    \begin{equation}
        \overline{\dist}_g(x,\overline{y}) \leq \frac{r}{2},
    \end{equation}
    thus by triangle inequality
    \begin{equation}
        B_{\frac{r}{2}}(\overline{y})
    \subset 
        B_r(x).
    \end{equation}
    For $r \leq r_0$, by \eqref{eq:AhlforsOnSingularity} we have 
    \begin{equation}
        \frac{\Tilde{c}_B}{2^m}
        r^m
    \leq
        \mu(B_{\frac{r}{2}}(\overline{y}))
    \leq
        \mu(B_r(x)),
    \end{equation}
    and so the inequality is true for $\overline{r}_0 \leq r_0$ and $\Tilde{c}_\vol := \min\left\{\frac{\tilde{c}_B}{2^m},\frac{\mathcal{L}^m\left(\mathbb{B}^m_1(0)\right)}{2(C_{\alpha,\dist_g}2)^m}\right\}$, finishing the proof.
\end{proof}
The following lemma shows how to obtain a second order form for the Sobolev inequality from the first order version.
    \begin{lemma}
    \label{lem:sobolevIneq2derivatives}
        Suppose that the space $(M,g)$ satisfies property \eqref{eq:SobolevIneq} with exponent $1\leq p < m$ and $\vol_g\left(M\right) < \infty$. Then given $p^* > 0$ such that $\frac{1}{p^*} = \frac{1}{p}-\frac{1}{m}$ there exists a constant $C_{\mathcal{S}, p} > 0$ such that for all $f \in W^{2,p}(M,g)$ we have
        \begin{equation}
            \left\|
                \left\|  
                    \nabla f
                \right\|_g
            \right\|_{L^{p^*}(M,\vol_g)}
        \leq
            C_{\mathcal{S}, p}
            \left(
            \left\|
                \left\|
                    Hf
                \right\|_g
            \right\|_{L^p(M,\vol_g)}
            +
            \left\|
                \left\|
                \nabla f
                \right\|_g
            \right\|_{L^p(M,\vol_g)}
            \right).
        \end{equation}  
    \end{lemma}
    \begin{proof}
         Given $f \in C^\infty(M) \cap W^{2,p}(M,g)$ and $\epsilon>0$, define $h_\epsilon(x) := \sqrt{g_x\left(\nabla f(x), \nabla f(x)\right) + \epsilon^2}$ and $h(x) := \sqrt{g_x\left(\nabla f(x), \nabla f(x)\right)}$. Notice that $h_\epsilon \in C^\infty(M)$ and, since $\vol_g(M) < \infty$ and $h \in L^p(M,\vol_g)$, $h_\epsilon \in L^p\left(M, \vol_g\right)$ as well for all $\epsilon>0$. In fact, given $u \in T_xM$ such that $g_x(u,u) = 1$ we have
         \begin{align}
             \left|
                \nabla_u h_\epsilon(x)
             \right|
        =&
            \left|
            \frac
            {
                g_x\left(
                    \nabla_u
                    \nabla f(x)
                    ,
                    \nabla f(x)
                \right)
            }
            {\sqrt{g_x\left(\nabla f(x),\nabla f(x)\right) + \epsilon^2}
            }
            \right|\\
        \leq&
            \sup_{\{v \in T_x M: g_x(v,v) = 1\}}
            g_x(\nabla_u\nabla f(x), v)
        \leq
            \|H f(x)\|_g.
            \label{eq:ineqWithHessian}
         \end{align}     
        This implies that $\left\|\nabla h_\epsilon\right\|_g \leq \|Hf\|_g \in L^p(M,\vol_g)$. Since $h\in L^p(M,\vol_g)$, it follows that $h_\epsilon\in W^{1,p}(M,g)$. Furthermore, $h_\epsilon$ converges pointwise to $h$. Defining for each $x \in M$ and $u \in T_xM$
        \begin{equation}
            A_u(x)
            :=
            \begin{cases}
                \frac
            {
                g_x\left(
                    \nabla_u
                    \nabla f(x)
                    ,
                    \nabla f(x)
                \right)
            }
            {\sqrt{g_x\left(\nabla f(x),\nabla f(x)\right)}
            }
            &
            \text{ if }
            \left\|\nabla f(x)\right\|_g > 0\\
            0
            &
            \text{ if }
            \left\|\nabla f(x)\right\|_g = 0
            \end{cases}
        \end{equation}
        we have $\nabla_u h_\epsilon(x)\rightarrow A_u(x)$. Since $T_xM$ is finite-dimensional, this implies that $\sup_{u \in T_xM : g_x(u,u) = 1}\left| A_u(x) - \nabla_u h_\epsilon (x) \right|\rightarrow 0$ for all $x \in M$. Therefore, \eqref{eq:ineqWithHessian} and the dominated convergence theorem imply that $h_\epsilon$ is a Cauchy sequence for the $W^{1,p}(M,g)$ norm as $\epsilon \rightarrow 0$ whose limit has representative $h \in W^{1,p}(M,g)$ since $h_\epsilon$ converges to $h$ in $L^p\left(M, \vol_g\right)$. Finally, we apply \eqref{eq:SobolevIneq} with exponent $p$ to conclude that, given $p^*>0$ such that $\frac{1}{p^*} = \frac{1}{p} - \frac{1}{m}$
        \begin{align}
            \left\|
                \left\|
                    \nabla f
                \right\|_{g}
            \right\|_{L^{p^*}(M, \vol_g)}
        =&
            \left\|
                h
            \right\|_{L^{p^*}(M,\vol_g)}\\
        \leq&
            C_{\mathcal{S},p}
            \left(
            \left\|
                \left\|
                    \nabla h
                \right\|_g
            \right\|_{L^p(M,\vol_g)}
            +
            \left\|
                h
            \right\|_{L^p(M,\vol_g)}
            \right)\\
        \leq&
            C_{\mathcal{S},p}
            \left(
            \left\|
                \left\|
                    H f
                \right\|_g
            \right\|_{L^p(M,\vol_g)}
            +
            \left\|
            \left\|
                \nabla f
            \right\|_g
            \right\|_{L^p(M,\vol_g)}
            \right)
            <\infty.
        \end{align}
        The result then follows from the density of $C^\infty(M)$ in $W^{2,p}(M,g)$. This completes the proof.
    \end{proof}
    \begin{lemma}
    \label{lem:finiteMeasureLemma}
        Suppose that $(M,g)$ is a Riemannian manifold such that the completion $(\overline{M},\overline{\dist}_g, \mu)$ is compact and there exists $r_0 > 0 $ such that
        \begin{equation}
            \mu(T_{r_0}(S)) < \infty,
        \end{equation}
        then $\mu$ is a finite measure.
    \end{lemma}
    \begin{proof}
        Since $g$ is a $C^\infty$ metric on $M$ and $M \setminus T_{r_0}(S)$ is compact, we have
        \begin{equation}
            \mu(M \setminus T_ {r_0}(S))<\infty.
        \end{equation}
        Thus we have
        \begin{equation}
            \mu(M)
        =
            \mu(M \setminus T_{r_0}(S))
        +
            \mu(T_{r_0}(S))
        <
            \infty.
        \end{equation}
    \end{proof}

\section{AMV Laplacian convergence}

Now we will define and study the Hardy--Littlewood maximal function. A reference for these operators on doubling metric measure spaces can be found in \cite{Heinonen2001}. This function is closely related to the AMV Laplacian operators.
\begin{definition}
    Let $(X,d,\nu)$ be a locally doubling metric measure space in which every ball has positive finite measure. Given $r_0>0$ and $f \in L^1_{\text{loc}}(X,\nu)$, define the restricted Hardy-Littlewood maximal function as
    \begin{equation}
        M_{r_0}\left(
            f
        \right)
        (x)
    =
        \sup_{0 < r < r_0}
        \frac{1}{\nu(B_r^d(x))}
        \int_{B_r^d(x)}
            \left|
                f(y)
            \right|
            d\nu(y).
    \end{equation}  
\end{definition}
\begin{theorem}
\label{thm:hardyLittleWoodThm}
    Suppose that $(X,d,\nu)$ is a locally doubling metric measure space in which every ball has positive finite measure. Then there exist $R_0,C_{p,\max}>0$ such that for all $1<p<\infty$ and $0<r<R_0$
    \begin{equation}
    \label{eq:HardyLittlewoodIneq}
        \left\|
            M_{r}(f)
        \right\|_{L^p(X,\nu)}
    \leq
        C_{p,\max}
        \|f\|_{L^p(X,\nu)}.
    \end{equation}
\end{theorem}
\begin{prop}
\label{prop:AMVOperatorIsBounded}
    Suppose that $(X,d,\nu)$ is a compact locally doubling metric measure space with $0<\nu(B_r(x))<\infty$ for all $r>0$ and $x \in X$. Then there exists $r_0>0$ such that for all $0< r <r_0$ the operators 
    $\Delta_r : L^p(X,\nu) \rightarrow L^p(X,\nu)$ as defined in \ref{def:amvLaplacian} are bounded for $1<p<\infty$.
\end{prop}
\begin{proof}

By Theorem~\ref{thm:hardyLittleWoodThm}, there exist $r_0$ and $C_{p,\max}>0$ such that for all $r \leq r_0$
    \begin{equation}
        \left\|
            M_{r}(f)
        \right\|_{L^p(X,\nu)}
    \leq
        C_{p,\max}
        \|f\|_{L^p(X,\nu)}
    \end{equation}
    Thus for $r < r_0$ we have
    \begin{align}
        \left\|
        \Delta_r f
    \right\|_{L^p(X,\nu)}
=&
\left(
    \int_{X}
        \left|
            \frac{f(x)}{r^2}
        -
            \frac{1}{r^2\nu(B_r^d(x))}
            \int_{B_r^d(x)}
            f(y)
            d\nu(y)
        \right|^p
        d\nu(x)
\right)^\frac{1}{p}\\
\leq&
    \frac{\|f\|_{L^p(X,\nu)}}{r^2}
    +
    \left(
    \int_{X}
        \left|
            \frac{1}{r^2\nu(B_r^d(x))}
            \int_{B_r^d(x)}
            f(y)
            d\nu(y)
        \right|^p
        d\nu(x)
    \right)^\frac{1}{p}\\
    \leq&
        \frac{\|f\|_{L^p(X,\nu)}}{r^2}
    +
    \left(
    \int_{X}
        \left|
            \frac{1}{r^2}
            M_{r}f(x)
        \right|^p
        d\nu(x)
    \right)^\frac{1}{p}\\
    \leq&
        \frac{\|f\|_{L^p(X,\nu)}}{r^2}
    +
        \frac{C_{p,\max} \|f\|_{L^p(X,\nu)}}{r^2},
    \end{align}
    This completes the proof.
\end{proof}

Now we study bounds for the AMV Laplacian operators on the metric measure space $(\overline{M},\overline{\dist}_g, \mu)$.
\begin{prop}
\label{prop:boundOnAMVLaplacian}
    Suppose $(M,g)$ is a Riemannian manifold satisfying \eqref{eq:curvatureCondition} and \eqref{eq:InjectivityCondition}, and suppose that its completion $(\overline{M},\overline{\dist}_g,\mu)$ is a compact locally doubling metric measure space.

    Let $f:\overline{M}\rightarrow\mathbb{R}$ satisfy $f\in C^\infty(M)$, and, for some $0<\beta\leq1$, suppose $f$ is $\beta$-H\"older continuous with H\"older constant $L_\beta>0$. Then there exist $C_{\alpha,\dist_g},C_{\Delta},r_0>0$ such that for all $0<r<r_0$ and $x \in M$ we have
    \begin{align}
    \label{eq:mainInequality}
        \left|
            \Delta_r f(x)
        \right|
    &\leq
        C_{\Delta}
        \left(
            \frac
            {L_\beta}
            {r^{2-\beta}}
            \chi_{
            T_{C_{\alpha,\dist_g} r^{\min\{1,2/\alpha\}}}(S)}(x)
        +
        \right.\\
        &\quad
        \left.\chi_{T_{C_{\alpha,\dist_g} r^{\min\{1,2/\alpha\}}}(S)^c}(x)
            \left(
            \frac{r M_{r}\left(\left\|\nabla f\right\|_g\right)(x)}
            {\overline{\dist}_g(x ,S)^\alpha}
        +
            M_{r}\left(
                \left\|
                    Hf
                \right\|_g
            \right)
            (x)
            \right)
        \right).
    \end{align}
    In particular we have
    \begin{equation}
    \label{eq:holderIneq}
        \left|
            \Delta_r f(x)
        \right|
    \leq
        C_{\Delta}
    \left(
        \frac{L_\beta C_{\alpha,\dist_g}^{
            \max\left\{
            \frac{(2-\beta)\alpha}{2},
            (2-\beta)
            \right\}}}{\overline{\dist}_g(x,S)^{
            \max\left\{
            \frac{(2-\beta)\alpha}{2},
            (2-\beta)
            \right\}}}
    +
        \frac{1}{C_{\alpha,\dist_g}^{\frac{\alpha}{2}}}
            \frac{1}{\overline{\dist}_g(x,S)^{\frac{\alpha}{2}}}
        M_{r_0}\left(
        \left\|
            \nabla f
        \right\|_g
        \right)
        (x)
        +
        M_{r_0}\left(
            \left\|Hf\right\|_g
        \right)
        (x)
    \right)
    \end{equation}
    and if $f$ is Lipschitz then
    \begin{equation}
    \label{eq:lipschitzIneq}
        \left|
            \Delta_r f(x)
        \right|
    \leq
        C_{\Delta}
        \left(
            \frac{L_1 C_{\alpha,\dist_g}^{
            \max\left\{
            1,
            \frac{\alpha}{2}
            \right\}
            }}{\overline{\dist}_g(x,S)^{
            \max\left\{
            1,
            \frac{\alpha}{2}
            \right\}      
            }}
    +
        \frac{1}{C_{\alpha,\dist_g}^{\frac{\alpha}{2}}}
            \frac{1}{\overline{\dist}_g(x,S)^{\frac{\alpha}{2}}}L_{1}
        +
        M_{r_0}\left(
            \left\|Hf\right\|_g
        \right)
        (x)
        \right).
    \end{equation}
\end{prop}
\begin{proof}
        We consider $C_{\alpha,\dist_g}>0$ and $0<r_0<R_0$ so that Theorem \ref{thm:hardyLittleWoodThm} and Lemmas \ref{lem:hessianAndInjecRadius}, \ref{lem:volumeFormExp} are satisfied. For $x \in M$ such that for $0<r<r_0 $ and $\overline{\dist}_g(x,S) \geq C_{\alpha,\dist_g} r^{\min\left\{1,\frac{2}{\alpha}\right\}}$ we have
        $\exp_x : B_r(0_x) \rightarrow B_r(x)$ is a diffeomorphism and the volume form
        \begin{equation}
            \left(\left(\exp_x\right)^*d\vol_g\right)_y
        =
            J_x(y)dy_1\wedge\cdots \wedge dy_m,
        \end{equation}
        satisfies equations \eqref{eq:ineqForVolumeForm} and \eqref{eq:ineqForVolumeForm2}. By a change of variables and the fundamental theorem of calculus we obtain
        \begin{align}
            \Delta_r f(x)
    &=
        \frac{1}{r^2\mu(B_r(x))}
        \int_{B_r(0_x)}
            \left(
                f\left(\exp_x\left(y\right)\right)
            -
                f\left(\exp_{x}\left(0_x\right)\right)
            \right)
            J_{x}(y)
            d\mathcal{L}^m(y)\\
    &=
        \frac{1}{r^2\mu(B_r(x))}
        \int_{B_r(0_x)}
            \left(
                \int_0^1
                d(f\circ \exp_{x})_{ty}
                \cdot y
                dt
            \right)
            J_{x}(y)
            d\mathcal{L}^m(y).
        \end{align}
        We now define the auxiliary functions
        \begin{align}
            G_r(x) 
    &= 
        \frac{1}{r^2\mu(B_r(x))}
        \int_{B_r(0_x)}
            \left(
                \int_0^1
                d(f\circ \exp_{x})_{ty}
                \cdot y
            -
                d(f\circ \exp_{x})_{0_x}
                \cdot y
            \right)
            dt
            d\mathcal{L}^m(y),\\
        F_r(x) 
    &= 
        \frac{1}{r^2\mu(B_r(x))}
        \int_{B_r(0_x)}
            \left(
                \int_0^1
                d(f\circ \exp_{x})_{ty}
                \cdot y
            \right)
            \left(
                J_{x}(y)
            -
                1
            \right)
            dt
            d\mathcal{L}^m(y),
        \end{align}
        and observe that, since
        \begin{equation}
            \frac{1}{r^2\mu(B_r(x))}
        \int_{B_r(0_x)}
            \left(
                \int_0^1
                d(f\circ \exp_{x})_{0_x}
                \cdot y
                dt
            \right)
            d\mathcal{L}^m(y)
            =
            0,
        \end{equation}
        we have
        \begin{equation}
            \Delta_r f(x) = G_r(x) + F_r(x).
        \end{equation}
        First we bound $G_r(x)$. Using the fundamental theorem of calculus, a change of variables, and Lemmas \ref{lem:hessianAndInjecRadius} and \ref{lem:volumeFormExp}
        \begin{align}
        \left|
            G_r(x)
        \right|
        &= 
        \left|
        \frac{1}{r^2\mu(B_r(x))}
        \int_{B_r(0_x)}
            \left(
                \int_0^1
                d(f\circ \exp_{x})_{ty}
                \cdot y
            -
                d(f\circ \exp_{x})_{0_x}
                \cdot y
            dt
            \right)
            d\mathcal{L}^m(y)
        \right|\\
        &= 
        \left|
        \frac{1}{r^2\mu(B_r(x))}
        \int_{B_r(0_x)}
            \left(
                \int_0^1
                \int_0^1
                \frac{d}{ds}
                d
                \left(
                    f\circ \exp_x
                \right)_{sty}
                \cdot y
                ds
            dt
            \right)
            d\mathcal{L}^m(y)
        \right|\\
        &=\left|
        \frac{1}{r^2\mu(B_r(x))}
        \int_{B_r(0_x)}
            \left(
                \int_0^1
                \int_0^1
                t
                \left.
                \frac{d^2}{d\alpha^2}\right|_{\alpha = st}
                \left(
                    f\circ \exp_x
                \right)(\alpha y)
            ds
            dt
            \right)
            d\mathcal{L}^m(y)
        \right|\\
        &\leq
        \frac{1}{r^2\mu(B_r(x))}
        \int_{B_r(0_x)}
            \left(
                \int_0^1
                \int_0^1
                t
                \|Hf(\exp_x(sty))\|_g
                g_x(y,y)
                ds
                dt
            \right)
            d\mathcal{L}^m(y)\\
        &=
            \frac{1}{r^2\mu(B_r(x))}
            \int_0^1
            \int_0^1
            \frac{1}{(st)^{m+2}}
            \int_{B_{str}(0_x)}
            t
                \|Hf(\exp_x(y))\|_g
                g_x(y,y)
            d\mathcal{L}^m(y)
            dsdt\\
        &\leq
            \frac{2}
            {\mathcal{L}^m\left(
                \mathbb{B}^m_1(0)
            \right)}
            \int_0^1
            \int_0^1
            \frac{1}{(str)^{m}(str)^2}
            \int_{B_{str}(0_x)}
                \|Hf(\exp_x(y))\|_g
                g_x(y,y)
            d\mathcal{L}^m(y)
            dsdt\\
        &\leq
            \frac{6}{2}
            \int_0^1
            \int_0^1
            \frac{(str)^2}{\mu(B_{str}(x))(str)^2}
            \int_{B_{str}(0_x)}
                \|Hf(\exp_x(y))\|_g
            d\mathcal{L}^m(y)
            dsdt\\
        &\leq
            6
            \int_0^1
            \int_0^1
            \frac{1}{\mu(B_{str}(x))}
            \int_{B_{str}(0_x)}
                \|Hf(\exp_x(y))\|_g
            J_x(y)
            d\mathcal{L}^m(y)
            dsdt\\
        &=
            6
            \int_0^1
            \int_0^1
            \frac{1}{\mu(B_{str}(x))}
            \int_{B_{str}(x)}
                \|Hf(y)\|_g
            d\mu(y)
            dsdt\\
        &\leq 
            6\int_0^1\int_0^1 M_r\left(\left\|Hf\right\|_g\right)(x) dsdt
        \leq
            6M_r\left(\left\|Hf\right\|_g\right)(x).
        \end{align}
        Now we bound $F_r(x)$,
        \begin{align}
            \left|
                F_r(x)
            \right|
        &=
        \left|
            \frac{1}{r^2\mu(B_r(x))}
            \int_{B_r(0_x)}
                \left(
                    \int_0^1
                    d(f\circ \exp_{x})_{ty}
                    \cdot y
                    dt
                \right)
                \left(
                    J_{x}(y)
                -
                    1
                \right)
                d\mathcal{L}^m(y)
            \right|\\
        &\leq
            \frac{1}{r^2\mu(B_r(x))}
            \int_{B_r(0_x)}
                \left(
                    \int_0^1
                    \|\nabla f\left(\exp_{x}\left(ty\right)\right)\|_g
                    g_x(y,y)^\frac{1}{2}
                    dt
                \right)
                \left|
                    J_{x}(y)
                -
                    1
                \right|
                d\mathcal{L}^m(y)\\
            &\leq
            \frac{1}{r\mu(B_r(x))}
            \int_{B_r(0_x)}
                \left(
                    \int_0^1
                    \|\nabla f\left(\exp_{x}\left(ty\right)\right)\|_g
                    dt
                \right)
                \frac{C_J}{\overline{\dist}_g(x,S)^\alpha}
                g_x(y,y)
                d\mathcal{L}^m(y)\\
            &\leq
                \frac{r}{\mu(B_r(x))}
            \int_{B_r(0_x)}
                \left(
                    \int_0^1
                    \|\nabla f\left(\exp_{x}\left(ty\right)\right)\|_g
                    dt
                \right)
                \frac{C_J}{\overline{\dist}_g(x,S)^\alpha}
                d\mathcal{L}^m(y)\\
            &=
                \frac{r}{\mu(B_r(x))}
                \frac{C_J}{\overline{\dist}_g(x,S)^\alpha}
            \int_{0}^1
            \frac{1}{t^m}
            \int_{B_{tr}(0_x)}
                    \|\nabla f\left(\exp_{x}\left(y\right)\right)\|_g
                d\mathcal{L}^m(y)dt\\
            &\leq
                \frac{2rC_J}{\mathcal{L}^m\left(
                \mathbb{B}^m_1(0)
            \right)\overline{\dist}_g(x,S)^\alpha }
            \int_{0}^1
            \frac{1}{(rt)^m}
            \int_{B_{tr}(0_x)}
                    \|\nabla f\left(\exp_{x}\left(y\right)\right)\|_g
                d\mathcal{L}^m(y)
                dt\\
            &\leq
                \frac{6rC_J}{\overline{\dist}_g(x,S)^\alpha }
            \int_{0}^1
            \frac{1}{\mu(B_{tr}(x))}
            \int_{B_{tr}(0_x)}
                    \|\nabla f\left(\exp_{x}\left(y\right)\right)\|_g
                J_x(y)
                d\mathcal{L}^m(y)
                dt\\
            &\leq
                \frac{6rC_J}{\overline{\dist}_g(x,S)^\alpha }
            \int_{0}^1
            \frac{1}{\mu(B_{tr}(x))}
            \int_{B_{tr}(x)}
                    \|\nabla f\left(y\right)\|_g
                d\mu(y)
                dt\\
            &\leq
                \frac{6 r C_J}
                {\overline{\dist}_g(x,S)^\alpha}
                M_r\left(\left\|\nabla f\right\|_g\right)(x).
        \end{align}
        Therefore, for $x \in M$ such that $\overline{\dist}_g(x,S) \geq C_{\alpha,\dist_g}r^{\min\left\{1,\frac{2}{\alpha}\right\}}$,
        \begin{equation}
            \left|
                \Delta_r f(x)
            \right|
        \leq
            \frac{6 r C_J}
                {\overline{\dist}_g(x,S)^\alpha}
                M_r\left(\left\|\nabla f\right\|_g\right)(x)
        +
            6M_r\left(\left\|Hf\right\|_g\right)(x).
        \end{equation}      
        On the other hand, for points $x \in M$ such that $\overline{\dist}_g(x,S) < C_{\alpha,\dist_g}r^{\min\left\{1,\frac{2}{\alpha}\right\}}$, since $f$ is $\beta$-H\"older
        \begin{align}
            \left|
                \Delta_r f(x)
            \right|
        =&
            \left|
                \frac{1}{r^2\mu(B_r(x))}
                \int_{B_r(x)}
                \left( f(y) - f(x)\right)d\mu(y)
            \right|\\
        \leq&
            \frac{L_\beta}{r^{2-\beta}},
        \end{align} 
        These inequalities prove \eqref{eq:mainInequality} by taking $C_\Delta  \geq 6C_J + 6$.

        Now notice that
        \begin{equation}
            \frac{L_\beta}{r^{2-\beta}}
            \chi_{
                T_{
                C_{\alpha,\dist_g}
                r^{\min\left\{1,\frac{2}{\alpha} \right\}}
                }
                (S)
                }
        \leq
            \frac{L_\beta C_{\alpha,\dist_g}^{
            \max\left\{
            \frac{(2-\beta)\alpha}{2},
            (2-\beta)
            \right\}}}{\overline{\dist}_g(x,S)^{
            \max\left\{
            \frac{(2-\beta)\alpha}{2},
            (2-\beta)
            \right\}}}.
        \end{equation}

        Furthermore, observe that for $\overline{\dist}_g(x,S) \geq C_{\alpha,\dist_g}r^{\min\left\{1,\frac{2}{\alpha}\right\}}$ we have
        \begin{equation}
            \frac{r}{\overline{\dist}_g(x,S)^\alpha}
        \leq
            \frac{1}{C_{\alpha,\dist_g}^{\frac{\alpha}{2}} \overline{\dist}_g(x,S)^{\alpha - \max\left\{1,\frac{\alpha}{2}\right\}}}\\
        \leq
            \frac{1}{C_{\alpha,\dist_g}^{\frac{\alpha}{2}}}
            \frac{\max\left\{1,\text{diam}(\overline{M})^{\frac{\alpha}{2}-(\alpha-1)}\right\}}{\overline{\dist}_g(x,S)^{\frac{\alpha}{2}}},
        \end{equation}
        From these inequalities and equation \eqref{eq:mainInequality} we conclude equation \eqref{eq:holderIneq} by taking $C_\Delta \geq 6C_J + 6 + 6C_J\max\left\{1,\text{diam}(\overline{M})^{\frac{\alpha}{2}-(\alpha-1)}\right\}$. Furthermore, if $f$ is Lipschitz we obtain inequality \eqref{eq:lipschitzIneq}, finishing the proof.
\end{proof}

\begin{lemma}
    \label{lem:1OverDistIsIntegrable}
    Let $0\leq k<m$, $1\leq p<\infty$, and $0\leq\alpha\leq m$. Suppose $(M,g)$ is a Riemannian manifold such that $(\overline{M},\overline{\dist}_g)$ is compact, and let $S=\overline{M}\setminus M$ have bounded upper $(m,k)$-dimensional Minkowski content. Then the function
    \begin{equation}
        y \in M \mapsto \frac{1}{\overline{\dist}_g(y,S)^\alpha}
    \end{equation}
    is in $L^p\left(\overline{M},\mu\right)$ if $m-k-\alpha p > 0$.
\end{lemma}
\begin{proof}
    The case $\alpha = 0$ is trivial since by Lemma \ref{lem:finiteMeasureLemma} $\overline{M}$ has finite volume. Thus we assume that $\alpha > 0$. Since $S$ has finite $(m,k)$ Minkowski content, there exists $r_0>0$ such that $\mu\left(T_{r_0}\left(S\right)\right) < \infty$, so Lemma~\ref{lem:finiteMeasureLemma} gives $\mu\left(\overline{M}\right)<\infty$.

    Consider the sequence of natural numbers $N_n \in \mathbb{N}$ given by $N_0 = 0$ and $N_n = 2^n$ for $n \geq 1$. Then we have
    \begin{equation}
        \sum_{n=1}^\infty
            {N_n^{-\frac{m-k-\alpha p}{\alpha}}}
        =
        \sum_{n=1}^\infty
            {2^{-n\frac{m-k-\alpha p}{\alpha}}}
        <
            \infty.
    \end{equation}
    For $n\geq1$, let
    \begin{align}
        E_n
    :=
        \left\{
            x \in M:
            \frac{1}{\overline{\dist}_g(x,S)^\alpha}
            >
            N_n
        \right\}
    =
        \left\{
            x \in M:
            \overline{\dist}_g(x,S)
            <
                \frac{1}{N_n^\frac{1}{\alpha}}
        \right\}
    =
        T_{{N_n^{-\frac{1}{\alpha}} }}(S) \cap M,
    \end{align}
    and
    \begin{equation}
        E_0 := \overline{M}.
    \end{equation}
    Define the function
    \begin{equation}
        f(x)
    =
        \sum_{n=0}^\infty
        N_{n+1}
        \chi_{E_n \setminus E_{n+1}}(x).
    \end{equation}
    If $x \in E_n\setminus E_{n+1}$, then
    \begin{equation}
        \frac{1}{\overline{\dist}_g(x,S)^\alpha}
    \leq
        N_{n+1}
    =
        N_{n+1}\chi_{E_n \setminus E_{n+1}}(x)
    \leq
        f(x).
    \end{equation}
    Since $\overline{M} = \cup_{n=0}^\infty \left(E_{n}\setminus E_{n+1}\right)$,
    it follows that, for all $x \in M$, $\frac{1}{\overline{\dist}_g(x,S)^\alpha} \leq f(x)$. 

    Since $S \subset \overline{M}$ has bounded upper $(m,k)$-dimensional Minkowski content, there exist $C_{\mathcal{M}}>0$ and $K\in\mathbb{N}$ such that for $n \geq K$
    \begin{equation}
        \mu(E_n)
    =
        \mu\left(T_{N_n^{-\frac{1}{\alpha}}}\left(S\right)\right)
    \leq
        C_{\mathcal{M}} N_n^{-\frac{m-k}{\alpha}}.
    \end{equation}
    Since $\mu(S) = 0$ we have
    \begin{align}
        \int_{\overline{M}}
            \frac{1}{\overline{\dist}_g(x,S)^{\alpha p}}
            d\mu(x)
    \leq&
        \int_{M}
            f(x)^p
            d\mu(x)
    =
        \int_{M}
            \sum_{n=0}^\infty
            (N_{n+1})^p
            \chi_{E_n \setminus E_{n+1}}(x)
            d\mu(x)\\
    \leq&
        \sum_{n=1}^K N_n^p
        \mu(M)
    +
        \sum_{n=K}^\infty
        (N_{n+1})^p
        \mu(E_n)\\
    \leq&
        \sum_{n=1}^K N_n^p
        \mu(M)
        +
        \sum_{n=K}^\infty
        N_{n+1}^p
        \frac{C_{\mathcal{M}}}{N_n^\frac{m-k}{\alpha}}\\
    =&
        \sum_{n=1}^K N_n^p
        \mu(M)
        +
        \sum_{n=K}^\infty
        2^p
        N_{n}^p
        \frac{C_{\mathcal{M}}}{N_n^\frac{m-k}{\alpha}}\\
    =&
       \sum_{n=1}^K N_n^p
        \mu(M)
        +
        \sum_{n=K}^\infty
        2^p
        C_{\mathcal{M}}
        {N_n^{-\frac{m-k-\alpha p}{\alpha}}}
    <
        \infty,
    \end{align}
    which proves that the function belongs to $L^p\left(\overline{M},\mu\right)$.
\end{proof}

\begin{lemma}
\label{lem:weakConvergenceLemma}
    Let $1<p<\infty$, and suppose that $\mu\left(\overline{M}\right)< \infty$. Consider a sequence of functions $g_n : \overline{M} \rightarrow \mathbb{R}$ such that
    \begin{equation}
        \sup_{n \in \mathbb{N}}
        \|g_n\|_{L^p\left(\overline{M},\mu\right)}
    <
        \infty
    \text{  and }
    \lim_n 
    g_n(x) = g(x)
    \quad
    \mu\text{-a.e.}\quad x \in M
    \end{equation}
    then
    \begin{equation}
        g_n \rightharpoonup_{L^p\left(\overline{M},\mu\right)}
        g.
    \end{equation}
\end{lemma}
\begin{proof}
    Suppose by contradiction that
    \begin{equation}
        g_n \not\rightharpoonup_{L^p\left(\overline{M},\mu\right)}g.
    \end{equation}
    Then there exist $\epsilon>0$ and $\psi \in L^q\left(\overline{M},\mu\right)$ with $\frac{1}{p}+\frac{1}{q}=1$ such that
    \begin{equation}
        \limsup_{n}
        \left|
            \int_{\overline{M}}
            g_n\psi d\mu
        -
            \int_{\overline{M}}
            g
            \psi d\mu
        \right|
    >
        \epsilon.
    \end{equation}
    This implies that there exists a subsequence $n_j$ such that
    \begin{equation}
        \inf_{j \in \mathbb{N}}
        \left|
            \int_{\overline{M}}
            g_{n_j}\psi d\mu
        -
            \int_{\overline{M}}
            g
            \psi d\mu
        \right|
    >
        \epsilon.
    \end{equation}
    Since the sequence $g_n$ is bounded in $L^p\left(\overline{M},\mu\right)$, we can apply the Banach--Alaoglu theorem, the reflexivity of $L^p\left(\overline{M},\mu\right)$, and the Eberlein--Šmulian theorem to conclude that there exists a subsequence $\{\overline{n}_j\}_{j \in \mathbb{N}} \subset \{n_j\}_{j \in \mathbb{N}}$ and $h \in L^p\left(\overline{M},\mu\right)$ such that
    \begin{equation}
        g_{\overline{n}_j} \rightharpoonup_{L^p\left(\overline{M}, \mu\right)}
        h.
    \end{equation}
    Notice also that, by Fatou's lemma, $g \in L^p\left(\overline{M},\mu\right)$.

    \vspace{10mm}

    Now we show that $h = g$ $\mu$-a.e. Let $E_{j,\epsilon}:=\left\{
                x \in \overline{M}:
                \left|
                    g_{\overline{n}_j}
                    -
                    g
                \right|
                >
                \epsilon
            \right\}$. Given $\epsilon > 0$, since $\mu\left(\overline{M}\right)<\infty$ and $g_{\overline{n}_j}$ converges pointwise to $g$, we must have
    \begin{equation}
        \lim_j
        \mu\left(E_{j,\epsilon}\right)
        =
        \lim_{j}
        \mu\left(
            \left\{
                x \in \overline{M}:
                \left|
                    g_{\overline{n}_j}
                    -
                    g
                \right|
                >
                \epsilon
            \right\}
        \right)
        =
        0.
    \end{equation}
    In particular, observe that for any fixed $\psi \in L^q\left(\overline{M}, \mu \right)$, we must have $\left\|\psi\right\|_{L^q\left(E_{j,\epsilon}, \mu\right)} \rightarrow 0$. Thus for any $\psi \in L^q\left(\overline{M}, \mu\right)$, we must have
    \begin{align}
    \left|
        \int_{\overline{M}}
            \left(
                g
                -
                h
            \right)\psi
            d\mu
    \right|
    &=
    \lim_{j}
    \left|
        \int_{\overline{M}}
        \left(
        g
        -
        g_{\overline{n}_j}
        \right)
        \psi
        d\mu
    \right|\\
    &\leq
    \limsup_{j}
    \left|
        \int_{E_{j,\epsilon}}
        \left(
        g
        -
        g_{\overline{n}_j}
        \right)
        \psi
        d\mu
    \right|
    +
    \left|
        \int_{\overline{M} \setminus E_{j,\epsilon}}
        \left(
        g
        -
        g_{\overline{n}_j}
        \right)
        \psi
        d\mu
    \right|\\
    &\leq
    \limsup_j
    \left(
        \left\|
            g
        \right\|_{L^p(\overline{M}, \mu)}
        +
        \left\|
            g_{\overline{n}_j}
        \right\|_{L^p(\overline{M}, \mu)}
    \right)
    \left\|
        \psi
    \right\|_{L^q\left(E_{j,\epsilon}, \mu\right)}
    +
    \epsilon\mu\left(\overline{M}\right)^\frac{1}{p}\left\|\psi\right\|_{L^q\left(\overline{M},\mu\right)}\\
    &=
    \epsilon\mu\left(\overline{M}\right)^\frac{1}{p}\left\|\psi\right\|_{L^q\left(\overline{M},\mu\right)}.
    \end{align}

    Since $\epsilon>0$ is arbitrary, we conclude that for all $\psi \in L^q\left(\overline{M}, \mu\right)$ we have $\int_{\overline{M}}
            \left(
                g
                -
                h
            \right)\psi
            d\mu = 0$, and so we conclude $g = h$ $\mu$-a.e.
    
    \vspace{10mm}

    However this weak convergence contradicts
    \begin{equation}
        \inf_{j \in \mathbb{N}}
            \left|
                \int_{\overline{M}}
                g_{\overline{n}_j}\psi d\mu
            -
                \int_{\overline{M}}
                g
                \psi d\mu
            \right|
    \geq
    \inf_{j \in \mathbb{N}}
            \left|
                \int_{\overline{M}}
                g_{n_j}\psi d\mu
            -
                \int_{\overline{M}}
                g
                \psi d\mu
            \right|
        >
            \epsilon,
    \end{equation}
    This completes the proof.
\end{proof}

\begin{lemma}
    \label{lem:ptWiseConvergence}
    Suppose $f \in C^\infty(M)$ and $x \in M$. Then we have
    \begin{equation}
        \lim_{r \rightarrow 0}
        \Delta_r f(x)
    =
        C_m \Delta_g f(x)
    \end{equation}
    where $C_m = \frac{1}{2(m+2)}$. The convergence holds $\mu$-a.e. in $\overline{M}$.
\end{lemma}
\begin{proof}
    For fixed $x \in M$, since $S$ is closed in $\overline{M}$ and $x \notin S$,
one has $\dist_g(x,S) > 0$. Hence for all sufficiently small $r$, $B_r(x) = B_r^M (x)$ and the usual
normal-coordinate AMV expansion on the smooth manifold applies. See \cite[Proposition 5.3]{DiasTewodroseAMV} for the proof. The convergence holds $\mu$-a.e. in $\overline{M}$ since $\mu\left(\overline{M} \setminus M\right) = 0$.
\end{proof}

\begin{proof}[Proof of point 1 of Theorem \ref{thm:convergenceOfAMV}]
    First, suppose that $f:\overline{M}\rightarrow\mathbb{R}$ is Lipschitz and belongs to $W^{2,p}(M,g)$ and that $m-k-p\max\left\{1,\frac{\alpha}{2}\right\}>0$.
    By Theorem~\ref{thm:hardyLittleWoodThm} and Proposition \ref{prop:boundOnAMVLaplacian} there exist $r_0,C_\Delta>0$ such that for all $0<r \leq r_0$ \eqref{eq:HardyLittlewoodIneq} holds and
    \begin{equation}
        \left|
            \Delta_r f(x)
        \right|
    \leq
        C_\Delta\left(
            \frac{L_1 C_{\alpha,\dist_g}^{
            \max\left\{
            1,
            \frac{\alpha}{2}
            \right\}
            }}{\overline{\dist}_g(x,S)^{
            \max\left\{
            1,
            \frac{\alpha}{2}
            \right\}      
            }}
        +
            \frac{L_1}{C_{\alpha,\dist_g}^\frac{\alpha}{2} \overline{\dist}_g(x,S)^\frac{\alpha}{2}}
            +
            M_{r_0}\left(\
            \left\|
                Hf
            \right\|_g
            \right)(x)
        \right)
        =:h(x).
    \end{equation}
    Since $m-k-p\max\left\{1,\frac{\alpha}{2}\right\}>0$, by Lemma \ref{lem:1OverDistIsIntegrable} we have that $\overline{\dist}_g(\cdot, S)^{-\max\left\{1,\frac{\alpha}{2}\right\}}, \overline{\dist}_g(\cdot, S)^{-\frac{\alpha}{2}} \in L^p(\overline{M}, \mu)$. On the other hand, since $f\in W^{2,p}(M,g)$, $\|Hf\|_g \in L^p(\overline{M},\mu)$, and so by \eqref{eq:HardyLittlewoodIneq} we have that $M_{r_0}\left(\
            \left\|
                Hf
            \right\|_g
            \right) \in L^p(\overline{M},\mu)$. From this we have $h \in L^p(\overline{M},\mu)$. Thus using Lemma \ref{lem:ptWiseConvergence} and the dominated convergence theorem, we conclude
    \begin{equation}
        \Delta_r f(x)
    \rightarrow_{L^p(\overline{M},\mu)}
        C_m \Delta_g f(x).
    \end{equation}

    If $\alpha < 2$ and $m-k-p = 0$, then define the functions
    \begin{align}
        h_r(x)
    :=
        C_{\Delta}
        \left(
            \frac
            {L_1}
            {r}
            \chi_{T_{C_{\alpha,\dist_g} r}(S)}(x)
        +
            \frac{L_1}
            {\overline{\dist}_g(x ,S)^{\frac{\alpha}{2}}}
        +
            M_{r_0}\left(
                \left\|
                    Hf
                \right\|_g
            \right)
            (x)
        \right).
    \end{align}
    For $r \leq r_0$ we have that $\left|\Delta_r f\right| \leq h_r$. Again, by Lemma~\ref{lem:1OverDistIsIntegrable}, since $m-k-\frac{\alpha}{2}p > 0$, $\overline{\dist}_g(\cdot, S)^{-\frac{\alpha}{2}} \in L^p(\overline{M},\mu)$, and by \eqref{eq:HardyLittlewoodIneq} also $M_{r_0}\left(
                \left\|
                    Hf
                \right\|_g
            \right) \in L^p(\overline{M},\mu)$. Furthermore, since $S$ has finite upper $(m,k)$-dimensional Minkowski content, we have
        \begin{align}
            \limsup_{r \rightarrow 0}
            \int_{\overline{M}}
            \left(
            \frac
            {1}
            {r}
            \chi_{T_{C_{\alpha,\dist_g} r}(S)}(x)
            \right)^p
            d\mu(x)
        =&
            \limsup_{r \rightarrow 0}
            \frac{1}{r^p}
            \mu(T_{C_{\alpha,\dist_g} r}(S))\\
        =&
            \limsup_{r \rightarrow 0}
            \frac{1}{r^{m-k}}
            \mu(T_{C_{\alpha,\dist_g} r}(S))
        <
            \infty.
        \end{align}
        This implies $\limsup_{r \rightarrow 0} \left\| \Delta_r f \right\|_{L^p(\overline{M},\mu)}\leq \limsup_{r \rightarrow 0} \left\| h_r \right\|_{L^p(\overline{M},\mu)}< \infty$. Applying Lemmas~\ref{lem:ptWiseConvergence} and \ref{lem:weakConvergenceLemma}, we obtain
        \begin{equation}
            \Delta_r f
        \rightharpoonup_{L^p(\overline{M},\mu)}
            C_m \Delta_g f.
        \end{equation}
    \end{proof}

    \begin{proof}[Proof of point 2 of Theorem \ref{thm:convergenceOfAMV}]
    Suppose $f:\overline{M}\rightarrow\mathbb{R}$ is $\beta$-H\"older continuous, and let $q\geq 1$ such that $m-k-\frac{\alpha pq}{2(p-q)}>0$ and $m-k-q\max\left\{\frac{(2-\beta)\alpha}{2},(2-\beta)\right\}>0$. Define the function
    \begin{align}
        \tilde{h}(x)
    :=
    C_\Delta
    \left(
       \frac{L_\beta C_{\alpha,\dist_g}^{
            \max\left\{
            \frac{(2-\beta)\alpha}{2},
            (2-\beta)
            \right\}}}{\overline{\dist}_g(x,S)^{
            \max\left\{
            \frac{(2-\beta)\alpha}{2},
            (2-\beta)
            \right\}}}
    +
        \frac{1}{C_{\alpha,\dist_g}^{\frac{\alpha}{2}}}
            \frac{1}{\overline{\dist}_g(x,S)^{\frac{\alpha}{2}}}
        M_{r_0}\left(
        \left\|
            \nabla f
        \right\|_g
        \right)
        (x)
        +
        M_{r_0}\left(
            \left\|Hf\right\|_g
        \right)
        (x)
    \right).
    \end{align}

    We know by Proposition \ref{prop:boundOnAMVLaplacian} that for $0< r \leq r_0$ we have $\left| \Delta_r f(x)\right| \leq \tilde{h}(x)$ for all $x \in M$. Choose $\gamma,\sigma>0$ such that $\frac{1}{\sigma} + \frac{1}{\gamma} = 1$ and $q\sigma = p$, so $\gamma= \frac{p}{p-q}$. Thus we have
    \begin{align}
        &
        \int_{\overline{M}}
        \left(
            \frac{1}{\overline{\dist}_g(x,S)^{\frac{\alpha}{2}}}
        M_{r_0}\left(
        \left\|
            \nabla f
        \right\|_g
        \right)
        (x)
        \right)^q
        d\mu(x)\\
    \leq&
        \left\|
            \frac{1}{\overline{\dist}_g(\cdot,S)^{\frac{q\alpha}{2}}}
        \right\|_{L^\gamma(\overline{M}, \mu)}
        \left\|
            M_{r_0}\left(
        \left\|
            \nabla f
        \right\|_g
        \right)^q
        \right\|_{L^\sigma(\overline{M}, \mu)}.
        \label{eq:holderArgumentSobolev}
    \end{align}
    By Lemma \ref{lem:1OverDistIsIntegrable} we know that $\left\|{\overline{\dist}_g(\cdot,S)^{-\frac{q\alpha}{2}}}
        \right\|_{L^\gamma(\overline{M}, \mu)}^\frac{1}{q} < \infty$ since
    \begin{equation}
        m - k -\frac{\alpha q\gamma}{2}
    =
        m - k -\frac{\alpha qp}{2(p-q)}
    >
        0.
    \end{equation}
    Furthermore, since $f\in W^{2,p}(M,g)$, applying equation \eqref{eq:HardyLittlewoodIneq} we have \begin{equation}
    \left\|
            M_{r_0}\left(
        \left\|
            \nabla f
        \right\|_g
        \right)^q
        \right\|_{L^\sigma(\overline{M}, \mu)}
        =
        \left\|
            M_{r_0}\left(
        \left\|
            \nabla f
        \right\|_g
        \right)
        \right\|_{L^p(\overline{M}, \mu)}^q
        < \infty,
    \end{equation}
    and so $\frac{M_{r_0}\left(
        \left\|
            \nabla f
        \right\|_g
        \right)
        (\cdot)}{\overline{\dist}_g(\cdot,S)^{\frac{\alpha}{2}}} \in L^q(\overline{M}, \mu)$. We have $\overline{\dist}_g(\cdot ,S)^{-\max\left\{
            \frac{(2-\beta)\alpha}{2},
            (2-\beta)
            \right\}} \in L^q(\overline{M}, \mu)$ since $m-k-q\max\left\{
            \frac{(2-\beta)\alpha}{2},
            (2-\beta)
            \right\} > 0$. Furthermore by Lemma \ref{lem:finiteMeasureLemma}, the measure $\mu$ is finite. Since $q < p$ and $\left\|Hf\right\|_g \in L^p\left(\overline{M},\mu\right)$, by Theorem \ref{thm:hardyLittleWoodThm}, $M_{r_0}\left(\left\|H f\right\|_g\right) \in L^p(\overline{M},\mu)$, and so $M_{r_0}\left(\left\|H f\right\|_g\right) \in L^q(\overline{M},\mu)$. Thus $\tilde{h} \in L^q(\overline{M}, \mu)$, and by applying Lemma \ref{lem:ptWiseConvergence} and dominated convergence we have
        \begin{equation}
            \Delta_r f \rightarrow_{L^q(\overline{M}, \mu)}
            C_m \Delta_g f.
        \end{equation}

        On the other hand, suppose now that $m-k-q\max\left\{
            \frac{(2-\beta)\alpha}{2},
            (2-\beta)
            \right\}=0$ and $m-k -\frac{\alpha pq}{2(p-q)}>0$. Define the functions given by
            \begin{equation}
                \tilde{h}_{r}(x)
            :=
             C_{\Delta}
        \left(
            \frac
            {L_\beta}
            {r^{2-\beta}}
            \chi_{T_{C_{\alpha,\dist_g} r^{\min\{1,2/\alpha\}}}(S)}(x)
        +
        \frac{M_{r_0}\left(\left\|\nabla f\right\|_g\right)(x)}
            {\overline{\dist}_g(x ,S)^{\frac{\alpha}{2}}}
        +
            M_{r_0}\left(
                \left\|
                    Hf
                \right\|_g
            \right)
            (x)
        \right).
            \end{equation}
        By Proposition \ref{prop:boundOnAMVLaplacian}, for $r \leq r_0$ we have $\left|\Delta_r f(x)\right| \leq \tilde{h}_r(x)$ for $x \in M$. By the previous arguments we know that $\frac{M_{r_0}\left(
        \left\|
            \nabla f
        \right\|_g
        \right)
        (\cdot)}{\overline{\dist}_g(\cdot,S)^{\frac{\alpha}{2}}} \in L^q\left(\overline{M}, \mu\right)$ and $M_{r_0}\left(\left\|\nabla f\right\|_g\right) \in L^q(\overline{M},\mu)$. Also, since $1<q\leq p$, $\left\|Hf\right\|_g \in L^p\left(\overline{M},\mu\right)$, and $\mu$ is finite, $\left\|Hf\right\|_g \in L^q\left(\overline{M},\mu\right)$, and Theorem~\ref{thm:hardyLittleWoodThm} gives $M_{r_0}\left(\left\|Hf\right\|_g\right) \in L^q\left(\overline{M}, \mu\right)$. Furthermore, since $S$ has finite upper $(m,k)$-dimensional Minkowski content, there exists $C_{\mathcal{M}}>0$ such that
        \begin{align}
            &\limsup_{r \rightarrow 0}
            \int_{\overline{M}}
            \left(
                \frac
                {1}
                {r^{2-\beta}}
                \chi_{T_{C_{\alpha,\dist_g} r^{\min\{1,2/\alpha\}}}(S)}(x)
            \right)^q
            d\mu(x)\\
        =&
            \limsup_{r \rightarrow 0}
            \frac{\mu(T_{C_{\alpha,\dist_g} r^{\min\left\{1,2/\alpha\right\}}}(S))}
            {r^{(2-\beta)q}}\\
        \leq&
            \limsup_{r \rightarrow 0}
            \frac{\mu(T_{C_{\alpha,\dist_g} r^{\min\left\{1,2/\alpha\right\}}}(S))}
            {r^{(2-\beta)q}}\\
        \leq&
            \limsup_{r \rightarrow 0}
            C_{\mathcal{M}}
            \frac{\left(C_{\alpha,\dist_g} r^{\min\{1,2/\alpha\}}\right)^{m-k}}{r^{(2-\beta)q}}
        <
            \infty,
        \end{align}
        Since $m-k-q\max\left\{
            \frac{(2-\beta)\alpha}{2},
            (2-\beta)
            \right\}=0$, it follows that $\limsup_{r \rightarrow 0}\|\Delta_r f\|_{L^q(\overline{M}, \mu)} \leq \limsup_{r \rightarrow 0}\|\Tilde{h}_r\|_{L^q(\overline{M}, \mu)} < \infty$.  Thus applying Lemmas \ref{lem:ptWiseConvergence} and \ref{lem:weakConvergenceLemma} we have for any sequence $r_n \rightarrow 0$ that
            \begin{equation}
                \Delta_{r_n} f \rightharpoonup_{L^q(\overline{M}, \mu)} C_m \Delta_g f,
            \end{equation}
            Thus, by the sequential characterization of limits, we obtain $\Delta_{r} f \rightharpoonup_{L^q(\overline{M}, \mu)} C_m \Delta_g f$ as $r \rightarrow 0$. Finally, if we assume that $(M,g)$ satisfies \eqref{eq:SobolevIneq} with exponent $p$, by Lemma \ref{lem:sobolevIneq2derivatives} we have
        \begin{equation}
        \label{eq:SobolevIneqOnProof}
            \left\|
            \left\|
            \nabla f
            \right\|_g
            \right\|_{L^{p^*}(\overline{M}, \mu)}
        \leq
        C_{\mathcal{S},p}
        \left(
            \left\|
            \left\|
             Hf
            \right\|_g
            \right\|_{L^{p}(\overline{M}, \mu)}
            +
            \left\|
            \left\|
             \nabla f
            \right\|_g
            \right\|_{L^{p}(\overline{M}, \mu)}
        \right)
        <\infty.
        \end{equation}
        Using \eqref{eq:HardyLittlewoodIneq}, we have $M_{r_0}\left(\left\|\nabla f\right\|_g\right) \in L^{p^*}(\overline{M},\mu)$. We can thus apply the same argument as in \eqref{eq:holderArgumentSobolev} with $\sigma = \frac{p^*}{q}$ and $\gamma = \frac{p^*}{p^*-q}$, obtaining the same conclusion as before with $p$ replaced by $p^*$ in the condition involving $\frac{\alpha qp}{2(p-q)}$ for $1< q \leq p$.
        \end{proof}

        \section{Examples of completions of manifolds}

        In this section we study the tuples $(N,S,g,\{\psi_i\}_{i=1}^l)$ that satisfy definition \ref{def:riemannianManifoldWithSingularitySet} and prove Theorem \ref{thm:riemannianManifoldWithSingularityThm}, that is we show that these spaces will satisfy the conditions of Theorem \ref{thm:convergenceOfAMV}. This will give us a large family of spaces satisfying the conditions of Theorem \ref{thm:convergenceOfAMV}.
        
        \subsection{Topology of completion of $(N\setminus S, g)$}

        In this subsection we start by exploring several notions of dimension of a set from geometric measure theory and the relations among them. These notions are the Minkowski content and dimension already defined in \ref{def:minkowskyContent} and the Hausdorff measure and dimension. For a reference on the study of the Hausdorff measure and Minkowski content in Euclidean space, we suggest \cite{AFP}.

        \begin{definition}
\label{def:hausdorffMeasureAndDimension}
    Given $A\subset \mathbb{R}^m$, the $k$-dimensional Hausdorff measure is defined by
    \begin{equation}
        \mathcal{H}^k(A) 
    = 
        \lim_{\delta \rightarrow 0}
        \mathcal{H}^k_\delta(A),
    \end{equation}
    where $\mathcal{H}^k_\delta$ is defined by
    \begin{equation}
        \mathcal{H}^k_\delta(A)
    :=
        \frac{1}{2^k}
        \inf
        \left\{
            \sum_{i \in I}
            r_i^k
            :
            r_i \leq \delta,
            A \subset \cup_{i \in I}\mathbb{B}^m_{r_i}(x_i)
        \right\}.
    \end{equation}
    We define the Hausdorff dimension of $A$ as
    \begin{equation}
        \text{dim}_{\mathcal{H}}(A)
    :=
        \inf\{k \geq 0: \mathcal{H}^k(A) = 0\}.
    \end{equation}
    \end{definition}

    \begin{lemma}
    \label{lem:dimensionBoundHausdorffMinkowsky}
    Let $A\subset \mathbb{R}^m$ be a bounded set. Then, for every $0\leq k \leq m$, in the Euclidean metric measure space $\mathfrak{E} = (\mathbb{R}^m, \dist_\euc, \mathcal{L}^m)$, we have
    \begin{equation}
        \mathcal{H}^k(A)
    \leq
        \frac{1}{\mathcal{L}^m(\mathbb{B}^m_1(0))}\mathcal{M}^{(m,k,-)}_{\mathfrak{E}}(A).
    \end{equation}
    In particular, $\text{dim}_{\mathcal{H}}(A) \leq \text{dim}_{\mathcal{M},\mathcal{E}}^{(m,-)}(A)$.
    \end{lemma}
    \begin{proof}
        Given $r > 0$, choose a maximal pairwise disjoint collection of balls $\{\mathbb{B}^m_{r}(x_i)\}_{i=1}^{N_r}$, $x_i \in A$, with centers in $A$. We have that
    \begin{equation}
    \label{eq:inclusionOfAinBalls}
        A \subset \cup_{i=1}^{N_r}\mathbb{B}^m_{2r}(x_i),
    \end{equation}
    since if there exists $y \in A \setminus \cup_{i=1}^{N_r}\mathbb{B}^m_{2r}(x_i)$, then by triangle inequality for all $z \in \mathbb{B}_r^m(y)$ and $x_i$ we have
    \begin{equation}
        \|
            z - x_i
        \|
    \geq
    -
        \|
            z-y
        \|
    +
    \|
        y-x_i
    \|
    >
        -r+2r
    = r,
    \end{equation}
    showing that $\mathbb{B}^m_r(y)$ is disjoint from the collection, contradicting maximality. It is clear that
    \begin{equation}
        \cup_{i=1}^{N_r}\mathbb{B}^m_r(x_i)
    \subset
        T_r(A),
    \end{equation}
    so by disjointness of the balls we have
    \begin{equation}
    \label{eq:lowerBoundMinkowski}
    \sum_{i=1}^{N_r}
        \mathcal{L}^m(\mathbb{B}^m_1(0))
        r^m
    =
    \mathcal{L}^m
    \left(
        \cup_{i=1}^{N_r}\mathbb{B}^m_r(x_i)
    \right)
    \leq
    \mathcal{L}^m\left(
        T_r\left(A\right)
    \right).
    \end{equation}
    Combining \eqref{eq:lowerBoundMinkowski} and inclusion \eqref{eq:inclusionOfAinBalls}, we obtain
    \begin{equation}
        \mathcal{L}^m(\mathbb{B}^m_1(0))\mathcal{H}^k_{2r}(A)
    \leq
        \frac{\mathcal{L}^m(\mathbb{B}^m_1(0))}{2^k}\sum_{i=1}^{N_r}
        (2r)^{k}
    \leq
       \frac{\mathcal{L}^m(T_r(A))}{r^{m-k}}.
    \end{equation}
    Taking the $\liminf$ with $r \rightarrow 0$ yields
    \begin{equation}
        \mathcal{H}^k(A)
    \leq
        \frac{1}
        {\mathcal{L}^m(\mathbb{B}^m_1(0))}
        \mathcal{M}^{(m,k,-)}_{\mathfrak{E}}(A).
    \end{equation}
    \end{proof}
    Let $x \in \mathbb{R}^m$ and a unit vector $v \in \partial \mathbb{B}^m_1(0)$ let
    \begin{equation}
        \pi_{x,v}
    :=
        \left\{
            y \in \mathbb{R}^m:
            \left\langle
                y-x,
                v
            \right\rangle
            =
            0
        \right\}
    \end{equation}
    be the hyperplane orthogonal to $v$ passing through the point $x$. Consider the projection map to the hyperplane $P_{x,v}:\mathbb{R}^m\rightarrow \pi_{x,v}$, defined by
    \begin{align}
        P_{x,v}(w)
    =
        w - \langle w, v \rangle v.
    \end{align}
    \begin{lemma}
    \label{lem:fullMeasureLemma}
    Suppose $\tilde{S} \subset \mathbb{R}^m$ is a set that satisfies $\mathcal{H}^{m-1}(\tilde{S}) = 0$. Given $x \in \mathbb{R}^m$, then the image
    \begin{align}
        P_{x,v}\left(\tilde{S}\right)
    \end{align}
    has zero $\mathcal{H}^{m-1}$ measure.
\end{lemma}
\begin{proof}
     Using Proposition 2.4.9 of \cite{AFP}, since $P_{x,v}$ is a $1$-Lipschitz maps, we have that
    \begin{equation}
        \mathcal{H}^{m-1}(P_{x,v}(\tilde{S}))
        =
        0.
    \end{equation}
\end{proof}

\begin{lemma}
    \label{lem:invertibleProjections}
    Let $x,y \in \mathbb{R}^m$ with $x \neq y$ and $\epsilon > 0$, and let $v_{x,y} 
    := \frac{x-y}{\|x-y\|}$. Define the sets
    \begin{align}
        A_{x,v_{x,y},\epsilon}
    :=&
        \left\{
            z \in \partial \mathbb{B}^{m}_{\epsilon}(x):
            \left\langle
                z-x,v_{x,y}
            \right\rangle
            > 
            0
        \right\},\\
        A_{y,v_{x,y},\epsilon}
    :=&
        \left\{
            z \in \partial \mathbb{B}^{m}_{\epsilon}(y):
            \left\langle
                z-y,v_{x,y}
            \right\rangle
            < 
            0
        \right\}.
    \end{align}
    Then the maps given by \begin{align}\tilde{P}_{x,v_{x,y},\epsilon} &:= \left(\left.P_{x,v_{x,y}}\right|_{A_{x,v_{x,y},\epsilon}}\right)^{-1}-x:
    \mathbb{B}^m_\epsilon(x) \cap \pi_{x,v_{x,y}}
    \rightarrow
    \partial \mathbb{B}^m_\epsilon(0)
    \\
    \tilde{P}_{y,v_{x,y},\epsilon} &:= \left(\left.P_{x,v_{x,y}}\right|_{A_{y,v_{x,y},\epsilon}}\right)^{-1}-y
    :
    \mathbb{B}^m_\epsilon(x) \cap \pi_{x,v_{x,y}}
    \rightarrow
    \partial \mathbb{B}^m_\epsilon(0)
    \end{align}
    are well defined.
\end{lemma}
\begin{proof}
    First we show that the map $\left.P_{x,v_{x,y}}\right|_{A_{x,v_{x,y},\epsilon}}$ is injective. Given $z_1, z_2 \in A_{x,v_{x,y},\epsilon}$, if $P_{x,v_{x,y}}(z_1) = P_{x,v_{x,y}}(z_2) = w \in \pi_{x,v_{x,y}}$, then
    \begin{align}
        z_1 &= w + t_1v_{x,y}\\
        z_2 &= w + t_2v_{x,y}.
    \end{align}
    Since $z_1,z_2 \in A_{x,v_{x,y},\epsilon}$, it follows that $t_1, t_2 > 0$ and
    \begin{align}
        \|x-w\|^2
        +
        t_1^2 =& \epsilon^2.\\
        \|x-w\|^2
        +
        t_2^2 =& \epsilon^2.
    \end{align}
    Because $t_1, t_2>0$, then $t_1 = t_2$.

    Now we show that the map $P_{x,v_{x,y}}$ is also surjective onto $\mathbb{B}^m_\epsilon(x) \cap \pi_{x,v_{x,y}}$. Given $w \in \mathbb{B}^m_\epsilon(x) \cap \pi_{x,v_{x,y}}$, there exists $t = \sqrt{\epsilon^2 - \|x-w\|^2}>0$ such that $z = w + tv_{x,y} \in A_{x,v_{x,y},\epsilon}$. Thus
    \begin{equation}
        P_{x,v_{x,y}}(w+tv_{x,y})
    =
        w,
    \end{equation}
    showing the surjectivity. This shows that the map $\left.P_{x,v_{x,y}}\right|_{A_{x,v_{x,y},\epsilon}}$ is invertible.

    A similar applies to $\left.P_{x,v_{x,y}}\right|_{A_{y,v_{x,y},\epsilon}}$, showing it is injective and surjective onto $\mathbb{B}^m_\epsilon(x) \cap \pi_{x,v_{x,y}}$ by using that $v_{x,y} = \frac{x-y}{\|x-y\|}$.
\end{proof}

\begin{lemma}
\label{lem:euclideanMetricByCurves}
    Let $R>0$ and $\tilde{S} \subset \mathbb{R}^m$ be a closed subset such that $\mathcal{H}^{m-1}\left(\tilde{S}\right) = 0$. Given the metric on $\mathbb{B}^m_R(0)\setminus \tilde{S}$ defined by
    \begin{equation}
        \tilde{\dist_g}_{\euc,\tilde{S}}(x,y)
    =
        \inf
        \left\{
            \int_0^1
            \|
                \gamma'(t)
            \|dt:
            \gamma(0) = x,
            \gamma(1) = y,
                \gamma:[0,1] \rightarrow \mathbb{B}^m_R(0)\setminus \tilde{S}
                \text{ is piecewise } C^1
        \right\},
    \end{equation}
    then
    \begin{equation}
        \tilde{\dist_g}_{\euc,\tilde{S}}(x,y) = \|x-y\| = \dist_\euc(x,y).
    \end{equation}
\end{lemma}
\begin{proof}
First we handle the case $R=1$. Given $x,y \in \mathbb{B}^m_1(0) \setminus \tilde{S}$ with $x \neq y$ define $v_{x,y} = \frac{x-y}{\|x-y\|}$. Let $\pi_{x,v_{x,y}}$ be the $m-1$ dimensional hyperplane containing $x$ orthogonal to $v_{x,y}$. Since $\tilde{S}$ is closed, consider $\epsilon>0$ such that $\overline{\mathbb{B}^m_\epsilon(x)}, \overline{\mathbb{B}^m_\epsilon(y)} \subset \mathbb{B}^m_1(0) \setminus \tilde{S}$.

    Now define the set given by
    \begin{align}
        V_{x,y,\epsilon}
    :=&
        \left\{
            w \in \mathbb{B}^m_\epsilon(x) \cap \pi_{x,v_{x,y}}:
            \forall_{t \in \mathbb{R}}
            w+tv_{x,y} \notin \tilde{S}
        \right\}.
    \end{align}
    We have that the set 
    \begin{equation}
        V_{x,y,\epsilon}
        =
        \mathbb{B}^m_\epsilon(x) \cap \pi_{x,v_{x,y}} \setminus P_{x,v_{x,y}}(\tilde{S}),
    \end{equation}
    has full measure $\mathcal{H}^{m-1}(\mathbb{B}^m_\epsilon(x)\cap\pi_{x,v_{x,y}})$ by Lemma \ref{lem:fullMeasureLemma}, thus $V_{x,y,\epsilon} \neq \emptyset$.
    Taking $w \in V_{x,y,\epsilon}$ and defining the curve $\gamma_{\epsilon}:[0,1] \rightarrow \mathbb{B}^m_1(0)$ by 
    \begin{equation}
        \gamma_\epsilon(t)
    =
        \begin{cases}
          3\left(\frac{1}{3}-t\right)x + 3t\left(
          \tilde{P}_{x,v_{x,y}, \epsilon}
            (w)
            +
            x
            \right)& t \in  [0,\frac{1}{3}]\\
          3(\frac{2}{3}-t)\left(\tilde{P}_{x,v_{x,y}, \epsilon}
            (w)
            +
            x
            \right)+ 3(t-\frac{1}{3})
            \left(\tilde{P}_{y,v_{x,y}, \epsilon}
            (w)
            +
            y
            \right) & t \in [\frac{1}{3},\frac{2}{3}]\\
          3(t-\frac{2}{3})y + 3(1-t)\left(\tilde{P}_{y,v_{x,y}, \epsilon}
            (w)
            +
            y
            \right)& t \in [\frac{2}{3},1]
        \end{cases}
    \end{equation}
    is such that $\gamma_\epsilon([0,1]) \subset \mathbb{B}^m_1(0) \setminus \tilde{S}$ and satisfies
    \begin{equation}
        l_\euc(\gamma_\epsilon)
        \leq
        4\epsilon + \|x-y\|,
    \end{equation}
    since $\epsilon>0$ is arbitrary we conclude that
    \begin{equation}
        \|x-y\|
    \leq
        \tilde{\dist_g}_{\euc,\tilde{S}}(x,y)
    \leq
        \|x-y\|.
    \end{equation}

    Since the result is scale invariant, given $R>0$, the result is also true for any ball $\mathbb{B}^m_R(0)$.

\end{proof}

\begin{prop}
\label{prop:completionMetric}
    Suppose the tuple $(N,S,g,\{\psi_i\}_{i=1}^l)$ satisfies Definition \ref{def:riemannianManifoldWithSingularitySet} with dimension $m \in \mathbb{N}$ and that for all parametrizations $\mathcal{H}^{m-1}\left(\psi_i^{-1}(S)\right) = 0$. Then the completion of $(N\setminus S, \dist_g)$ is homeomorphic to $N$ with the topology given by the parametrizations $\{\psi_i\}_{i=1}^l$.
\end{prop}
\begin{proof}
To prove the result, we first prove the following equivalence
\begin{align}
\label{eq:equivalenceOfCauhySequences}
    &x_n \in N \setminus S \text{ is a Cauchy sequence for distance } \dist_g\\
\Leftrightarrow\quad
    &
    \exists_{i \in \{1,\cdots,l\},n_0
    \in \mathbb{N}} : \forall_{ n\geq n_0} 
    \quad
    x_n \in \psi_i\left(\mathbb{B}^m_{\frac{3}{2}}(0)\right) \text{ and }\\ 
    &\quad\quad
    \left\{\psi_i^{-1}(x_n)\right\}_{n \geq n_0}
    \text{ is a Cauchy sequence for the Euclidean metric}.
\end{align}

Suppose first that $x_n \in N \setminus S$ is a Cauchy sequence for the metric $\dist_g$. Consider $c_{LB} > 0$ from equation \eqref{eq:boundBellowMetric}. There exists $n_0 \in \mathbb{N}$ such that for all $n,m \geq n_0$ we have
    \begin{equation}
        \label{eq:CauchySmallForParametrization}
        \forall n,m \geq n_0 \quad
        \dist_g(x_n,x_m) < \frac{\sqrt{c_{LB}}}{4}.
    \end{equation}
    By equation \eqref{eq:NcoveredFiniteParametrizations}, there exists $i \in \{1,\cdots,l\}$ such that $x_{n_0} \in \psi_i(\mathbb{B}^m_1(0))$. Thus by  equation \eqref{eq:boundBellowMetric} and \eqref{eq:CauchySmallForParametrization} we conclude that  
\begin{equation}
\label{eq:pointsInParametrization}
\forall n \geq n_0 \quad\quad 
x_n \in \psi_i\left(\mathbb{B}^m_{\frac{5}{4}}(0)\right).
\end{equation}
Furthermore by equations \eqref{eq:boundBellowMetric},  \eqref{eq:CauchySmallForParametrization} and \eqref{eq:pointsInParametrization} if for $n,m \geq n_0$ the curve $\gamma_{n,m}:[0, 1] \rightarrow N \setminus S$ is such that
\begin{align}
\label{eq:smallLengthForCurve}
    \int_0^1 g(\gamma_{n,m}'(t),\gamma_{n,m}'(t))^\frac{1}{2}dt
    <
    \frac{\sqrt{c_{LB}}}{4},\quad
    \gamma_{n,m}(0) = x_n,\quad
    \gamma_{n,m}(1) = x_m,
\end{align}
then $\gamma_{n,m}\left(\left[0,1\right]\right) \subset \psi_i\left(\mathbb{B}^m_{\frac{3}{2}}\left(0\right)\right)$. Because $\dist_g(x_n,x_m)< \frac{\sqrt{c_{LB}}}{4}$, there are curves with length below $\frac{\sqrt{c_{LB}}}{4}$ arbitrarily close to the infimum, and since for $n \geq n_0$ we have $x_n \in \psi_i\left(\mathbb{B}^m_{\frac{5}{4}}(0)\right)$, every such curve remains in $\psi_i\left(\mathbb{B}^m_{\frac{3}{2}}(0)\right)$. This shows that for $n,m \geq n_0$
\begin{align}
    \dist_g(x_n,x_m)
=&
    \inf
    \Bigg\{
        \int_0^1g\left(\left(d\psi_i\right)_{\tilde{\gamma}(t)}\tilde{\gamma}'(t), \left(d\psi_i\right)_{\tilde{\gamma}(t)}\tilde{\gamma}'(t)\right)^\frac{1}{2}
        dt:
        \psi_i(\tilde{\gamma}(0))=x_n,\\
    &\quad
        \psi_i(\tilde{\gamma}(1))=x_m,
        \text{ and }\tilde{\gamma}:[0,1] \rightarrow \mathbb{B}^m_{\frac{3}{2}}(0)\setminus \psi_i^{-1}(S) 
    \Bigg\}\\
\geq&
        \sqrt{c_{LB}}
        \inf
    \Big\{
        l_\euc(\tilde{\gamma}):
        \psi_i\left(\tilde{\gamma}\left(0\right)\right)=x_n, \psi_i\left(\tilde{\gamma}\left(1\right)\right)=x_m,\\
    &\quad
        \quad\quad\quad\quad\text{ and }\tilde{\gamma}:[0,1] \rightarrow \mathbb{B}^m_{\frac{3}{2}}(0)\setminus \psi_i^{-1}(S) 
    \Big\}.
\end{align}
Since $\overline{\mathbb{B}^m_{\frac{3}{2}}(0)} \cap \psi_i^{-1}\left(S\right)$ is closed and $\mathcal{H}^{m-1}\left(\overline{\mathbb{B}^m_{\frac{3}{2}}(0)} \cap \psi_i^{-1}\left(S\right)\right) = 0$ we can apply Lemma \ref{lem:euclideanMetricByCurves} to conclude that for $n,m \geq n_0$
\begin{equation}
    \dist_g(x_n,x_m)
\geq
    \sqrt{c_{LB}}
    \left\|
        \psi_i^{-1}(x_n)
        -
        \psi_i^{-1}(x_m)
    \right\|.
\end{equation}
This shows that $\left\{\psi_i^{-1}(x_n)\right\}_{n}$ is a Cauchy sequence for the Euclidean metric since $x_n$ is a Cauchy sequence for $\dist_g$, proving the first implication.

On the other hand suppose $x_n \in N \setminus S$ is a sequence such that there exist $i \in \{1,\cdots, l\}$ and $n_0 \in \mathbb{N}$ such that for all $n \geq n_0$ then $x_n \in \psi_i\left(\mathbb{B}^m_{\frac{3}{2}}(0)\right)$ and $\{\psi_i^{-1}(x_n)\}_{n \geq n_0}$ is Cauchy for the Euclidean metric. Since $\overline{\mathbb{B}^m_{\frac{3}{2}}(0)} \cap \psi_i^{-1}\left(S\right)$ is closed and $\mathcal{H}^{m-1}\left(\overline{\mathbb{B}^m_{\frac{3}{2}}(0)} \cap \psi_i^{-1}\left(S\right)\right) = 0$ we can apply Lemma \ref{lem:euclideanMetricByCurves} and inequality \eqref{eq:boundAboveMetric} to conclude for $n,m \geq n_0$
\begin{align}
    &
    \|\psi_i^{-1}(x_n)-\psi_i^{-1}(x_m)\|\\
    =
    &
    \inf
    \Big\{
       l_{\text{euc}}(\tilde{\gamma}):
        \psi_i(\tilde{\gamma}(0))=x_n, \psi_i(\tilde{\gamma}(1))=x_m
        \text{ and }\tilde{\gamma}:[0,1] \rightarrow \mathbb{B}^m_{\frac{3}{2}}(0)\setminus \psi_i^{-1}(S)
    \Big\}\\
    \geq&
    \frac{1}{\sqrt{C_{UB}}}
    \inf
    \Bigg\{
       \int_0^1g\left(\left(d\psi_i\right)_{\tilde{\gamma}(t)}\tilde{\gamma}'(t), \left(d\psi_i\right)_{\tilde{\gamma}(t)}\tilde{\gamma}'(t)\right)^\frac{1}{2}
        dt:
        \psi_i(\tilde{\gamma}(0))=x_n,\\
    &\quad\quad\quad\quad\quad\quad\quad\quad
        \psi_i(\tilde{\gamma}(1))=x_m
        \text{ and }\tilde{\gamma}:[0,1] \rightarrow \mathbb{B}^m_{\frac{3}{2}}(0)\setminus \psi_i^{-1}(S)
    \Bigg\}\\
    \geq&
        \frac{1}{\sqrt{C_{UB}}}
        \dist_g(x_n, x_m),
        \label{eq:boundOfMetric}
\end{align}
thus $x_n$ is a Cauchy sequence for $\dist_g$.

This allows us to define a map given by $\Phi:\left(\overline{N\setminus S}, \overline{\dist}_g\right) \rightarrow N$, where given an equivalence class of Cauchy sequences $[\{x_n\}_{n \in \mathbb{N}}] \subset N \setminus S$, we send it to the limit $\Phi([\{x_n\}_{n \in \mathbb{N}}]) = \psi_i\left(\lim_n \psi_i^{-1}(x_n)\right) \in N$.

The limit exists since $\psi_{i}^{-1}(x_n) \subset \mathbb{B}^m_{\frac{3}{2}}(0)$, $\psi_i:\mathbb{B}^m_2(0)\rightarrow N$ is continuous and the limit $\lim_n \psi_i^{-1}(x_n) \in \overline{\mathbb{B}^m_{\frac{3}{2}}(0)}$ exists because the sequence is Cauchy and $\overline{\mathbb{B}^m_{\frac{3}{2}}(0)}$ is complete.

The map is well defined because if $\{x_n\}_{n \in \mathbb{N}}, \{y_n\}_{n \in \mathbb{N}} \subset [\{x_n\}_{n \in \mathbb{N}}]$ are two Cauchy sequences in the same equivalence class, we can define the sequence
\begin{equation}
    z_n:=
    \begin{cases}
        x_{\frac{n}{2}}
        &
        \text{ if }
        n \text{ is even}\\
        y_{\frac{n+1}{2}}
        & \text{ if } n \text{ is odd}
    \end{cases}.
\end{equation}
Since $\{x_n\}_{n \in \mathbb{N}}, \{y_n\}_{n \in \mathbb{N}} \subset [\{x_n\}_{n \in \mathbb{N}}]$ are two Cauchy sequences in the same equivalence class, the sequence $\left\{z_n\right\}_{n \in \mathbb{N}}$ is also Cauchy and in the same class. Thus, by \eqref{eq:equivalenceOfCauhySequences} there exists $i \in \{1,...,l\}$ and $n_0 \in \mathbb{N}$ such that
\begin{equation}
    \forall_{n \geq n_0}
\quad
    z_n \in \psi_{i}\left(\mathbb{B}^m_{\frac{3}{2}}(0)\right),
    \quad
    \left\{
        \psi_i^{-1}(z_n)
    \right\}_{n \geq n_0} \text{ is Cauchy for the Euclidean metric}.
\end{equation}
Thus we have that $z_n$ converges in the topology of $N$ to some point $z \in N$. On the other hand by the definition of $\Phi$ we have that $x_n \rightarrow \Phi\left(\left[\left\{x_n\right\}_{n \in \mathbb{N}}\right]\right)$ and $y_n \rightarrow \Phi\left(\left[\left\{y_n\right\}_{n \in \mathbb{N}}\right]\right)$ in the topology of $N$, thus since $N$ is Hausdorff, we conclude by uniqueness of limits
\begin{equation}
    \Phi\left(\left[\left\{x_n\right\}_{n \in \mathbb{N}}\right]\right)
    =
    z
    =
    \Phi\left(\left[\left\{y_n\right\}_{n \in \mathbb{N}}\right]\right),
\end{equation}
showing it is well defined.

Now we show the map is surjective. Given $y \in N$, there exists $i \in \{1,...,l\}$ and $x \in \mathbb{B}^m_1(0)$ such that $\psi_i(x) = y$. Furthermore since $\mathcal{H}^{m-1}\left(\psi_i^{-1}(S) \cap \mathbb{B}^m_{\frac{3}{2}}(0)\right) = 0$, we have that $\mathbb{B}^m_{\frac{3}{2}}(0) \setminus \psi_i^{-1}\left(S\right)$ is dense in $\mathbb{B}^m_\frac{3}{2}(0)$, thus there exists a sequence $x_n \in \mathbb{B}^m_{\frac{3}{2}}(0)\setminus \psi_i^{-1}\left(S\right)$ such that $x_n \rightarrow x$, and so $x_n$ is a Cauchy sequence. By the equivalence \eqref{eq:equivalenceOfCauhySequences}, we have that $\psi_i(x_n) = y_n$ is a Cauchy sequence for the distance $\dist_g$, and $\Phi\left(\left[\left\{y_n\right\}_{n \in \mathbb{N}}\right]\right) = \psi_i\left(\lim_n\psi_i^{-1}(y_n)\right)
=
\psi_i(\lim_n x_n)
= \psi_i\left(x\right) 
= y$, showing surjectivity.

To prove injectivity, suppose that \begin{equation}
\Phi
\left(
\left[
\left\{
x_n
\right\}_{n \in \mathbb{N}}
\right]
\right)
 = 
\Phi\left(
\left[
\left\{y_n
\right\}_{n \in \mathbb{N}}
\right]
\right)
=
q.  
\end{equation}
Consider $i \in \{1,...,l\}$ such that $q \in \psi_i\left(\mathbb{B}^m_1(0)\right)$. Both sequences $x_n$ and $y_n$ converge in the topology of $N$ to the same point $q$, thus there exists some $n_0 \in \mathbb{N}$ such that for all $n \geq n_0$ we have
\begin{equation}
    x_n, y_n \in \psi_i\left(\mathbb{B}^m_{\frac{3}{2}}(0)\right),
    \quad
    \psi_i^{-1}(x_n) \rightarrow \psi_i^{-1}(q),
    \quad
    \psi_i^{-1}(y_n) \rightarrow
    \psi_i^{-1}(q).
\end{equation}
Doing a calculation similar to the estimate \eqref{eq:boundOfMetric}, we obtain
\begin{equation}
    \dist_g\left(
        x_n, y_n
    \right)
    \leq
    \sqrt{C_{UB}}
    \left\|
        \psi_i^{-1}\left(x_n\right)
    -
        \psi_i^{-1}\left(y_n\right)
    \right\|
    \rightarrow
    0.
\end{equation}  
This shows that
\begin{equation}
    \left[\left\{x_n\right\}_{n \in \mathbb{N}}\right]
    =
    \left[\left\{y_n\right\}_{n \in \mathbb{N}}\right],
\end{equation}
showing injectivity.

Now we show that $\Phi^{-1} : N \rightarrow \left(\overline{N\setminus S}, \overline{\dist}_g\right)$ is continuous. Since the chart topology of $N$ is first countable, we only need to prove sequential continuity. Given $x \in N$ and $y_n \in N$ such that $y_n \rightarrow x$, then there exists $i \in \{1,...,l\}$, $\tilde{x} \in \mathbb{B}^m_1(0)$ such that $x = \psi_i\left(\tilde{x}\right)$ and $n_0 \in \mathbb{N}$ such that for all $n \geq n_0$
\begin{equation}
    y_n \in \psi_i\left(\mathbb{B}^m_{\frac{5}{4}}(0)
    \right),
    \quad
    \psi_i^{-1}(y_n)
    \rightarrow
    \tilde{x}
    =
    \psi_i^{-1}(x).
\end{equation}
Consider the Cauchy sequences 
\begin{equation}
\label{eq:cauchySequencesFromN}
\Phi^{-1}(x) = \left[\left\{x^j\right\}_{j \in \mathbb{N}}\right], \quad \Phi^{-1}(y_n) = \left[\left\{y^j_n\right\}_{j \in \mathbb{N}}\right].
\end{equation}
This means that for each fixed $n \in \mathbb{N}$,  $\lim_j y^j_n = y_n$ and $\lim_j x^j = x$ in the chart topology. Using these convergences, the fact the sequences are Cauchy, and equations \eqref{eq:equivalenceOfCauhySequences} and \eqref{eq:cauchySequencesFromN} for each $n \geq n_0$, there exists $j_n \in \mathbb{N}$ such that
\begin{equation}
    \lim_j
    \dist_g\left(
        x^j, y_n^j
    \right)
    \leq
    \dist_g\left(
        x^{j_n}, y_n^{j_n}
    \right)
    +
    \frac{1}{n}.
\end{equation}
\begin{align}
    y_n^{j_n} \in \psi_i\left(\mathbb{B}^m_{\frac{3}{2}}(0)\right),&
    \quad
    \left\|
    \psi_i^{-1}\left(
        y_n^{j_n}
    \right)
    -
    \psi_i^{-1}(y_n)
    \right\|
    \leq
    \frac{1}{n},\\
    x^{j_n}
    \in \psi_i\left(\mathbb{B}^m_\frac{3}{2}(0)\right),&
    \quad
    \left\|
    \psi_i^{-1}\left(
        x^{j_n}
    \right)
    -
    \psi_i^{-1}(x)
    \right\|
    \leq
    \frac{1}{n}.
\end{align}
Define for $n \geq 2n_0$
\begin{equation}
    w_n:=
    \begin{cases}
        y_{\frac{n}{2}}^{j_{\frac{n}{2}}}&\text{ if } $n$ \text{ is even}\\
        x^{j_\frac{n+1}{2}}
        &\text{ if } $n$ \text{ is odd}
    \end{cases}.
\end{equation}
Since $\psi_i^{-1}(y_n) \rightarrow \psi_i^{-1}(x)$, we have $\psi_i^{-1}(w_n) \rightarrow \psi_i^{-1}(x)$ and so by \eqref{eq:equivalenceOfCauhySequences} we have that $w_n$ is Cauchy for $\dist_g$ and so $\limsup_n \dist_g(x^{j_n}, y_n^{j_n}) = 0$. Thus we conclude
\begin{align}
    \limsup_n
    \overline{\dist}_g
    \left(
        \Phi^{-1}(x),
        \Phi^{-1}(y_n)
    \right)
=&
    \limsup_n
    \lim_{j}
    \dist_g\left(
        x^j,
        y^j_n
    \right)\\
\leq&
    \limsup_n
    \dist_g
    \left(
        x^{j_n}, y_n^{j_n}
    \right) = 0.
\end{align} 
This shows that $\Phi^{-1}$ is continuous. Since $N$ is compact and $\overline{N\setminus S}$ is Hausdorff, the map $\Phi^{-1}$ is a homeomorphism, finishing the proof.

\end{proof}

\begin{rem}
    When $\mathcal{H}^{m-1}(S) > 0$, the conclusion Proposition \ref{prop:completionMetric} may fail.
    If $S = S^1$ then in $\mathbb{R}^2$ the completion of the distance $\dist_{g_\euc}$ induced by curves in $\mathbb{R}^2 \setminus S$ gives a space homeomorphic to
    \begin{equation}
        \overline{\mathbb{R}^2\setminus S}
    \cong
        (\mathbb{R}^2 \setminus \mathbb{B}^2_1(0))
    \coprod
        \overline{\mathbb{B}^2_1(0)},
    \end{equation}
    where $\coprod$ denotes the disjoint union of the two sets. This is not homeomorphic to $\mathbb{R}^2$ since it contains two different components. Similarly if $S = [-1,1]\times \{0\}$ one can show the homeomorphism
    \begin{equation}
        \overline{\mathbb{R}^2\setminus S}
    \cong
        \mathbb{R}^2
        \setminus
        \mathbb{B}^2_1(0),
    \end{equation}
\end{rem}

\begin{prop}
\label{prop:completionOfMetricIsParametrizationTopology}
    Let the tuple $(N,S,g,\{\psi_i\}_{i=1}^l)$ satisfy definition \eqref{def:riemannianManifoldWithSingularitySet} with dimension $m$, and given $0\leq k<m-1$ suppose that for all parametrizations $\mathcal{M}^{(m,k,+)}_{\mathfrak{E}}\left(\psi_i^{-1}\left(S\right)\right) < \infty$. Then the completion of $(N\setminus S, \dist_g)$ is homeomorphic to $N$ with the topology given by the parametrizations $\{\psi_i\}_{i=1}^l$. If $\nu(E) = \vol_g(E \setminus S)$, the metric measure space $\mathfrak{N} = (N, \overline{\dist}_g, \nu)$ satisfies $\mathcal{M}^{(m,k,+)}_{\mathfrak{N}}(S)<\infty$. Furthermore $\mathfrak{N}$ is a locally doubling metric measure space.
\end{prop}
\begin{proof}
    Because $\mathcal{M}_{\mathfrak{E}}^{(m,k,+)}(\psi_i^{-1}(S)) < + \infty$ and $k< m-1$, we have $\mathcal{M}_{\mathfrak{E}}^{(m,m-1,+)}(\psi_i^{-1}(S)) = 0$, and by Lemma \ref{lem:dimensionBoundHausdorffMinkowsky} we have $\mathcal{H}^{m-1}\left(\psi_i^{-1}\left(S\right)\right) = 0$. Thus we can apply Proposition \ref{prop:completionMetric} to conclude the completion of $(N\setminus S, \dist_g)$ is homeomorphic to $(N,\overline{\dist}_g)$. If $0<r<\frac{\sqrt{c_{LB}}}{2}$, using equation \eqref{eq:boundBellowMetric}, for a curve to start at $\psi_i\left(\mathbb{B}^m_1(0)\right)$ and leave $\psi_i\left(\mathbb{B}^m_2(0)\right)$ it must have a length strictly greater than $r$. Thus using equation \eqref{eq:boundBellowMetric} we have for $0<r < \frac{\sqrt{c_{LB}}}{2}$
    \begin{align}
        \psi_{i}^{-1}\left(
            T_{r,\overline{\dist}_g}\left(S\right)
        \right)
        \cap
        \mathbb{B}^m_1(0)
    \subset
        T_{\frac{r}{\sqrt{c_{LB}}}, \dist_\euc}
        \left(
            \psi_i^{-1}\left(S\right)
        \right).
    \end{align}
    Consider $\{e_j\}_{j=1}^m$ is an orthonormal basis for $\mathbb{R}^m$. For each $i \in \{1,...,l\}$ and $x \in \psi_i\left(\mathbb{B}^m_1(0)\right) \setminus S$, given the matrix $\left(g_{i,\psi_i^{-1}(x)}\right)_{jk} = g_x\left(\left(d\psi_i\right)e_j, \left(d\psi_i\right)e_k\right)$, by inequalities \eqref{eq:boundBellowMetric}, \eqref{eq:boundAboveMetric} we have
    \begin{equation}
        c_{LB}^\frac{m}{2}
    \leq
        \sqrt{\det g_{i,\psi_i^{-1}(x)}}
    \leq
        C_{UB}^\frac{m}{2}.
    \end{equation}
    Using the fact that for all $i \in \{1,\cdots,l\}$ we have $\mathcal{M}^{(m,k,+)}_{\mathfrak{E}}\left(\psi_i^{-1}\left(S\right)\right) < \infty$ and $\mathcal{L}^{m}\left(\psi_i^{-1}(S)\right) = 0$, there exist $r_0, C_{\mathcal{M}}>0$ such that
    \begin{align}
        \nu\left(
            T_{r,\overline{\dist}_g}
            \left(
                S
            \right)
        \right)
    \leq&
        \sum_{i=1}^l
        \nu\left(
            T_{r,\overline{\dist}_g}
            \left(
                S
            \right)
            \cap\psi_i\left(\mathbb{B}^m_1\left(0\right)\right)
        \right)\\
    =&
        \sum_{i=1}^l
        \nu\left(
            T_{r,\overline{\dist}_g}
            \left(
                S
            \right)
            \cap\psi_i\left(\mathbb{B}^m_1\left(0\right)\right)
            \setminus S
        \right)\\
    =&
        \sum_{i=1}^l
        \int_{\psi_{i}^{-1}\left(
            T_{r,\overline{\dist}_g}\left(S\right) \setminus S
        \right) \cap \mathbb{B}^m_1(0)}
            1
            \cdot
            \sqrt{\det(g_{i,y})}
            d\mathcal{L}^m(y)\\
    \leq&
        \sum_{i=1}^l
        \int_{T_{\frac{r}{\sqrt{c_{LB}}}, \dist_\euc}
        \left(
            \psi_i^{-1}\left(S\right)
        \right)}C_{UB}^{\frac{m}{2}}
            d\mathcal{L}^m(y)\\
    \leq&
        \sum_{i=1}^l
        C_{UB}^{\frac{m}{2}}
        \mathcal{L}^m\left(
            T_{\frac{r}{\sqrt{c_{LB}}}, \dist_\euc}
        \left(
            \psi_i^{-1}\left(S\right)
        \right)
        \right)\\
    \leq&
        \sum_{i=1}^l
        C_{UB}^{\frac{m}{2}}
        C_{\mathcal{M}}
        \frac{r^{m-k}}{c_{LB}^{\frac{m-k}{2}}}.
    \end{align}

    Consider $r_0$ such that $\frac{r_0}{c_{LB}^\frac{1}{2}} \leq \frac{1}{2}$. By Proposition \ref{prop:completionMetric}, we have that the distance $\overline{\dist}_g(x,y)$ can be computed by taking sequences $x_n,y_n \in N \setminus S$ such that $x_n \rightarrow x$ and $y_n \rightarrow y$, and  computing the $\lim_n \dist_g\left(x_n, y_n\right)$, where $\dist_g\left(x_n, y_n\right)$ is given by the infimum of curves in $N \setminus S$ from $x_n$ to $y_n$. Fix $x \in N$ and $i \in \{1,...,l\}$ such that $x \in \psi_i\left(\mathbb{B}^m_1(0)\right)$. Using the previous fact, equation
    \eqref{eq:boundBellowMetric} and Proposition \ref{prop:completionMetric} and the choice of $r_0$ we can conclude
    \begin{equation}
        B_{r}(x) \subset \psi_i\left(\mathbb{B}^m_2(0)\right).
    \end{equation}
    We can argue similarly using
    \eqref{eq:boundBellowMetric},
    \eqref{eq:boundAboveMetric} and Proposition \ref{prop:completionMetric} to conclude
    \begin{equation}
        \mathbb{B}^m_{r/C_{UB}^\frac{1}{2}}
        \left(\psi_i^{-1}\left(x\right)\right)
    \subset
        \psi_i^{-1}\left(B_{r}(x)\right)
    \subset
        \mathbb{B}^m_{r/c_{LB}^\frac{1}{2}}
        \left(\psi_i^{-1}\left(x\right)\right).
    \end{equation}

    Thus
    \begin{align}
        \nu\left(B_r(x)\right)
    \leq&
        \int_{\psi_i^{-1}\left(B_r(x)\right)}
        |\det g_{i,y}|^\frac{1}{2}
        d\mathcal{L}^m(y)\\
    \leq&
        C_{UB}^\frac{m}{2}
        \mathcal{L}^m\left(\mathbb{B}^m_{r/c_{LB}^\frac{1}{2}}
        \left(\psi_i^{-1}\left(x\right)\right)\right)\\
    \leq&
        \frac{C_{UB}^\frac{m}{2}}{c_{LB}^\frac{m}{2}}
        \mathcal{L}^m\left(\mathbb{B}^m_1(0)\right)
        r^m.
    \end{align}
    Similarly we can prove that
    \begin{equation}
        \nu(B_r(x))
    \geq
        \frac{c_{LB}^{\frac{m}{2}}}
        {C_{UB}^\frac{m}{2}}
        \mathcal{L}^m\left(
            \mathbb{B}^m_1(0)
        \right)
        r^m.
    \end{equation}
    Hence the space is locally Ahlfors regular, thus it is locally doubling for $r<\frac{r_0}{2}$, finishing the proof.
\end{proof}

\subsection{Injectivity radius of the exponential map of $(N\setminus S, g)$}

In the following subsection, we assume the conditions of Proposition \ref{prop:completionOfMetricIsParametrizationTopology}, so that the completion $\left(\overline{N\setminus S}, \overline{\dist}_g\right)$ can be identified with $N$. The goal of this subsection is to prove the following proposition.

\begin{prop}
\label{prop:injectivityRadiusProp}
    Let the tuple $(N,S,g,\left\{\psi_i\right\}_{i=1}^l)$ satisfy Definition \ref{def:riemannianManifoldWithSingularitySet} with dimension $m$. Furthermore, suppose that $(N\setminus S, g)$ satisfies \eqref{eq:curvatureCondition}. Then there exists $r_0, C_{\alpha,\inj}>0$ such that for all $x \in N\setminus S$
    \begin{equation}
        \inj(x)
    \geq
        \min\left\{
            r_0,
            \left(
                \frac{\overline{\dist}_g(x,S)}{C_{\alpha, \inj}}
            \right)^{\max\{1,\frac{\alpha}{2}\}}
        \right\}.
    \end{equation}
\end{prop}

We first prove an auxiliary lemma

\begin{lemma}
\label{lem:homotopyContinuation}
    Let $\tau, \sigma > 0$, and suppose that $\exp_x: B_{\tau}(0_x) \rightarrow B_\tau(x)$ is such that for all $z \in B_\tau(0_x)$ and $u \in \partial B_1(0_x)$
    \begin{equation}
        \sigma
    \leq
        g_{\exp_x(z)}
        \left(
            (d\exp_x)_{z}
            u,
            (d\exp_x)_{z}
            u
        \right)
    \leq
        \frac{1}{\sigma}.
    \end{equation}
    Let $0<s<\tau$, $0<\delta\leq 1$ and consider a smooth map $H:[0,1]\times [0,\delta] \rightarrow \exp_x\left(\overline{B_s(0_x)}\right)$ and a lift $\tilde{H}:[0,1]\times [0,\delta[ \rightarrow \overline{B_s(0_x)}$ such that 
    \begin{equation}
        \exp_x\left(\tilde{H}(\nu,t)\right)
        =
        H(\nu,t).
    \end{equation}
    The lift then extends continuously to $[0,1]\times [0,\delta]$ with image on $\overline{B_s(0_x)}$.

    Furthermore if $H:[0,1] \times [0,1] \rightarrow \exp_{x}\left(B_\tau(0_x)\right)$ is a homotopy and there exist $0\leq \delta < 1$ and a lift $\tilde{H}:[0,1]\times [0,\delta] \rightarrow B_\tau(0_x)$, then the lift extends continuously to a lift in $[0,1]\times [0,\delta + \epsilon]$, for some $\epsilon>0$ with image in $B_\tau(0_x)$.
\end{lemma}
\begin{proof}
    Since $H$ is smooth, there exists $0\leq L<\infty$ such that $\sup_{(\nu,t) \in [0,1]\times [0,\delta]}\left|\left(d H\right)_{(\nu,t)}\right| = L$. 

    Because $d\exp_x$ has a uniform lower bound, $\exp_x$ is a local diffeomorphism. Thus, the inverse function theorem implies that every local branch of the lift $\tilde{H}$ is smooth, and so $\tilde{H}$ is smooth. In particular, for all $(\nu,t) \in [0,1]\times [0,\delta[$
    \begin{equation}
        \left|
            \left(d
            \tilde{H}\right)_{(\nu,t)}
        \right|
        \leq
        \frac{1}{\sqrt{\sigma}}
        \left|
            \left(
                d \exp_x
            \right)_{\tilde{H}(\nu,t)}
            \cdot
            \left(d
            \tilde{H}\right)_{(\nu,t)}
        \right|
        =
        \frac{1}{\sqrt{\sigma}}\left|
            (d H)_{(\nu,t)}
        \right|
        \leq
         \frac{L}{\sqrt{\sigma}}.
    \end{equation}
    Therefore $\tilde{H}$ is uniformly Lipschitz. Because $[0,1]\times [0,\delta[$ is dense in $[0,1]\times[0,\delta]$, and both $[0,1]\times [0,\delta]$, and $\overline{B_s(0_x)}$ are complete, there is a unique Lipschitz extension of $\tilde{H}$ to $[0,1]\times [0,\delta]$. Since $\exp_x$ and $H$ are continuous, this extension is also a lift, showing the first part of the proof.

    For the second part, note that $\tilde{H}\left([0,1]\times [0,\delta]\right)\subset B_\tau(0_x)$ is compact. Thus, we can find a finite open cover $\left\{U_n\right\}_{n=1}^k$ of $\tilde{H}\left([0,1]\times [0,\delta]\right)$ such that each restriction $\left.\exp_x\right|_{U_n} : U_n \rightarrow \exp_x\left(U_n\right)$ is a diffeomorphism onto the image. Furthermore, since $\tilde{H}$ is continuous and $[0,1]$ is compact, consider $\sigma \in \mathbb{N}$ such that for all $i \in \{0,...,\sigma-1\}$ there exists $n_i \in \{1,...,k\}$ such that
    \begin{equation}
        \tilde{H}\left(
            \left[\frac{i}{\sigma},\frac{i+1}{\sigma}\right]\times \{\delta\}
        \right)
        \subset
        U_{n_i},
        \quad
        H\left(
            \left[\frac{i}{\sigma},\frac{i+1}{\sigma}\right]\times \{\delta\}
        \right)
        \subset 
        \exp_x\left(
            U_{n_i}
        \right).
    \end{equation}
    Since the cover is finite and $H$ is continuous, there exists $\epsilon>0$ such that $\delta+\epsilon \leq 1$, and for all $i \in \{0,...,\sigma-1\}$ we have
    \begin{equation}
        H\left(
            \left[\frac{i}{\sigma},\frac{i+1}{\sigma}\right]\times \left[\delta, \delta + \epsilon\right]
        \right)
        \subset
        \exp_x\left(
            U_{n_i}
        \right).
    \end{equation}
    Define the extension $\hat{H}:[0,1]\times [0,\delta + \epsilon]$ by
    \begin{equation}
        \hat{H}(\nu,t)
    =
        \begin{cases}
            \left.
                \tilde{H}(\nu,t)
            \right|_{[0,1]\times[0,\delta]}
            &
            \text{ if } t \in [0,\delta],\\
            \left(\left.\exp_x\right|_{U_{n_i}}\right)^{-1}
            \circ
            H(\nu,t)
            &
            \text{ if }
            t \in [\delta, \delta+\epsilon]
            \text{ and }
            \nu \in \left[\frac{i}{\sigma},\frac{i+1}{\sigma}\right]
        \end{cases}
    \end{equation}
    Because $\exp_x$ is local diffeomorphism, path lifts with the same initial value are unique. Hence
    \begin{equation}
        \left(\left.\exp_x\right|_{U_{n_{i-1}}}\right)^{-1}
            \circ
            H\left(\frac{i}{\sigma},t\right)
        =
        \left(\left.\exp_x\right|_{U_{n_i}}\right)^{-1}
            \circ
            H\left(\frac{i}{\sigma},t\right).
    \end{equation}
    At $t = \delta$, because $\left.\tilde{H}(\nu,\delta)\right|_{[0,1]\times [0,\delta]} \in U_{n_i}$, we must have
    \begin{equation}
        \left.\tilde{H}(\nu,\delta)\right|_{[0,1]\times [0,\delta]}
        =
        \left(\left.\exp_x\right|_{U_{n_i}}\right)^{-1}
            \circ
            H(\nu,\delta).
    \end{equation}
    Thus, this map is a continuous extension of $\tilde{H}$ to $[0,1]\times [0,\delta+\epsilon]$, and is still a lift of $H$ by $\exp_x$, concluding the proof.
\end{proof}

Throughout the proof take the convention that for $\alpha = 0$, then $\min\left\{1,\frac{2}{\alpha}\right\} = 1$.

\begin{lemma}
\label{lem:exponentialContainsBall}
    Let the tuple $(N,S,g,\left\{\psi_i\right\}_{i=1}^l)$ satisfy definition \ref{def:riemannianManifoldWithSingularitySet} with dimension $m$ and assume that \eqref{eq:curvatureCondition} holds. Set $t_r:=80\sqrt{\frac{C_{UB}}{c_{LB}}}r > r$. Then there exist $r_0, C_{\alpha,\inj}>0$ such that for all $0<r<r_0$, $i \in \{1,...,l\}$, and $x \in \psi_i\left(\mathbb{B}^m_1(0)\right)$ satisfying $\overline{\dist}_g(x,S) > C_{\alpha,\inj} r^{\min\left\{1,\frac{2}{\alpha}\right\}}$ we have
    \begin{equation}
    \label{eq:ExpContainsBall}
        \psi_i\left(
            \mathbb{B}^m_{\frac{4r}{\sqrt{c_{LB}}}}\left(
                \psi_i^{-1}\left(x\right)
            \right)
        \right)
    \subset
        \exp_x
        \left(
            \overline{B_{t_r}(0_x)}
        \right).
    \end{equation}
\end{lemma}
\begin{proof}
    By Lemma \ref{lem:exponentialMapDerivative}, there exist $0<C_{\alpha,\dist_g}$ and $0<r_0$, and for $0<r<r_0$ such that for al $x \in N\setminus S$ satisfying $\overline{\dist}_g(x,S) > C_{\alpha,\dist_g}\left(160\sqrt{\frac{C_{UB}}{c_{LB}}} r\right)^{\min\left\{1,\frac{2}{\alpha}\right\}}$ we have that \eqref{eq:ineqsForExponentialJacobian} is satisfied. In particular \eqref{eq:ineqsForExponentialJacobian} implies that $\exp:B_{2t_r}(0_x) \rightarrow N\setminus S$ does not have critical points. Furthermore, by decreasing $r_0$ so that $4t_{r_0} < \frac{\sqrt{c_{LB}}}{2}$, for all $i \in \{1,...,l\}$, $0<r<r_0$ and $x \in \psi_i\left(\mathbb{B}^m_1(0)\right)$ then $B_{4t_r}(x) \subset \psi_i\left(\mathbb{B}^m_\frac{3}{2}(0)\right)$. Take $C_{\alpha,\inj} \geq C_{\alpha,\dist_g}\left(160\sqrt{\frac{C_{UB}}{c_{LB}}}\right)^{\min\left\{1,\frac{2}{\alpha}\right\}}$. Let $\tilde{x} = \psi_i^{-1}\left(x\right)$.
    Suppose, for contradiction, that
    \begin{equation}
        \psi_i\left(
            \mathbb{B}^m_{\frac{4r}{\sqrt{c_{LB}}}}\left(
                \tilde{x}
            \right)
        \right)
    \not\subset
        \exp_x
        \left(
            \overline{B_{t_r}(0_x)}
        \right),
    \end{equation}
    and let $z_{\text{out}} \in \mathbb{B}^m_{\frac{4r}{\sqrt{c_{LB}}}}
    \left(\tilde{x}\right)$ such that $\psi_i(z_{\text{out}}) \notin \exp_x
        \left(
            \overline{B_{t_r}(0_x)}
        \right)$. Consider the curve $\eta(t) = (1-t)\tilde{x} + tz_{\text{out}}$. Define the first exit time
        \begin{equation}
            t_\delta := 
            \inf
            \left\{
                \delta \in [0,1]:
                \psi_i\left(\eta(\delta)\right)
                \notin
                \exp_x\left(
                    \overline{B_{t_r}(0_x)}
                \right)
            \right\}
            < 1.
        \end{equation}
    Let $\tilde{z} = \eta(t_\delta) \in \mathbb{B}^m_{\frac{4r}{\sqrt{c_{LB}}}}\left(
                \tilde{x}
            \right)$, Since $\exp_x\left(\overline{B_{t_r}(0_x)}\right)$ is compact and $\exp_x$ is a local diffeomorphism around $0_x$, then
            the curve $\tilde{\gamma}(t) = (1-t)\tilde{x} + t\tilde{z}$ satisfies
    \begin{align}
        \psi_i\left(
            \tilde{\gamma}\left(\left[0,1\right]\right)
        \right)
    \subset&
        \exp_x
        \left(
            \overline{B_{t_r}(0_x)}
        \right),\\
        \psi_i\left(\tilde{\gamma}(1)\right) = \psi_i\left(\tilde{z}\right)
        \in&
        \partial
        \left(
            \exp_x\left(
                \overline{B_{t_r}(0_x)}
            \right)
        \right).
        \label{eq:pointInBoundary}
    \end{align}

        Let $t_r<s< \frac{3t_r}{2}$. Since $\left.\exp_x\right|_{B_{2t_r}(0_x)}$ has no critical points, it is a local diffeomorphism. Thus, there exists $0<\delta \leq 1$ and a unique lift $\hat{\gamma}:[0,\delta] \rightarrow B_{s}(0_x)$ such that
        \begin{align}
        \label{eq:propertiesOfLift1}
        \gamma(t)
        :=
        \psi_i\left(
            \tilde{\gamma}(t)
        \right)
        =
        \exp_x\left(
            \hat{\gamma}(t)
        \right),
        \quad
        \hat{\gamma}(0) = 0_x.
        \end{align}
        Now let
        \begin{equation}
        \label{eq:maxExtension}
            \delta_{\text{max}}
        :=
        \sup
            \left\{
                \delta \in ]0,1]:
                \text{There exists a unique lift } 
                \hat{\gamma}:[0,\delta] \rightarrow B_s(0_x)
                \text{ satisfying }
                \eqref{eq:propertiesOfLift1}
            \right\}.
        \end{equation}
        First suppose $\delta_{\text{max}} < 1$. We must have a lift $\hat{\gamma}:[0,\delta_{\text{max}}[ \rightarrow B_s(0_x)$ satisfying \eqref{eq:propertiesOfLift1},, since the lifts obtained for smaller $\delta$ all agree by uniqueness of the lift, since $\exp_x$ is a local diffeomorphism in $B_s(0_x)$ and the lift has the same initial point.

        Since $\gamma$ is smooth and in $B_{\frac{3t_r}{2}}(0_x)$ we have \eqref{eq:ineqsForExponentialJacobian}, by taking $H(\nu,t) = \gamma(t)$ and $\hat{H}(\nu,t) = \hat{\gamma}(t)$, we can use Lemma \ref{lem:homotopyContinuation} to extend the curve $\hat{\gamma}$ continuously to $[0,\delta_{\max}]$, preserving the fact that it is a lift of $\gamma$. We must have that $\hat{\gamma}(\delta_{\max}) \in \partial B_s(0_x)$, otherwise $\hat{\gamma}([0,\delta_{\max}]) \subset B_s(0_x)$, and by Lemma \ref{lem:homotopyContinuation}, we can extend the lift $\hat{\gamma}$ to $[0,\delta_{\max} + \epsilon]$ for some $\epsilon>0$ while the image remains in $B_s(0_x)$, contradicting the definition of $\delta_{\max}$.

        Since $\hat{\gamma}(0) = 0_x$ and $\hat{\gamma}(\delta_{\max}) \in \partial B_s(0_x)$, we must have
        \begin{equation}
            l_\euc\left(
                \hat{\gamma}
                \left(
                    [0,\delta_{\max}]
                \right)
            \right)
            \geq
            s
            \geq
            t_r,
        \end{equation}      
        where $l_\euc$ is the Euclidean length of a path. On the other hand, using equations \eqref{eq:boundAboveMetric}, \eqref{eq:ineqsForExponentialJacobian} we obtain
    \begin{align}
        l_{\text{euc}}\left(\tilde{\gamma}\left([0,\delta_{\text{max}}[\right)\right)
    &=
        \int_0^{\delta_{\text{max}}}
            \|\tilde{\gamma}'(t)\|
            dt\\
    &=
        \int_0^{\delta_{\text{max}}}
            \left\|d(\psi_i^{-1}\circ \exp_x)_{\hat{\gamma}(t)}\cdot\hat{\gamma}'(t)\right\|dt\\
    &\geq
        \int_0^{\delta_{\text{max}}}
        \frac{1}{\sqrt{C_{UB}}}
        g\left(
            \left(d\exp_x\right)_{\hat{\gamma}(t)}
            \hat{\gamma}'(t),
            \left(d\exp_x\right)_{\hat{\gamma}(t)}
            \hat{\gamma}'(t)
        \right)^\frac{1}{2}
        dt
        \\
    &\geq
        \int_0^{\delta_{\text{max}}}
            \frac{1}{\sqrt{2C_{UB}}}
            \|\hat{\gamma}'(t)\|dt
    =
        \frac{1}{\sqrt{2C_{UB}}}
        l_{\text{euc}}\left(\hat{\gamma}\left([0,\delta_\text{max}[\right)\right)
        \label{eq:lengthComparison1}
    \end{align}
    Since $\tilde{\gamma}$ is a line segment, starting at $\tilde{x}$ contained in $\mathbb{B}^m_{\frac{4r}{\sqrt{c_{LB}}}}(\tilde{x})$, then $l_\euc\left(\tilde{\gamma}\right) \leq \frac{4r}{\sqrt{c_{LB}}}$. Therefore,
    \begin{equation}
    \label{eq:lengthComparison2}
            80\sqrt{\frac{2C_{UB}}{c_{LB}}}r
        =
            \sqrt{2} t_r
        \leq
            \sqrt{2}
            l_{\text{euc}}(\hat{\gamma}\left(
                [0,\delta_{\max}]
            \right)
            )
        \leq
            2\sqrt{C_{UB}}
            l_{\text{euc}}\left(
            \tilde{\gamma}\left(
                [0,\delta_{\max}]
            \right)
            \right)
        \leq
            8\sqrt{\frac{C_{UB}}{c_{LB}}}r,
        \end{equation}
    which is a contradiction. Thus we must have $\delta_{\text{max}} = 1$.

    Since $\exp_x$ is a local diffeomorphism, there exists a unique lift $\hat{\gamma}:[0,1[ \rightarrow B_s(0_x)$ for all $t_r < s < \frac{3 t_r}{2}$, that agrees for different $s$. Since the intersection of these sets is $\overline{B_{t_r}(0_x)}$, the unique lift must be contained in $\overline{B_{t_r}(0_x)}$. Again we can use Lemma \ref{lem:homotopyContinuation} to extend this lift to the interval $[0,1]$ and with image $\hat{\gamma}([0,1])$ contained in $\overline{B_{t_r}(0_x)}$.

    We must have $\hat{\gamma}(1) \in \partial B_{t_r}(0_x)$, otherwise if $\hat{\gamma}(1) \in B_{t_r}(0_x)$, since $\exp_x$ is an open map, and $\hat{\gamma}$ a lift of $\gamma$, we would have $ \psi_i\left(\tilde{\gamma}(1)\right) = \exp_x\left(\hat{\gamma}(1)\right) \in \text{int}\left(\exp_x\left(\overline{B_{t_r}(0_x)}\right)\right)$, contradicting equation \eqref{eq:pointInBoundary}. Thus $\hat{\gamma}(1) \in \partial B_{t_r}(0_x)$.  We can then use the same argument as in the inequalities \eqref{eq:lengthComparison1} and \eqref{eq:lengthComparison2} , to obtain a contradiction. This finishes the proof.
       
\end{proof}

\begin{lemma}
\label{lem:injectivityLemma}
    Let the tuple $(N,S,g,\left\{\psi_i\right\}_{i=1}^l)$ satisfy definition \ref{def:riemannianManifoldWithSingularitySet} with dimension $m$ and \eqref{eq:curvatureCondition}. There exist $r_0, C_{\alpha,\inj}>0$ such that for all $0<r<r_0$ and $x \in N \setminus S$ satisfying $\overline{\dist}_g(x,S) > C_{\alpha,\inj} r^{\min\left\{1,\frac{2}{\alpha}\right\}}$, the map $\exp_x:B_r(0_x) \rightarrow B_r(x)$ is injective.
\end{lemma}
\begin{proof}
    Consider $C_{\alpha,\inj}>0$ large enough and $r_0$ small enough so that by Lemmas \ref{lem:exponentialMapDerivative} and \ref{lem:exponentialContainsBall}, for all $0<r<r_0$, if $\overline{\dist}_g(x,S) > C_{\alpha,\dist_g}\left(160\sqrt{\frac{C_{UB}}{c_{LB}}} r\right)^{\min\left\{1,\frac{2}{\alpha}\right\}}$ then 
    \eqref{eq:ExpContainsBall} is satisfied, and \eqref{eq:ineqsForExponentialJacobian} is satisfied for the ball $B_{2t_r}(0_x)$.

    Furthermore, by equation \eqref{eq:boundBellowMetric}, by choosing $r_0$ so that $4t_{r_0} < \frac{\sqrt{c_{LB}}}{2}$, for all $i \in \{1,...,l\}$ and $x \in \psi_i\left(\mathbb{B}^m_1(0)\right)$, then $B_{4t_r}(x) \subset \psi_i\left(\mathbb{B}^m_{\frac{3}{2}}(0)\right)$.

    Fix $i \in \{1,...,l\}$ and $x \in N\setminus S$ such that $x \in \psi_i\left(\mathbb{B}^m_1\left(0\right)\right)$ and $\overline{\dist}_g(x,S) > C_{\alpha,\inj}r^{\min\left\{1,\frac{2}{\alpha}\right\}}$. Now consider $\hat{y}_0, \hat{y}_1 \in B_r(0_x)$ such that $\exp_x(\hat{y}_0) = \exp_x(\hat{y}_1)$ and define the curves $\hat{\gamma}_j(t) = t\hat{y}_j$ and $\tilde{\gamma}_j(t) = \psi_i^{-1}\circ \exp_x\left(\hat{\gamma}_j(t)\right)$ for $j \in \{0,1\}$. Using equation \eqref{eq:boundBellowMetric} we have
    \begin{equation}
        \sqrt{c_{LB}}
        l_{\euc}\left(
            \tilde{\gamma}_j
        \right)
        \leq
        l_{\euc}
        \left(
            \hat{\gamma}_j
        \right)
        \leq
        r.
    \end{equation}
    Thus we have
    \begin{equation}
        \tilde{\gamma}_j\left(\left[0,1\right]\right)
    \subset
        \mathbb{B}^m_{\frac{4 r}{\sqrt{c_{LB}}}}\left(\psi_i^{-1}\left(x\right)\right).
    \end{equation}
    and notice that $\tilde{\gamma}_0(0) = \tilde{\gamma}_1(0) = \psi_i^{-1}(x)$ and $\tilde{\gamma}_0(1) = \tilde{\gamma}_1(1)$. Define the homotopy $\tilde{H}:[0,1]^2 \rightarrow \mathbb{B}^m_{\frac{4r}{\sqrt{c_{LB}}}}\left(\psi_i^{-1}\left(x\right)\right)$ by
    \begin{align}
        \tilde{H}(t,s) = (1-s)\tilde{\gamma}_0(t)+s\tilde{\gamma}_1(t).
    \end{align}
    Notice that by \eqref{eq:ExpContainsBall} we have $\psi_i\left(\tilde{H}\left([0,1]^2\right)\right) \subset \exp_x\left(\overline{B_{t_r}(0_x)}\right)$ with $t_r = \frac{80\sqrt{C_{UB}}}{\sqrt{c_{LB}}}r$. We now prove that there exists a unique lift $\hat{H}:[0,1]^2\rightarrow \overline{B_{\frac{5t_r}{4}}(0_x)}$ such that
    \begin{align}
    \label{eq:homotopyLiftProperty}
    \exp_x\left(\hat{H}\left(t,s\right)\right)
    =
        \psi_i\left(\tilde{H}\left(t,s\right)\right),\quad
        \hat{H}\left(t,0\right)
    =
        \hat{\gamma}_0(t)
    =
        t\hat{y}_0.
    \end{align}
    If such a lift exists, it is unique because $\exp_x$ is a local diffeomorphism in $B_{\frac{3t_r}{2}}(0_x)$.

    Because estimate \eqref{eq:ineqsForExponentialJacobian} holds on the larger ball $B_{\frac{3t_r}{2}}(0_x)$, Lemma \ref{lem:homotopyContinuation} applied with $\delta = 0$ gives $t_\delta>0$ and a unique lift $\hat{H}:[0,1]\times [0,t_\delta] \rightarrow B_{\frac{5t_r}{4}}(0_x)$, satisfying \eqref{eq:homotopyLiftProperty}. We can thus define
    \begin{equation}
    \label{eq:maxTimeForHomotopy}
        \delta_{\text{max}}
    :=
        \sup
        \left\{
            \delta \in ]0,1]:
            \text{There exists a lift } \hat{H}:[0,1]\times [0,\delta] \rightarrow B_{\frac{5t_r}{4}}(0_x)
            \text{ satisfying }
            \eqref{eq:homotopyLiftProperty}
        \right\}.
    \end{equation}
    Suppose first that $\delta_{\text{max}}<1$. By compatibility of the lifts on nested domains (since $\exp_x$ is a local diffeomorphism in $B_{\frac{5t_r}{4}}(0_x)$ and the initial condition is the same), they define a unique lift on the half-open maximal domain $\hat{H}:[0,1]\times [0,\delta_{\text{max}}[ \rightarrow B_{\frac{5t_r}{4}}(0_x)$.

    By Lemma \ref{lem:homotopyContinuation} the lift extends continuously to $\hat{H}:[0,1]\times[0,\delta_{\max}] \rightarrow \overline{B_{\frac{5t_r}{4}}(0_x)}$. There must also exist some $t_{\delta_{\max}} \in [0,1]$ such that $\hat{H}(t_{\delta_{\max}},\delta_{\max}) \in \partial B_{\frac{5t_r}{4}}(0_x)$, otherwise $\hat{H}\left([0,1]\times [0,\delta_{\max}]\right) \subset B_{\frac{5 t_r}{4}}(0_x)$, and we can use Lemma \ref{lem:homotopyContinuation} to extend the lift to $\hat{H}:[0,1]\times [0,\delta_{\max}+\epsilon] \rightarrow B_{\frac{5t_r}{4}}(0_x)$, contradicting the maximality of \eqref{eq:maxTimeForHomotopy}.

    Thus there exists $t_{\delta_{\max}} \in [0,1]$ such that $\hat{H}\left(t_{\delta_{\max}}, \delta_{\max}\right) \in \partial B_{\frac{5t_r}{4}}(0_x)$. Since $\hat{H}(t_{\delta_{\max}}, 0) = \hat{\gamma}_0(t_{\delta_{\max}}) = t_{\delta_{\max}} \hat{y}_0 \in \overline{B_r(0_x)} \subset \overline{B_{t_r}(0_x)}$, we must have
    \begin{equation}
        l_\euc
        \left(
            \hat{H}
            \left(
                t_{\delta_{\max}}
                ,
                [0,\delta_{\max}]
            \right)
        \right)
        \geq
        \frac{t_r}{4}.
    \end{equation}      
    Using equations \eqref{eq:boundAboveMetric} and \eqref{eq:ineqsForExponentialJacobian}, we obtain
    \begin{align}
        &l_{\text{euc}}\left(\tilde{H}\left(t_{\delta_{\max}},[0,\delta_{\text{max}}[\right)\right))\\
    =&
        \int_0^{\delta_{\text{max}}}
            \left\|
            \frac{d}{ds}\tilde{H}(t_{\delta_{\max}},s)
            \right\|
            ds\\
    =&
        \int_0^{\delta_{\text{max}}}
            \left\|d(\psi_i^{-1}\circ \exp_x)_{\hat{H}(t_{\delta_{\max}},s)}\cdot
            \frac{d}{ds}\hat{H}(t_{\delta_{\max}},s)\right\|ds\\
    \geq&
        \int_0^{\delta_{\text{max}}}
        \frac{1}{\sqrt{C_{UB}}}
        g\left(
            \left(d\exp_x\right)_{\hat{H}(t_{\delta_{\max}},s)}
            \cdot
            \frac{d}{ds}
            \hat{H}(t_{\delta_{\max}},s),
            \left(d\exp_x\right)_{\hat{H}(t_{\delta_{\max}},s)}
            \cdot
            \frac{d}{ds}
            \hat{H}(t_{\delta_{\max}},s)
        \right)^\frac{1}{2}
        ds
        \\
    \geq&
        \int_0^{\delta_{\text{max}}}
            \frac{1}{\sqrt{2C_{UB}}}
            \left\|\frac{d}{ds}\hat{H}(t_{\delta_{\max}},s)
            \right\|ds
    =
        \frac{1}{\sqrt{2C_{UB}}}
        l_{\text{euc}}\left(\hat{H}\left(t_{\delta_{\max}},[0,\delta_\text{max}[\right)\right).
    \end{align}
    The curve $s \in [0,\delta_{\text{max}}] \mapsto \tilde{H}(t_{\delta_{\max}},s)$, is a line segment contained in $\mathbb{B}^m_{\frac{4r}{\sqrt{c_{LB}}}}\left(\psi_i^{-1}\left(x\right)\right)$. Hence its length is at most
    \begin{equation}
        l_{\euc}\left(
            \tilde{H}(t_{\delta_{\max}},[0,\delta_{\text{max}}])
        \right)
        \leq
        \frac{8r}{\sqrt{c_{LB}}}.
    \end{equation}  
    This implies that
    \begin{align}
        \frac{10r\sqrt{C_{UB}}}{\sqrt{c_{LB}}}
    =
        \frac{t_r}{8}
    \leq
        \frac{1}{\sqrt{2}}
        l_{\euc}
        \left(
            \hat{H}\left(t_{\delta_{\max}},\left[0,\delta_{\text{max}}\right]\right)
        \right)
    &\leq
        \sqrt{C_{UB}}
        l_{\euc}\left(
            \tilde{H}(t_{\delta_{\max}},[0,\delta_{\text{max}}])
        \right)\\
    &\leq
        \frac{8r\sqrt{C_{UB}}}{\sqrt{c_{LB}}},
    \end{align}
    which is a contradiction. Thus we conclude that $\delta_{\max} = 1$. This shows that there exists a lift $\hat{H}:[0,1]\times [0,1[ \rightarrow B_{\frac{5t_r}{4}}(0_x)$, and using Lemma \ref{lem:homotopyContinuation}, we can extend the lift to $\hat{H}:[0,1]\times [0,1]\rightarrow \overline{B_{\frac{5t_r}{4}}(0_x)}$.

     The lifts are unique since $\exp_x$ is a local diffeomorphism. Given the initial point $\hat{H}(0,0) = 0_x$ and $\exp_x\left(\hat{H}\left(0,s\right)\right) = x $ and $\exp_x\left(\hat{H}(t,1)\right) = \exp_x(t\hat{y}_1)$ we must have
    \begin{equation}
        \hat{H}(0,s) = 0_x,
    \quad\quad
        \hat{H}(t,1) = t\hat{y}_1. 
    \end{equation}
    Furthermore since $\exp_x\left(\hat{H}\left(1, s\right)\right) = \exp_x\left(\hat{y}_0\right) = \exp_x\left(\hat{y}_1\right)$ by uniqueness of lifts of paths, since $\hat{H}(1,0) = \hat{y}_0$, this implies $\hat{H}(1,s) = \hat{y}_0$ for all $s \in [0,1]$. Thus, $\hat{y}_0 = \hat{H}(1,1) = \hat{y}_1$, showing that $\exp_x: B_r(0_x) \rightarrow B_r(x)$ is injective, finishing the proof.
    
\end{proof}

\begin{lemma}
\label{lem:surjectivityLemma}
    Let the tuple $(N,S,g,\left\{\psi_i\right\}_{i=1}^l)$ satisfy definition \ref{def:riemannianManifoldWithSingularitySet} with dimension $m$ and \eqref{eq:curvatureCondition}. There exist $r_0, C_{\alpha,\inj}>0$ such that for all $0<r<r_0$ and $x \in N \setminus S$ satisfying $\overline{\dist}_g(x,S) > C_{\alpha,\inj} r^{\min\left\{1,\frac{2}{\alpha}\right\}}$, the map $\exp_x:B_r(0_x) \rightarrow B_r(x)$ is surjective.
\end{lemma}
\begin{proof}
    Choose $r_0, C_{\alpha, \inj} > 0$ such that for all $0<r<r_0$ and $x \in N\setminus S$ such that $\overline{\dist}_g\left(x,S\right) > C_{\alpha,\inj} r^{\min\left\{1,\frac{2}{\alpha}\right\}}$ by Lemmas \ref{lem:exponentialMapDerivative} and \ref{lem:injectivityLemma} then the map $\exp_x : B_r(0_x) \rightarrow B_r(x)$ is injective and a local diffeomorphism. Suppose for contradiction $\exp_x: B_r(0_x) \rightarrow B_r(x)$ is not surjective. Thus there exists $0<s_1<r$ and $z \in B_{s_1}(x)$ such that $z \notin \exp_x(B_r(0_x))$. Choose $0<s_2$ such that $0<s_1<s_2<r$. This choice implies $z \notin \exp_x\left(\overline{B_{s_2}(0_x)}\right)$. Since $z \in B_{s_1}(x)$ there exists a smooth curve $\gamma: [0,1] \rightarrow N \setminus S$ such that $\gamma(0) = x$, $\gamma(1) = z$, and
    \begin{align}
        \int_0^1
            g\left(
                \gamma'(t),
                \gamma'(t)
            \right)^\frac{1}{2}
            dt
        <
        s_1,\quad 
        \gamma(t) \neq x
        \,
        \forall_{t > 0}.
    \end{align}
    Define
    \begin{equation}
        \delta 
    :=
        \inf
        \left\{
            t \in [0,1]:
            \gamma(t) \notin \exp_x\left(B_{s_2}\left(0_x\right)\right)
        \right\} \leq 1.
    \end{equation}
    By the definition of $\delta$ and the openness of $\exp_x\left(B_{s_2}(0_x)\right)$, we have
    \begin{align}
        \gamma([0,\delta[)
    \subset&
        \exp_x\left(B_{s_2}\left(0_x\right)\right)\\
        \gamma(\delta)
    \in&
        \partial \exp_x\left(B_{s_2}\left(0_x\right)\right),
        \label{eq:finalPointCurveInBoundary}
    \end{align}    
    By Lemma \ref{lem:injectivityLemma}, we can define $\hat{\gamma}:[0,\delta[ \rightarrow B_{s_2}(0_x)$ given by $\hat{\gamma}(t) = \exp_x^{-1}\left(\gamma(t)\right)$. By Lemma \ref{lem:homotopyContinuation}, this lift extends continuously to a curve $\hat{\gamma} :[0,\delta] \rightarrow \overline{B_{s_2}(0_x)}$.

    Since $\exp_x : {B}_s(0_x) \rightarrow B_s(x)$ is an open map for all $0\leq s<r$, we must have $\hat{\gamma}(\delta) \in \partial B_{s_2}(0_x)$, otherwise $\exp_x\left(\hat{\gamma}(\delta)\right) = \gamma(\delta) \in \text{int}\left(\exp_x\left(B_{s_2}\left(0_x\right)\right)\right)$, in contradiction with \eqref{eq:finalPointCurveInBoundary}.

    Given the decomposition of the path given by $\hat{\gamma}'(t) = \hat{\gamma}'_r(t) \partial r + \sum_{i=1}^{m-1}\hat{\gamma}'_{\theta_i}(t) \partial\theta_i$ for $t>0$, where $\partial r$ and $\partial \theta_i$ are orthogonal. Thus we have by Gauss' Lemma
    \begin{equation}
        g\left(
            \left(d\exp_x\right)\hat{\gamma}'(t),
            \left(d\exp_x\right)\hat{\gamma}'(t)
        \right)^\frac{1}{2}
        \geq
        \left|
            \hat{\gamma}'_r(t)
        \right|.
    \end{equation}  
    Since $\hat{\gamma}(0) = 0_x$ (only for $t=0$) and $\hat{\gamma}(\delta) \in \partial B_{s_2}(0_x)$ we must have
    \begin{align}
        s_2
    &\leq
        \int_0^\delta
            |\hat{\gamma}_r'(t)|
            dt\\
    &\leq
        \int_0^\delta
        g\left(
            \left(d\exp_x\right)\hat{\gamma}'(t),
            \left(d\exp_x\right)\hat{\gamma}'(t)
        \right)^\frac{1}{2}
        dt\\
    &=
        \int_0^\delta
        g\left(\gamma'(t), \gamma'(t)\right)^\frac{1}{2}
        dt\\
    &\leq
        \int_0^1
        g\left(\gamma'(t), \gamma'(t)\right)^\frac{1}{2}
        dt
    <
        s_1.
    \end{align}
    This contradicts $s_1<s_2$ and proves the lemma.
\end{proof}

\begin{proof}[Proof of Proposition \ref{prop:injectivityRadiusProp}]
    Choose constants $r_0>0$ and $C_{\alpha,\inj}>0$ so that Lemmas \ref{lem:exponentialMapDerivative}, \ref{lem:injectivityLemma} and \ref{lem:surjectivityLemma} are all satisfied. Fix $x \in N \setminus S$ and define $R_x = \min\left\{r_0, \left(
    \frac{\overline{\dist}_g\left(x,S\right)}{C_{\alpha,\inj}}
    \right)^{\max\{1,\frac{\alpha}{2}\}}\right\}$. For every $0<r< R_x$, we then have $\exp_x : B_r(0_x) \rightarrow B_r(x)$ is a bijective local diffeomorphism. Thus it is a diffeomorphism. The definition of injectivity radius gives $\inj(x) \geq R_x$.
\end{proof}

\subsection{Sobolev space on $N\setminus S$}

In this section, given a tuple $(N,S,g,\{\psi_i\}_{i=1}^l)$, we study some properties of the Sobolev spaces $W^{1,p}\left(N \setminus S\right)$ defined in Definition \ref{def:sobolevSpace}. In fact, we show a Sobolev inequality for the manifolds given by \ref{def:riemannianManifoldWithSingularitySet}. In much of the literature, such inequalities are proved under the assumption that the manifold is complete and satisfies a lower bound on the Ricci curvature, see, for example, \cite{Hebey1996}. Our space is not complete, and our Ricci curvature can blow up as we get close to the singularity $S$, so we will need another approach. We will try to obtain this by using the Sobolev inequality on Euclidean sets for the parametrizations. Sobolev inequalities on open sets of Euclidean hold, for example, when the set has $C^1$ or Lipschitz boundary (see \cite{EvansPDEs}). However we are interested in obtaining Sobolev bounds on sets of the form $\mathbb{B}^m_1(0) \setminus \psi_i^{-1}\left(S\right)$. For this we will need the notion of capacity of a set, which will allow us to ignore $\psi_i^{-1}\left(S\right)$ as long as this capacity is zero. The notion of capacity is studied in \cite{Turesson2000}.

\begin{definition}(\cite[Definition 2.35]{Turesson2000})
    For $1\leq p < \infty$, and a set $K \subset \mathbb{R}^m$, define the Sobolev $p$-capacity of $K$ by
    \begin{equation}
        C_{p}(K)
    :=
        \inf_{u \in W(K)}
        \int_{\mathbb{R}^m}
        \left(
            |u(x)|^p
            +
            |\nabla u(x)|^p
        \right)
        d\mathcal{L}^m(x)
    \end{equation}
    where the infimum is taken over the set
    \begin{equation}
        W(K)
    :=
        \{u \in W^{1,p}\left(\mathbb{R}^m\right): u|_U = 1, \text{ where } U \text{ is an open set containing } K\}.
    \end{equation}
\end{definition}

\vspace{10mm}

\begin{theorem}(\cite[Theorem 2.44]{Turesson2000})
\label{thm:sobolevSpacesCapacity}
    Let $\Omega \subset \mathbb{R}^m$ be an open subset, and let $K\subset \Omega$ is a relatively closed set in $\Omega$. Then if $C_p(K) = 0$ then the restriction map
    \begin{equation}
        W^{1,p}(\Omega)
    \rightarrow
        W^{1,p}(\Omega \setminus K)
    \end{equation}
    is an isometric isomorphism.
\end{theorem}

\begin{lemma}
\label{lem:capacityLemmaMinkowski}Let  $k \geq 0$ and $K \subset \mathbb{R}^m$ be a compact with finite upper $(m,k)$-dimensional Minkowski content, $\mathcal{M}^{(m,k,+)}_{\mathfrak{E}}(K)<\infty$. If $p \geq 1, k \geq 0$ and
    \begin{equation}
        m-k-p > 0,
    \end{equation}
    then $C_p(K) = 0$.
\end{lemma}
\begin{proof}

Given $\epsilon>0$, define $u_\epsilon$ by
\begin{equation}
    u_\epsilon(x)
    :=
    \left(\frac{\epsilon - \dist_\euc\left(x,T_{\epsilon,\dist_\euc}\left(K\right)\right)}
          {\epsilon}
    \right)_+.
\end{equation}
Since $K$ is compact, the function $u_\epsilon$ has compact support and is bounded by $1$, thus, $u_\epsilon \in L^p\left(\mathbb{R}^m\right)$.
 The distance function and
the positive-part truncation are globally $1$-Lipschitz, hence $u_\epsilon$ is globally $\frac{1}{\epsilon}$-Lipschitz and $\left|\nabla u_\epsilon(x)\right|\leq \frac{1}{\epsilon}$ for $\mathcal{L}^m$-a.e. $x \in \mathbb{R}^m$. Furthermore, its support is contained in $\overline{T_{2\epsilon, \dist_\euc}(K)}$. Thus we conclude $u_\epsilon \in W^{1,p}\left(\mathbb{R}^m\right)$. Because $u_\epsilon(x) = 1$ in the open set $T_{\epsilon,\dist_\euc}(K)$ which contains $K$, we have $u_\epsilon \in W(K)$.

Since $K$ has finite $(m,k)$-dimensional Minkowski content, there exists a constant $C_{\mathcal{M}} > 0$ such that for all $0<\epsilon < 1$ we have
    \begin{equation}
        \mathcal{L}^m(T_{\epsilon,\dist_\euc}(K))
    \leq
        C_{\mathcal{M}} \epsilon^{m-k}.
    \end{equation}
Thus we obtain for all $0<\epsilon < \frac{1}{3}$
    \begin{equation}
        \int_{\mathbb{R}^m}
            |u_\epsilon(x)|^p
            d\mathcal{L}^m(x)
        \leq
            \mathcal{L}^m
            (T_{3\epsilon, \dist_\euc}(K))
        \leq
            C_{\mathcal{M}}3^{m-k} \epsilon^{m-k},
    \end{equation}
    \begin{equation}
        \int_{\mathbb{R}^m}
            |\nabla u_\epsilon(x)|^p
            d\mathcal{L}^m(x)
        \leq
            \int_{T_{3\epsilon,\dist_\euc}(K)}
                \left(\frac{1}{\epsilon}\right)^p
                d\mathcal{L}^m(x)
        \leq
            C_{\mathcal{M}}3^{m-k} \epsilon^{m-k-p}.
    \end{equation}
    Therefore, for $0<\epsilon < \frac{1}{3}$, we have
    \begin{equation}
        C_p(K)
    \leq
        C_{\mathcal{M}}3^{m-k} \epsilon^{m-k}
    +
        C_{\mathcal{M}}3^{m-k} \epsilon^{m-k-p}.
    \end{equation}
    Since $m-k-p > 0$, letting $\epsilon \rightarrow 0$, we conclude $C_p(K) = 0$.

\end{proof}

\begin{lemma}
\label{lem:sobolevOnParametrizations}
    Let $(N,S,g,\{\psi_i\}_{i=1}^l)$ be a tuple satisfying \ref{def:riemannianManifoldWithSingularitySet}. Suppose there exist $k \geq 0$ and $p\geq 1$ such that, for every parametrization, we have $\mathcal{M}^{(m,k,+)}_{\mathfrak{E}}\left(\psi_i^{-1}\left(S\right)\right) < \infty$ and
    \begin{equation}
        m-k-p > 0,
    \end{equation}
    then there exists $C_{\mathcal{S},\euc} > 0$ such that for all $i \in \left\{1,\cdots, l\right\}$
    \begin{equation}
        \left\|u\circ \psi_i\right\|_{L^{p^*}\left(\mathbb{B}^m_1(0) \setminus \psi_i^{-1}(S), \mathcal{L}^m\right)}
    \leq
        C_{\mathcal{S}, \euc}
        \left(
        \left\| \left\|\nabla u\right\|_g \right\|_{L^{p}\left(\psi_i\left(\mathbb{B}^m_1\left(0\right)\right) \setminus S, \nu\right)}
        +
        \left\| u\right\|_{L^{p}\left(\psi_i\left(\mathbb{B}^m_1\left(0\right)\right) \setminus S, \nu\right)}
        \right),
    \end{equation}
    for $u \in W^{1,p}\left(N\setminus S, g\right) \cap C^\infty(N \setminus S)$, where $\frac{1}{p^*} = \frac{1}{p}-\frac{1}{m}$.
\end{lemma}
\begin{proof}
    Fix $u \in W^{1,p}(N\setminus S, g) \cap C^\infty(N \setminus S)$, and a parametrization $\psi_i \in \left\{\psi_j\right\}_{j=1}^l$. We notice that the bi-Lipschitz of the parametrizations \eqref{eq:boundAboveMetric} and the metric we have
    \begin{equation}
        \|\nabla (u\circ \psi_i)(x)\|
    \leq
        \sqrt{C_{UB}}
        \|\nabla u(\psi_i(x))\|_g.
    \end{equation}
    Let $\left\{e_1,\cdots,e_m\right\}$ be an orthonormal basis of $\mathbb{R}^m$, and consider the matrix for $y \in \mathbb{B}^m_2(0)$ $\left(g_{i,y}\right)_{jk} = g_{\psi_i(y)}\left(\left(d\psi_i\right)_ye_j, \left(d\psi_i\right)_ye_k\right)$. Hence, by inequality \eqref{eq:boundBellowMetric} 
    \begin{align}
        &\int_{\mathbb{B}^m_1(0) \setminus \psi_i^{-1}(S)}
        \|\nabla(u \circ \psi_i)(x)\|^pd\mathcal{L}^m(x)\\
    \leq&
        \frac{1}{\left(c_{LB}\right)^\frac{m}{2}}
        \int_{\mathbb{B}^m_1(0) \setminus \psi_i^{-1}(S)}
        \|\nabla(u \circ \psi_i)(x)\|^p
        \sqrt{\det g_{i, x}}d\mathcal{L}^m(x)\\
    \leq&
        \frac{C_{UB}^\frac{p}{2}}{\left(c_{LB}\right)^\frac{m}{2}}
        \int_{\mathbb{B}^m_1(0) \setminus \psi_i^{-1}(S)}
        \|\nabla u \left(\psi_i\left(x\right)\right)\|_g^p
        \sqrt{\det g_{i, x}}d\mathcal{L}^m(x)\\
    =&
        \frac{C_{UB}^\frac{p}{2}}{\left(c_{LB}\right)^\frac{m}{2}}
        \int_{\psi_i\left(\mathbb{B}^m_1(0)\right) \setminus S}
        \|\nabla u \left(x\right)\|^p_g d\nu(x) < \infty.
    \end{align}
    Similarly
    \begin{equation}
        \int_{\mathbb{B}^m_1(0) \setminus \psi_i^{-1}(S)}
        \left| u \circ \psi_i(x)\right|^pd\mathcal{L}^m(x)
    \leq
        \frac{1}{(c_{LB})^\frac{m}{2}}
        \int_{\psi_i\left(\mathbb{B}^m_1(0)\right) \setminus S}
        \left| u \left(x\right)\right|^p d\nu(x) < \infty,
    \end{equation}
    thus $u\circ \psi_i \in W^{1,p}\left(\mathbb{B}^m_1\left(0\right) \setminus \psi_i^{-1}\left(S\right), g_\euc\right)$.

    Notice that $K := \psi_i^{-1}(S) \cap \overline{\mathbb{B}^m_{\frac{3}{2}}\left(0\right)}$ is compact. Since $m-k-p > 0$ by Lemma \ref{lem:capacityLemmaMinkowski} we must have $C_p\left(K\right) = 0$. Thus by monotonicity of the capacity, this implies the set $E := \mathbb{B}^m_1\left(0\right) \cap K$ must also have capacity zero, $C_p(E) = 0$. $E$ is relatively closed in $\mathbb{B}^m_1(0)$, and so by Theorem \ref{thm:sobolevSpacesCapacity} we have that $W^{1,p}\left(\mathbb{B}^m_1\left(0\right)\setminus \psi_i^{-1}(S)\right) \cong W^{1,p}\left(\mathbb{B}^m_1\left(0\right)\setminus E\right) \cong W^{1,p}\left(\mathbb{B}^m_1\left(0\right)\right)$ are isometric spaces. Furthermore, the definition of capacity, $C_p(E)=0$ implies $\mathcal{L}^m(E)=0$, since every admissible function equals one on a neighborhood of $E$, so its $L^p$ term bounds the measure of that neighborhood from below. Thus $L^p\left(\mathbb{B}^m_1(0)\setminus E\right)\cong L^p\left(\mathbb{B}^m_1(0)\right)$ since $\mathcal{L}^m\left(E\right) = 0$. We know there exists  a constant $C_{\mathcal{S}, \mathbb{B}^m_1(0)}>0$ such that for all $f \in W^{1,p}\left(\mathbb{B}^m_1\left(0\right), g_\euc \right)$ we have
    \begin{align}
    \label{eq:sobolevIneqOnBall}
    \left\|
            f
        \right\|_{L^{p^*}\left(\mathbb{B}^m_1\left(0\right), \mathcal{L}^m\right)}
    \leq
        C_{\mathcal{S}, \mathbb{B}^m_1(0)}
        \left(
        \left\|
            \nabla f
        \right\|_{L^p\left(\mathbb{B}^m_1\left(0\right), \mathcal{L}^m\right)}
        +
        \left\|
            f
        \right\|_{L^p\left(\mathbb{B}^m_1\left(0\right), \mathcal{L}^m\right)}
        \right).
    \end{align}
    By the isometry of the Sobolev spaces, the Sobolev inequality \eqref{eq:sobolevIneqOnBall} is also true for $W^{1,p}(B_1^m(0)\setminus E)$. Hence, since $u\circ \psi_i \in W^{1,p}\left(\mathbb{B}^m_1\left(0\right) \setminus E\right)$, we have
    \begin{align}
        &\left\|
            u\circ \psi_i
        \right\|_{L^{p^*}\left(\mathbb{B}^m_1\left(0\right) \setminus \psi_i^{-1}\left(S\right)\right)}\\
    \leq&
        C_{\mathcal{S}, \mathbb{B}^m_1(0)}
        \left(
        \left\|
            \nabla\left( u\circ \psi_i\right)
        \right\|_{L^{p}\left(\mathbb{B}^m_1\left(0\right) \setminus \psi_i^{-1}\left(S\right), \mathcal{L}^m\right)}
        +
         \left\|
             u\circ \psi_i
        \right\|_{L^{p}\left(\mathbb{B}^m_1\left(0\right) \setminus \psi_i^{-1}\left(S\right), \mathcal{L}^m\right)}
        \right)\\
    \leq&
        \frac{C_{\mathcal{S}, \mathbb{B}^m_1(0)}
        }
        {\left(c_{LB}\right)^{\frac{m}{2p}}}
        \left(
        C_{UB}^\frac{1}{2}
        \left\|
            \left\|
                \nabla u
            \right\|_g
        \right\|_{L^{p}\left(\psi_i\left(\mathbb{B}^m_1\left(0\right) \right) \setminus S, \nu\right)}
        +
        \left\|
            u
        \right\|_{L^{p}\left(\psi_i\left(\mathbb{B}^m_1\left(0\right) \right) \setminus S, \nu\right)}
        \right).
    \end{align}     
    Thus the claim holds with, for example, $C_{\mathcal{S}, \euc} := \frac{C_{\mathcal{S}, \mathbb{B}^m_1(0)}
    }{\left(c_{LB}\right)^{\frac{m}{2p}}}\left(C_{UB}^\frac{1}{2} + 1\right)$.
\end{proof}

\begin{prop}
    Let $(N,S,g,\{\psi_i\}_{i=1}^l)$ be a tuple satisfying \ref{def:riemannianManifoldWithSingularitySet}. Suppose that there exist $k \geq 0$ and $p \geq 1$, for each parametrization we have $\mathcal{M}^{(m,k,+)}_{\mathfrak{E}}\left(\psi_i^{-1}\left(S\right)\right) < \infty$ and
    \begin{equation}
        m-k-p > 0.
    \end{equation}
    Then \eqref{eq:SobolevIneq} is satisfied by $(N\setminus S, g)$ with exponent $p$.
\end{prop}
\begin{proof}
    Given $i \in \{1,...,l\}$, and $\psi_i \in \left\{\psi_j\right\}_{j=1}^l$, define $\left(g_{i,\psi_i^{-1}(x)}\right)_{jk} = g_x\left(\left(d\psi_i\right)e_j, \left(d\psi_i\right)e_k\right)$. By equation \eqref{eq:boundAboveMetric} and Lemma \ref{lem:sobolevOnParametrizations} we obtain, for $u \in W^{1,p}\left(N\setminus S, g\right) \cap C^\infty(N\setminus S)$
    \begin{align}
        \left\|
            u
        \right\|_{L^{p^*}\left(\psi_i\left(\mathbb{B}^m_1\left(0\right)\right) \setminus S, \nu\right)}
    =&
    \left(
        \int_{\mathbb{B}^m_1\left(0\right) \setminus \psi_i^{-1}\left(S\right)}
            \left|
            u\circ \psi_i (x)
            \right|^{p^*}
            \sqrt{\det g_{i,x}}
            d\mathcal{L}^m(x)
    \right)^\frac{1}{p^*}\\
    \leq&
        \left(
            C_{UB}
        \right)^\frac{m}{2p^*}
        \left(
        \int_{\mathbb{B}^m_1\left(0\right) \setminus \psi_i^{-1}\left(S\right)}
            \left|
            u\circ \psi_i (x)
            \right|^{p^*}
            d\mathcal{L}^m(x)
    \right)^\frac{1}{p^*}\\
    \leq&
        \left(
            C_{UB}
        \right)^\frac{m}{2p^*}
        C_{\mathcal{S}, \euc}
        \left(
        \left\|
            \left\|
                \nabla u
            \right\|_g
            \right\|_{L^p\left(\psi_i\left(\mathbb{B}^m_1\left(0\right)\right)\setminus S, \nu\right)}
        +
            \left\|
                u
            \right\|_{L^p\left(\psi_i\left(\mathbb{B}^m_1\left(0\right)\right)\setminus S, \nu\right)}
        \right).
    \end{align}
    Thus, by \eqref{eq:NcoveredFiniteParametrizations}, we have
    \begin{align}
        \left\|
            u
        \right\|_{L^{p^*}\left(N \setminus S, \nu\right)}
    &\leq
        \sum_{i=1}^l
        \left\|
            u
        \right\|_{L^{p^*}\left(\psi_i\left(\mathbb{B}^m_1\left(0\right)\right) \setminus S, \nu\right)}\\
    &\leq
        \sum_{i=1}^l
        \left(
            C_{UB}
        \right)^\frac{m}{2p^*}
        C_{\mathcal{S}, \euc}
        \left(
        \left\|
            \left\|
                \nabla u
            \right\|_g
            \right\|_{L^p\left(\psi_i\left(\mathbb{B}^m_1\left(0\right)\right)\setminus S, \nu\right)}
        +
            \left\|
                u
            \right\|_{L^p\left(\psi_i\left(\mathbb{B}^m_1\left(0\right)\right)\setminus S, \nu\right)}
        \right)\\
    &\leq
        l \left(
            C_{UB}
        \right)^\frac{m}{2p^*}
        C_{\mathcal{S}, \euc}
        \left(
        \left\|
            \left\|
                \nabla u
            \right\|_g
            \right\|_{L^p\left(N \setminus S, \nu\right)}
        +
            \left\|
                u
            \right\|_{L^p\left(N \setminus S, \nu\right)}
        \right).
        \label{eq:sobolevIneForSmooth}
    \end{align}
Since this inequality is true for all $u \in W^{1,p}\left(N\setminus S, g\right) \cap C^\infty(N\setminus S)$, the result follows by density of this space in $W^{1,p}\left(N\setminus S, g\right) $. Given $u \in W^{1,p}\left(N\setminus S, g\right)$, let $u_n \in W^{1,p}(N\setminus S, g) \cap C^\infty(N\setminus S)$ be a sequence converging to $u$ in this norm. Inequality \eqref{eq:sobolevIneForSmooth} implies this sequence is Cauchy in $L^{p^*}\left(N\setminus S, \nu\right)$. This shows that $u_n$ converges in $L^{p^*}(N\setminus S, \nu)$ to $u \in L^{p^*}(N\setminus S,\nu)$, by uniqueness of limits, since it converged in $L^p(N\setminus S,\nu)$. Passing to the limit in the smooth estimate proves \eqref{eq:SobolevIneq} for $u \in W^{1,p}(N\setminus S, g)$.
\end{proof}

\section{Example for non-convergence}
    In this section, we construct an example for which the AMV Laplacian does not converge. We prove the following.

    \begin{prop}
    \label{prop:nonConvergenceResult}
        Let $z \in \mathbb{S}^2$. There exists a metric $\left( \mathbb{S}^2 \setminus \{z\}, g_\psi\right)$ and parametrizations $\left\{\Phi_i:\mathbb{B}^2_2(0) \rightarrow \mathbb{S}^2\right\}_{i=1}^2$ such that, $\left(\mathbb{S}^2,\{z\}, g_\psi, 
        \left\{\Phi_i\right\}_{i=1}^2\right)$ satisfies Definition \ref{def:riemannianManifoldWithSingularitySet}, and the curvature bound \ref{def:suppremumOfCurvatures} satisfies
        \begin{equation}
        \label{eq:curvatureBoundForExample}
            K(y) \leq \frac{C}{\overline{\dist}_{g_\psi}(y,z)^2}
        \end{equation}
        for some $C>0$. Furthermore, there exists $h : \mathbb{S}^2 \rightarrow \mathbb{R}$ such that $h \in C^\infty\left(\mathbb{S}^2\right)$, $\Delta_{g_\psi} h \in L^2\left(\mathbb{S}^2, \vol_{g_\psi}\right)$, $h$ is Lipschitz, $h \notin W^{2,2}\left(\mathbb{S}^2\setminus\{z\}, g_\psi\right)$ and
        \begin{equation}
            \Delta_r h \not \rightarrow_{L^p\left(\mathbb{S}^2, \vol_{g_\psi}\right)} C_m\Delta_{g_\psi} h
        \end{equation}
        for all $p \in [2,\infty[$, and for $p \in ]2,\infty[$ we have
        \begin{equation}
        \left\|
            \Delta_r h
        \right\|_{L^p\left(\mathbb{S}^2, \vol_{g_\psi}\right)}
        \rightarrow
        \infty.
        \end{equation}
    \end{prop}

    \begin{rem}
        In view of Theorem \ref{thm:convergenceOfAMV}, the completion of $\left(\mathbb{S}^2\setminus \left\{z\right\}, g_\psi\right)$ has a singularity set of dimension $0$ (a point), and so $k = 0$. Moreover, the curvature satisfies \eqref{eq:curvatureBoundForExample}, and so $\alpha = 2$. Thus, for all $p \geq 2$ we have
        \begin{equation}
            m-k-
            p
            \max
            \left\{1,\frac{\alpha}{2}
            \right\} \leq 0.
        \end{equation}
        Furthermore, although $\Delta_{g_\psi} h \in L^2\left(\mathbb{S}^2, \vol_{g_\psi}\right)$, we have $h \notin W^{2,2}\left(\mathbb{S}^2 \setminus \{z\},g_\psi\right)$. Thus, two hypotheses of Theorem \ref{thm:convergenceOfAMV} fail in this example, showing that they cannot both be omitted if convergence is to hold.
    \end{rem}

    \subsection{Calculation for discontinuous metric}

    Fix $\eta > 0$ and define the function $G_\eta : \mathbb{R}^2 \rightarrow \mathbb{R}$ given by
    \begin{equation}
        G_\eta(x,y)
    :=
        \begin{cases}
          1  & x \leq 0 \\
          e^\eta & x > 0
        \end{cases}.
    \end{equation}
    We define the following discontinuous metric tensor on $\mathbb{R}^2 \setminus \{(0,0)\}$:
    \begin{equation}
    \label{eq:discontinuousMetricDef}
        (g_\eta)_{(x,y)}
    :=
        G_\eta(x,y)g_\euc.
    \end{equation}
    This defines a distance on $\mathbb{B}^2_2(0,0)\setminus \{(0,0)\}$ given by
\begin{align}
    \dist_\eta\left(
        (x_0,y_0),
        (x_1,y_1)
    \right)
:=
    \inf&\Bigg\{
        \int_0^1 g_\eta(\gamma'(t),\gamma'(t))^\frac{1}{2}dt
    : \gamma(0) = (x_0,y_0), \gamma(1) = (x_1,y_1),\\
&\quad\quad\quad
    \gamma:[0,1] \rightarrow \mathbb{B}^2_2(0) \setminus \{(0,0)\}, \text { and } \gamma \text{ is piecewise } C^1
    \Bigg\}.
\end{align}
The completion of the metric space $\left(\mathbb{B}_2^2(0,0)\setminus\{(0,0)\}, \dist_\eta\right)$ is
$\left(\overline{\mathbb{B}^2_2(0,0)}, \overline{\dist}_\eta\right)$. We can also define a measure on $\mathbb{B}^2_2(0,0)$ given by
\begin{equation}
    d\mu_\eta = G_\eta d\mathcal{L}^2.
\end{equation}

\begin{lemma}
\label{lem:ballsOfDiscontinuousMetric}
    For $\eta>1$ and $0<r<\frac{1}{2}$, the balls $B^\eta_r(0,0)$ given by the distance $\overline{\dist}_\eta$ satisfy
    \begin{align}
    \label{eq:equality1discontinuousMetric}
        &B^\eta_r(0,0)
        \cap
        \left\{(x,y) \in \mathbb{R}^2: x\leq 0\right\}
    =
        \mathbb{B}^2_r(0,0)
        \cap
        \left\{(x,y) \in \mathbb{R}^2:
        x \leq 0\right\},\\
    \label{eq:equality2DiscontinuousMetric}
    &B^\eta_r(0,0)
    \cap\left\{(x,y)\in \mathbb{R}^2 : x > 0, |y| \geq \frac{x}{\sqrt{e^\eta-1}}\right\}\\
    =&
    \left\{(x,y)\in \mathbb{R}^2 : x > 0,
    |y| \geq \frac{x}{\sqrt{e^\eta-1}},
        |y|< r - \sqrt{e^\eta - 1}x
    \right\},\\
    \label{eq:equality3DiscontinuousMetric}
    &B^\eta_r(0,0)
    \cap\left\{(x,y)\in \mathbb{R}^2 : x > 0, |y| < \frac{x}{\sqrt{e^\eta-1}}\right\}\\
    =&
    \mathbb{B}^2_{r/e^{\eta/2}}(0,0)
    \cap
    \left\{(x,y)\in \mathbb{R}^2 : x > 0,  |y| < \frac{x}{\sqrt{e^\eta-1}}\right\}.
    \end{align}
    In particular
    \begin{equation}
    \label{eq:inclusionOfBallsDiscMetric}
        B^\eta_r(0,0)
    \cap
        \left\{(x,y) \in \mathbb{R}^2: x > 0\right\}
    \subset
        \left\{(x,y) \in \mathbb{R}^2: x > 0, |y|<r -\sqrt{e^\eta - 1}x
        \right\}.
    \end{equation}
\end{lemma}
\begin{proof}
    Given $(x,y) \in \mathbb{B}^2_2(0,0)$, we start by noticing that $\overline{\dist}_\eta\left((x,y),(0,0)\right)$ can be computed by the length of curves $\gamma:[0,1] \rightarrow \mathbb{B}^2_2(0,0)$ such that $\gamma([0,1[) \subset \mathbb{B}^2_2(0,0)$ and $\gamma(1) = (0,0)$. Notice that for any admissible curve $\gamma:[0,1] \rightarrow \mathbb{B}^2_2(0,0)$, we have
    \begin{equation}
        \int_0^1 g_\eta\left(\gamma'(t),\gamma'(t)\right)^\frac{1}{2}dt
    \geq
        \int_0^1 g_\euc\left(\gamma'(t),\gamma'(t)\right)^\frac{1}{2}dt
    =
        l_\euc(\gamma).
    \end{equation}
    Thus, for all $r > 0$
    \begin{equation}
        B^\eta_r(0,0)
    \subset
        \mathbb{B}^2_r(0,0).
    \end{equation}
    Furthermore, if $(x,y) \in \mathbb{B}^2_2(0,0)$ and $x \leq 0$, then the curve $\gamma(t) = t(x,y)$ satisfies
    \begin{equation}
        \overline{\dist}_\eta\left(
            \left(
                x,y
            \right),
            \left(
                0,0
            \right)
        \right)
        \leq
        \int_0^1 g_\eta\left(\gamma'(t),\gamma'(t)\right)^\frac{1}{2}dt
    =
        \int_0^1 g_\euc\left(\gamma'(t),\gamma'(t)\right)^\frac{1}{2}dt
    =
        \|(x,y)\|.
    \end{equation}
    This concludes the first equality \eqref{eq:equality1discontinuousMetric}.

    Now let $(x,y) \in \mathbb{B}^2_r(0,0)$ such that $x > 0$. By symmetry with respect to the $x$-axis of the metric, we have in general
    \begin{equation}
        \overline{\dist}_\eta
        \left(
            (x,y)
        ,
            (0,0)
        \right)
    =
        \overline{\dist}_\eta
        \left(
            (x,-y)
        ,
            (0,0)
        \right),
    \end{equation}
    so we will assume $y \geq 0$.
    To calculate the distance between $(x,y)$ and $(0,0)$, consider a curve $\gamma:[0,1] \rightarrow \mathbb{B}^2_2(0,0)$ joining $(x,y)$ to $(0,0)$. Using the map given by $H(x,y) = (x,|y|)$ we have that $H\circ \gamma$ is absolutely continuous (its length can be approximated by piecewise $C^1$ curves). Also
    \begin{equation}
        \int_0^1 g_\eta\left(\gamma'(t),\gamma'(t)\right)^\frac{1}{2}dt
    =
        \int_0^1 g_\eta\left((H\circ\gamma)'(t),(H\circ\gamma)'(t)\right)^\frac{1}{2}dt,
    \end{equation}
    and hence $H\circ \gamma$ has the same length as $\gamma$. Thus, we may assume that the curve $\gamma(t)=(\gamma_x(t),\gamma_y(t))$ satisfies $\gamma_y(t) \geq 0$. Since $\gamma(1) = (0,0)$, there exists a first time $t_0 \in ]0,1]$ such that $\gamma_x(t_0) = 0$. Then we have that the curve $\sigma_1 : [0,t_0] \rightarrow \mathbb{B}^2_2(0,0)$ given by $\sigma_1(t) = \frac{t_0-t}{t_0}(x,y) + \frac{t}{t_0}\gamma(t_0)$ satisfies
    \begin{align}
        \int_0^{t_0}
        g_\eta\left(\gamma'(t),\gamma'(t)\right)^\frac{1}{2}dt
    =&
        \int_0^{t_0} e^{\frac{\eta}{2}} g_\euc\left(\gamma'(t),\gamma'(t)\right)^\frac{1}{2}dt\\
    \geq&
        \int_0^{t_0} e^{\frac{\eta}{2}} g_\euc\left(\sigma_1'(t),\sigma_1'(t)\right)^\frac{1}{2}dt\\
    =&
        \int_0^{t_0} g_\eta\left(\sigma_1'(t),\sigma_1'(t)\right)^\frac{1}{2}dt.
    \end{align}
    Since the length of $\sigma_1([0,t_0])$ is no greater than that of $\gamma([0,t_0])$, we may assume that $\left.\gamma\right|_{[0,t_0]} = \sigma_1$; that is, $\left.\gamma\right|_{[0,t_0]}$ is a straight line from $(x,y)$ to $\gamma(t_0)$. If $t_0<1$, then $\gamma(t_0) = (0,\gamma_y(t_0))$. Define $\sigma_2:[t_0,1] \rightarrow \mathbb{B}^2_2(0)$ by $\sigma_2(t) = \left(1-\frac{t-t_0}{1-t_0}\right)\gamma(t_0)$. Then
    \begin{align}
        \int_{t_0}^{1}
        g_\eta\left(\gamma'(t),\gamma'(t)\right)^\frac{1}{2}dt
    \geq&
        \int_{t_0}^{1} g_\euc\left(\gamma'(t),\gamma'(t)\right)^\frac{1}{2}dt\\
    \geq&
        \int_{t_0}^{1} g_\euc\left(\sigma_2'(t),\sigma_2'(t)\right)^\frac{1}{2}dt\\
    =&
        \int_{t_0}^{1} g_\eta\left(\sigma_2'(t),\sigma_2'(t)\right)^\frac{1}{2}dt.
    \end{align}
    Thus, we may assume that $\gamma|_{[t_0,1]} = \sigma_2$, a straight line from $\gamma(t_0)$ to $(0,0)$. Therefore, to calculate the distance from $(x,y)$ to $(0,0)$, it suffices to minimize over curves that first follow a straight line from $(x,y)$ to $(0,s)$ for some $s \geq 0$ and then follow the vertical segment from $(0,s)$ to $(0,0)$. We must therefore solve the minimization problem
    \begin{align}
        \overline{\dist}_\eta\left(
            (x,y),
            (0,0)
        \right)
    =
        \min_{s\geq 0}
            e^{\eta/2}
            \sqrt{
                x^2 + (s-y)^2
            }
            +
            s
        =
            \min_{s\geq 0}
            h_\eta(s).
    \end{align}
    The derivatives of $h_\eta$ are
    \begin{align}
        h_\eta'(s)
    =
        e^{\eta/2}\frac{(s-y)}{\sqrt{x^2+(s-y)^2}}
    +
        1,\quad\quad
    \label{eq:secondDerivativeBound}
        h_\eta''(s)
    =
        e^{\eta/2}
        \frac{x^2}{\left(x^2+(s-y)^2\right)^\frac{3}{2}}
    >
        0.
    \end{align}
    We calculate the critical points:
    \begin{align}
        h_\eta'(s)
    =
        0&\Leftrightarrow\\
        e^{\eta/2}\frac{(s-y)}{\sqrt{x^2+(s-y)^2}}
    +
        1
    =
        0&\Rightarrow\\
        e^{\eta}(s-y)^2
    =
        x^2+(s-y)^2
        &\Leftrightarrow\\
        s
    =
        y
    \pm
        \frac{x}{\sqrt{e^\eta - 1}}.
    \end{align}
    In fact, we must have $s_{(x,y)} = y - \frac{x}{\sqrt{e^\eta-1}}$ so that $h_\eta'(s_{(x,y)}) = 0$. For $s_{(x,y)} \geq 0$, we must have
    \begin{equation}
    \label{eq:conditionForPositiveS}
        y \geq \frac{x}{\sqrt{e^\eta-1}}.
    \end{equation}
    Suppose that $(x,y)$ satisfies \eqref{eq:conditionForPositiveS}. By \eqref{eq:secondDerivativeBound}, $s_{(x,y)}$ is then the minimizer, and hence
    \begin{align}
        \overline{\dist}_\eta\left(\left(x,y\right),\left(0,0\right)\right)
    =&
        h_\eta(s_{(x,y)})\\
    =&
        e^{\eta/2}
        \sqrt{x^2 + \frac{x^2}{e^\eta-1}}
        +
        y-\frac{x}{\sqrt{e^\eta - 1}}\\
    =&
        y
    +
         \sqrt{e^{\eta}-1}x.
    \end{align}
    On the other hand, if $y < \frac{x}{\sqrt{e^\eta-1}}$, then $s_{(x,y)} < 0$. Since $h_\eta'(s_{(x,y)}) = 0$ and \eqref{eq:secondDerivativeBound} shows that $h_\eta'$ is strictly increasing, we have $h_\eta'(s) > 0$ for all $s \geq 0$. Thus, $\min_{s\geq 0} h_\eta(s) = h_\eta(0)$, and
    \begin{equation}
        \overline{\dist}_\eta((x,y),(0,0))
    =
        e^{\eta/2}
        \sqrt{x^2+y^2},
    \end{equation}
    which implies equation \eqref{eq:equality3DiscontinuousMetric}. The inclusion \eqref{eq:inclusionOfBallsDiscMetric} can be deduced from the previous equalities.
\end{proof}

Let the function $f: \mathbb{B}^2_2(0,0) \rightarrow \mathbb{R}$ be given by $f(x,y) = x$.

\begin{lemma}
\label{lem:integralOfXDiscontMetric}
    There exist $\delta >0$ and $\eta > 0$ such that for all $0<r < \frac{1}{2}$
    \begin{equation}
        \int_{B^\eta_r(0,0)}
            f
            d\mu_\eta
        \leq
            -\delta r^3.
    \end{equation}
\end{lemma}
\begin{proof}
    Using Lemma \ref{lem:ballsOfDiscontinuousMetric}, for $0<r<\frac{1}{2}$ we have
    \begin{align}
        \int_{B_{r}^\eta(0,0) \cap \{x \leq 0\}}
            f
            d\mu_\eta
    =&
        \int_{\mathbb{B}_{r}^2(0,0) \cap \{x \leq 0\}}
            x
            d\mathcal{L}^2(x,y)\\
    =&
        \int_{\frac{\pi}{2}}^{\frac{3\pi}{2}}
        \int_0^r
            \rho^2
            \cos(\theta)
            d\rho
        d\theta
    =
        -\frac{2}{3}r^3.
    \end{align}
    Also using \eqref{eq:inclusionOfBallsDiscMetric}, and the fact that $f(x,y) \geq 0$ for $x >0$, we obtain
    \begin{align}
        \int_{B_{r}^\eta(0,0) \cap \{x > 0\}}
            f
            d\mu_\eta
    &\leq
        \int_{\left\{0<x<\frac{r}{\sqrt{e^\eta-1}}, |y|<-\sqrt{e^\eta - 1}x + r
        \right\}}
            x
            e^\eta
            d\mathcal{L}^2(x,y)\\
    &=
        \int_{\left\{x > 0, -r +\sqrt{e^{\eta}-1}x < y <-\sqrt{e^\eta - 1}x + r
        \right\}}
            x
            e^\eta
            d\mathcal{L}^2(x,y)\\
    &=
        e^\eta
        \int_0^{\frac{r}{\sqrt{e^\eta-1}}}
        \int_{-r + \sqrt{e^\eta - 1}x}^{r-\sqrt{e^\eta - 1}x}
            x
            dy
            dx\\
    &=
        e^{\eta}
        \int_0^{\frac{r}{\sqrt{e^\eta-1}}}
        x
        \left(
            2r - 2\sqrt{e^\eta - 1}x
        \right)
        dx\\
    &=
        e^\eta
        \left(
            r
            \left(
                \frac{r}{\sqrt{e^\eta - 1}}
            \right)^2
            -
            \frac{2}{3}
            \sqrt{e^\eta - 1}
            \left(
                \frac{r}{\sqrt{e^\eta - 1}}
            \right)^3
        \right)\\
    &=
        r^3
        \frac{1}{3}
        \frac{e^\eta}
        {e^\eta -1 }.
    \end{align}
    Choose $\eta > 0$ so that $\frac{e^\eta}{e^\eta - 1} < \frac{3}{2}$. Taking $\delta = \frac{1}{6}$, we obtain
    \begin{align}
        \int_{B^\eta_r(0,0)}
        f
        d\mu_\eta
    =&
        \int_{B^\eta_r(0,0) \cap \{x \leq 0\}}
        f
        d\mu_\eta
    +
        \int_{B^\eta_r(0,0) \cap \{x > 0\}}
        f
        d\mu_\eta\\
    \leq&
        -\frac{2}{3}r^3
    +
        r^3\frac{1}{3}\frac{e^\eta}{e^\eta - 1}
    <
        -\delta r^3,
    \end{align}
    finishing the proof.
\end{proof}

\subsection{Metrics in $\mathbb{S}^2\setminus \{z\}$}

Let $z \in \mathbb{S}^2$. Now we introduce a family of metrics on $\mathbb{S}^2\setminus \{z\}$, where $z \in \mathbb{S}^2$. Let $\varphi\in C^\infty([-\pi,\pi])$ be a function such that $\partial_x^k\varphi(-\pi) = \partial_x^k\varphi(\pi)$ for all integers $k \geq 0$. Consider a metric $(\mathbb{S}^2\setminus \{z\}, g_\varphi)$, and a parametrization $\Phi:\mathbb{B}^2_2(0) \rightarrow \mathbb{S}^2$ with $\Phi(0) = z$ such that
\begin{equation}
\label{eq:metricOnSphere}
\Phi^*g_\varphi|_{\mathbb{B}^2_2(0)\setminus \{0\}}
=
    e^{\varphi(\theta)}
    g_\euc|_{\mathbb{B}^2_2(0)\setminus \{0\}},
\end{equation}
Such metrics can be constructed. Let $\Phi_1:= \Phi$. We can find another parametrization, $\Phi_2 : \mathbb{B}^2_2(0) \rightarrow \mathbb{S}^2$ such that 
\begin{equation}
\label{eq:twoParametrizationsCover}
    \mathbb{S}^2 \setminus \Phi_1\left(\mathbb{B}^2_1(0)\right) \subset \Phi_2\left(\mathbb{B}^2_1(0)\right)
\end{equation}

Similarly, this metric defines the Riemannian volume measure $\mu_\varphi = \vol_{g_\varphi}$, which satisfies
\begin{equation}
\label{eq:metricOnSphere2}
    \Phi^* d\mu_\varphi = e^{\varphi(\theta)}d\mathcal{L}^2.
\end{equation}
Since $\varphi$ is smooth on a compact interval, it is bounded, and the metric satisfies \eqref{eq:boundBellowMetric} and \eqref{eq:boundAboveMetric}. We denote the induced distance on $\mathbb{S}^2 \setminus \{z\}$ by $\dist_{\varphi}$.

\begin{lemma}
    Let $\varphi \in C^\infty\left(\left[-\pi,\pi\right]\right)$ such that for every integer $k \geq 0$, $\partial_x^k\varphi(-\pi) = \partial_x^k\varphi(\pi)$. Then $ \left(\mathbb{S}^2, \{z\}, g_\varphi, \left\{\Phi_i\right\}_{i=1}^2\right)$ satisfies definition \ref{def:riemannianManifoldWithSingularitySet}, and
    \begin{equation}
        K_\varphi
        \left(\Phi\left(x,y\right)\right) = \left|\frac{e^{-\varphi(\theta(x,y))}}{2r^2}\varphi''(\theta(x,y))\right|.
    \end{equation}
    In particular there exists $C>0$ depending on $\varphi$ such that
    \begin{equation}
    \label{eq:boundOnCurvature}
        K_\varphi
        \left(
            \Phi
            \left(
                x,y
            \right)
        \right)
    \leq
        \frac
        {C}
        {\left\|(x,y)\right\|^2}.
    \end{equation}
\end{lemma}
\begin{proof}
    There exists a parametrization $\Phi_2$ such that \eqref{eq:twoParametrizationsCover}. Furthermore by continuity of the metric  in $\mathbb{S}^2 \setminus \Phi\left(\mathbb{B}^2_1(0)\right)$ and compactness of this set, there exist constants $C_\Psi, \tilde{c}_\Psi > 0$ such that
    \begin{equation}
        \tilde{c}_\Psi
        \|v\|^2
    \leq
        \Psi^*g_{\varphi}(v,v)
    \leq
        C_\Psi \|v\|^2,
    \end{equation}
    using compactness of $\mathbb{S}^2 \setminus \Phi\left(\mathbb{B}^2_1(0)\right)$ and continuity of the metric $g_\varphi$. On the other hand we have that
    \begin{equation}
    \label{eq:metricBoundsSphere}
            \min_{\theta \in ]-\pi,\pi]}
                e^{\varphi(\theta)}
                \|v\|^2
            \leq
            \Phi^*
                g_\varphi(v,v)
            \leq
            \max_{\theta \in  ]-\pi,\pi]}
                e^{\varphi(\theta)}
                \|v\|^2.
    \end{equation}
    Since $\varphi$ is continuous, these constants are finite and strictly positive. Furthermore, $\{z\}$ has finite upper $(2,0)$-Minkowski content. Therefore, $\left(\mathbb{S}^2, \{z\}, g_\varphi, \left\{\Phi_i\right\}_{i=1}^2\right)$ satisfies Definition \ref{def:riemannianManifoldWithSingularitySet}. By Theorem \ref{thm:riemannianManifoldWithSingularityThm}, the completion of $\left(\mathbb{S}^2 \setminus \{z\}, \dist_\varphi\right)$ is $\left(\mathbb{S}^2, \overline{\dist_\varphi}\right)$.

    In polar coordinates, the nonzero Christoffel symbols under the parametrization $\Phi$ are
    \begin{align}
        \Gamma_{r\theta}^r
    =
        \Gamma_{\theta r}^r
    =&
        \frac{1}{2}\varphi'(\theta),
    \quad\quad
        \Gamma_{\theta\theta}^r
    =
        -r,
    \quad\quad
        \Gamma_{rr}^\theta
    =
        -
        \frac{1}{2 r^2}
        \varphi'(\theta),\\
        \Gamma_{r\theta}^\theta
    =
        \Gamma_{\theta r}^\theta
    =&
        \frac{1}{r}, \quad\quad\quad\quad
    \Gamma_{\theta\theta}^\theta
    =
    \frac{1}{2}
    \varphi'(\theta),
    \quad
    \Gamma_{rr}^r = 0.
    \end{align}
    Since this is a two-dimensional metric, to calculate $K_\varphi\left(\Phi
            \left(
                x,y
            \right)\right)$ from \eqref{eq:supremmumOfTheSectionalCurvatures}, it suffices to calculate the Gaussian curvature $k_{g_\varphi}\left(\Phi
            \left(
                x,y
            \right)\right)$, because $K_\varphi\left(\Phi
            \left(
                x,y
            \right)\right) = \left|k_{g_\varphi}\left(\Phi
            \left(
                x,y
            \right)\right)\right|$. We have
    \begin{align}
        &(g_\varphi)_{\Phi
            \left(
                x,y
            \right)}
        \left(
            R(\partial \theta, \partial r, \partial \theta),
            \partial r
        \right)\\
    =&
        k_{g_\varphi}\left(\Phi
            \left(
                x,y
            \right)\right)
        \left(
            g_\varphi(\partial \theta, \partial \theta)
            g_\varphi(\partial r, \partial r)
        -
            g_\varphi(\partial \theta, \partial r)
            g_\varphi(\partial \theta, \partial r)
        \right)\\
    =&
        r^2
        e^{2\varphi(\theta(x,y))}
        k_{g_\varphi}\left(\Phi
            \left(
                x,y
            \right)\right).
    \end{align}

    Writing $R_{\Phi
            \left(
                x,y
            \right)}\left(\partial \theta, \partial r, \partial \theta\right) = R_{\Phi
            \left(
                x,y
            \right)}\left(\partial \theta, \partial r, \partial \theta\right)^r \partial r + R_{\Phi
            \left(
                x,y
            \right)}\left(\partial \theta, \partial r, \partial \theta\right)^\theta \partial \theta$, we obtain
    \begin{align}
        R_{\Phi
            \left(
                x,y
            \right)}\left(\partial \theta, \partial r, \partial \theta
        \right)^r
    =&
        \partial_r \Gamma_{\theta \theta}^r
       -
       \partial_\theta \Gamma_{\theta r}^r
       +
        \Gamma_{\theta \theta}^m
        \Gamma_{rm}^r
        -
        \Gamma_{\theta r}^m
        \Gamma_{\theta m}^r\\
    =&
        -1
    -
        \frac{1}{2}\varphi''(\theta)
    +
        \left(\frac{1}{2}\varphi'(\theta)\right)^2
    -
    \left(\left(\frac{1}{2}\varphi'(\theta)\right)^2 - 1\right)\\
    =&
        -\frac{1}{2}\varphi''(\theta).
    \end{align}
    Here, the repeated index $m$ is summed over $r$ and $\theta$. Since $\left(g_\varphi\right)_{\Phi
            \left(
                x,y
            \right)}\left(\partial \theta, \partial r\right) = 0$,
    \begin{equation}
        k_{g_\varphi}\left(\Phi
            \left(
                x,y
            \right)\right)
    =
        -\frac{e^{-\varphi(\theta(x,y))}}{2r^2}\varphi''(\theta(x,y)).
    \end{equation}
    This implies
    \begin{equation}
        K_{\varphi}\left(\Phi
            \left(
                x,y
            \right)\right)
    \leq
        \max_{\theta \in [-\pi,\pi]}
        \left|
        \frac{e^{-\varphi(\theta)} \varphi''(\theta)}{2}
        \right|
        \frac{1}{r^2},
    \end{equation}
    finishing the proof.
\end{proof}

Choose a function $h:\mathbb{S}^2 \rightarrow \mathbb{R}$ such that
\begin{equation}
\label{eq:definitionOfh}
    h\circ \Phi|_{\mathbb{B}^2_1(0)}(x,y) = x,
    \quad
    h|_{\mathbb{S}^2 \setminus \Phi\left(\mathbb{B}^2_{\frac{3}{2}}(0)\right)} = 0,
    \quad
    h \in C^\infty\left(\mathbb{S}^2 \setminus \{z\}\right).
\end{equation}

\begin{lemma}
\label{lem:hIsHarmonicAndNotSobolev}
    Let $\varphi \in C^\infty\left(\left[-\pi,\pi\right]\right)$ such that for all integers $k \geq 0$, $\partial_x^k\varphi(-\pi) = \partial_x^k\varphi(\pi)$. Then
    \begin{equation}
        \Delta_{g_\varphi} h|_{\Phi\left(\mathbb{B}^2_1(0,0)\right) \setminus \{z\}} = 0,
    \end{equation}
    and $\Delta_{g_\varphi} h \in L^p\left(\mathbb{S}^2 \setminus \{z\}, \mu_\varphi\right)$ for all $1 \leq p \leq \infty$, and $h \in W^{2,p}\left(\mathbb{S}^2\setminus \{z\}, g_\varphi\right)$ for $p < 2$. Furthermore, if $\varphi'(\theta_0) \neq 0$ for some $\theta_0$, then $h \notin W^{2,p}\left(\mathbb{S}^2 \setminus \{z\}, g_\varphi\right)$ for $p \geq 2$.
\end{lemma}
\begin{proof}
    First we notice that under a conformal change of metric, and using \eqref{eq:definitionOfh} we have for $w \in \Phi\left(\mathbb{B}^2_1\left(0,0\right)\right) \setminus \{z\}$
    \begin{equation}
        \Delta_{g_\varphi} h (w)
    =
        \frac{1}{e^{\varphi \left(\theta\left(\Phi^{-1}(w)\right)\right)}}
        \Delta_{\euc}
        \left(
            h\circ \Phi
        \right)
        \left(\Phi^{-1}(w)\right)
    =
        0.
    \end{equation}
    From this we can easily deduce that $\Delta_{g_\varphi} h \in L^p\left(\mathbb{S}^2 \setminus \{z\}, \mu_\varphi\right)$ for all $1 \leq p \leq \infty$, using that $h$ is $C^\infty$ in the compact set $\mathbb{S}^2 \setminus \Phi\left(\overline{\mathbb{B}^2_{\frac{1}{2}}\left(0,0\right)}\right)$. Now we calculate the Hessian of $h$. Again since $h$ is $C^\infty$ in the compact set $\mathbb{S}^2 \setminus \Phi\left(\overline{\mathbb{B}^2_{\frac{1}{2}}(0,0)}\right)$, here the Hessian is uniformly bounded. For points $w \in \Phi\left(\mathbb{B}^2_1(0,0)\right)$ let $r = \|\Phi^{-1}(w)\|$ and $\theta_w = \theta\left(\Phi^{-1}(w)\right)$. Then
    \begin{align}
    \label{eq:radial2DerivativesOfh}
        \nabla_{\partial r}
        \nabla_{\partial r}
        h(w)
    &=
        \partial_r \partial_r
        \left(h\circ \Phi\right)
    -
        (g_\varphi)_{w}
        \left(
            \nabla h,
            \nabla_{\partial r}
            \partial r
        \right)\\
    &=
        -\frac{1}{2r^2}
        \varphi'(\theta_w)
        r
        \sin(\theta_w)
    =
        -\frac{\sin(\theta_w)\varphi'(\theta_w)}
        {2 r}.
    \end{align}
    Similarly we have
    \begin{align}
        \nabla_{\partial \theta}
        \nabla_{\partial \theta}
        h(w)
    =&
        \partial_\theta
        \partial_\theta
        \left(h\circ\Phi\right)
    -
        (g_\varphi)_{w}(\nabla h, \nabla_{\partial \theta}\partial \theta)\\
    =&
        -r\cos(\theta_w)
    -
        \left(
        -\partial_r\left(
                h\circ \Phi
            \right)
            r
        +
            \partial_\theta
            \left(
                h\circ \Phi
            \right)
            \frac{\varphi'(\theta_w)}{2}
        \right)
    =
        \frac
        {r\sin(\theta_w) \varphi'(\theta_w)}
        {2},
    \end{align}
    and
    \begin{align}
        \nabla_{\partial \theta}
        \nabla_{\partial r}
        h(w)
    =
        \nabla_{\partial r}
        \nabla_{\partial \theta}
        h(w)
    &=
        \partial_r
        \partial_\theta
        \left(
            h\circ \Phi
        \right)
    -
        \left(
            \partial_r(h\circ \Phi)
            \frac{\varphi'(\theta_w)}{2}
        +
            \partial_\theta(h\circ \Phi)
            \frac{1}{r}
        \right)\\
    &=
        -\sin(\theta_w)
    -
    \left(
        \cos(\theta_w)\frac{\varphi'(\theta_w)}{2}
    -
        \sin(\theta_w)
    \right)\\
    &=
        -\frac{\cos(\theta_w)\varphi'(\theta_w)}{2}.
    \end{align}
    Now notice that the vectors $\frac{\partial r}{e^{\varphi(\theta_w)/2}}$ and $\frac{\partial \theta}{r e^{\varphi(\theta_w)/2} }$ are orthonormal with respect to the metric $g_\varphi$ in $T_w\mathbb{S}^2$. Thus, if $v \in T_{w}\mathbb{S}^2$, $v = a \frac{\partial r}{e^{\varphi(\theta_w)/2}} + b\frac{\partial \theta}{r e^{\varphi(\theta_w)/2} }$ and $g_\varphi(v,v) = 1$ then $a^2 + b^2 = 1$. Then we have
    \begin{align}
        &\sup_{v_1,v_2 \in T_{w}\mathbb{S}^2 : g_\varphi(v_1,v_1) = g_\varphi(v_2,v_2) = 1}
            \left|
            \nabla_{v_1}
            \nabla_{v_2}
            h(w)
            \right|\\
        =&
            \sup_{a_1,b_1,a_2,b_2 : a_1^2+b_1^2 = 1, a_2^2+b_2^2 = 1}
            \left|
            \nabla_{a_1 \frac{\partial r}{e^{\varphi(\theta_w)/2}} + b_1\frac{\partial \theta}{r e^{\varphi(\theta_w)/2} }}
            \nabla_{a_2 \frac{\partial r}{e^{\varphi(\theta_w)/2}} + b_2\frac{\partial \theta}{r e^{\varphi(\theta_w)/2} }}
            h(w)
            \right|\\
    \leq&
            \left|
                \frac{1}{e^{\varphi(\theta_w)}}
                \nabla_{\partial r}
                \nabla_{\partial r}
                h(w)
            \right|
            +
            2
            \left|
                \frac{1}{e^{\varphi(\theta_w)}r}
                \nabla_{\partial r}
                \nabla_{\partial \theta}
                h(w)
            \right|
            +
             \left|
                \frac{1}{e^{\varphi(\theta_w)}r^2}
                \nabla_{\partial \theta}
                \nabla_{\partial \theta}
                h(w)
            \right|\\
    \leq&
        \frac{1}{e^{\varphi(\theta_w)}}
        \left(
            \frac{\left| \sin(\theta_w)\varphi'(\theta_w)\right|}{2r}
        +
            2
            \frac{\left|\cos(\theta_w) \varphi'(\theta_w)\right|}{2r}
        +
            \frac{r\left|\sin(\theta_w) \varphi'(\theta_w)\right|}{2r^2}
        \right)\\
    \leq&
        \frac{4 \left|\varphi'(\theta_w)\right|}{e^{\varphi(\theta_w)}r}
    \leq
        \left(\max_{\theta \in ]-\pi,\pi]}\frac{4 \left|\varphi'(\theta)\right|}{e^{\varphi(\theta)}}\right)
        \frac{1}{r}.
    \end{align}
    Thus, taking $C_H = \left(\max_{\theta \in ]-\pi,\pi]}\frac{4 \left|\varphi'(\theta)\right|}{e^{\varphi(\theta)}}\right)$, we have
    \begin{equation}
    \label{eq:hessianBoundBy1OverDist}
        \|H h (w)\|_{g_\varphi} \leq \frac{C_H}{\dist_\euc\left(\Phi^{-1}(w), \left(0,0\right)\right)}.
    \end{equation}

    Equation \eqref{eq:hessianBoundBy1OverDist} and Lemma \ref{lem:1OverDistIsIntegrable} imply that $\|(H h)\circ \Phi\|_{g_\varphi} \in L^p\left(\mathbb{B}^2_1(0,0), \mathcal{L}^2\right)$ for all $p < 2$, since the upper $(2,0)$ Minkowski content of a point in this space is finite and this space has dimension $2$. Thus for any $p < 2$
    \begin{align}
        \int_{\Phi\left(\mathbb{B}^2_1(0,0)\right)}
            \left\|
                (Hh)
            \right\|_{g_\varphi}^p
            d\mu_\varphi
    &=
        \int_{\mathbb{B}^2_1(0,0)}
            \left\|
                (Hh)(\Phi(x,y))
            \right\|_{g_\varphi}^p
            e^{\varphi(\theta_{(x,y)})}
            d\mathcal{L}^2(x,y)\\
    &\leq
        \max_{\theta\in ]-\pi,\pi]}
        e^{\varphi(\theta)}
        \int_{\mathbb{B}^2_1(0,0)}
        \left\|
                (Hh)(\Phi(x,y))
            \right\|_{g_\varphi}^p
            d\mathcal{L}^2(x,y)
    <
        \infty.
    \end{align}
    We have that $\|(\nabla h)_{w}\|_{g_\varphi} \leq \max_{\theta \in ]-\pi, \pi]}e^{-\varphi(\theta)/2}$ for $w \in \Phi\left(\mathbb{B}^2_1(0,0)\right)$. This implies that $\|\nabla h\|_{g_\varphi} \in L^{p}(\mathbb{S}^2\setminus \{z\}, \mu_\varphi)$, since $h$ is bounded and smooth away from $z$. Combining these estimates, we have $h \in W^{2,p}\left(\mathbb{S}^2\setminus \{z\}, g_\varphi\right)$, for all $1 \leq p < 2$.

    Now suppose there exists $\theta_0 \in ]-\pi, \pi]$ such that $|\varphi'(\theta_0)|>0$. By continuity of $\varphi'$, after replacing $\theta_0$ by a nearby point if necessary, we may assume that $\theta_0 \notin \{0, \frac{\pi}{2}, \pi, -\frac{\pi}{2}\}$. Hence, there exists $\sigma > 0$ such that
    \begin{equation}
        \left|\sin\left(\theta_0\right)\right| > \sigma, \quad\quad
        \left|\varphi'\left(\theta_0\right)\right| > \sigma.
    \end{equation}
    By continuity of these functions, there exists $\epsilon >0$ such that $]\theta_0-\epsilon,\theta_0+\epsilon[ \cap \{0, \frac{\pi}{2}, \pi, -\frac{\pi}{2}\} = \emptyset$, and for all $\theta \in ]\theta_0-\epsilon, \theta_0 + \epsilon[$
    \begin{equation}
        \left|\sin\left(\theta\right)\right| > \sigma, \quad\quad
        \left|\varphi'\left(\theta\right)\right| > \sigma.
    \end{equation}

    Equation \eqref{eq:radial2DerivativesOfh} and $\left\|\frac{\partial r}{e^{\varphi\left(\theta_w\right)/2}}\right\|_{g_\varphi} = 1$ imply that, for $w \in \Phi\left(\mathbb{B}^2_1(0,0)\right)$ with $\theta_{w} \in ]\theta_0-\epsilon, \theta_0 + \epsilon[$,
    \begin{align}
        \min_{\theta \in ]-\pi,\pi]}
        e^{-\varphi(\theta)}
        \frac{\sigma^2}
        {2r}
    <&
        \frac{\left|\varphi'\left(\theta_w\right) \sin\left(\theta_{w}\right)\right|}
        {2re^{\varphi(\theta_w)}}
    =
        \left|
        \nabla_{\frac{\partial r}{e^{\varphi(\theta_w)/2}}}
        \nabla_{\frac{\partial r}{e^{\varphi(\theta_w)/2}}}
        h\left(w\right)
        \right|\\
    \leq&
        \left\|
            Hh\left(w\right)
        \right\|_{g_\varphi}.
    \end{align}
    We have $c_h = \min_{\theta \in ]-\pi,\pi]}
        e^{-\varphi(\theta)}
        \frac{\sigma^2}
        {2} > 0$, and so
        \begin{align}
        \int_{\mathbb{S}^2\setminus \{z\}}
            \|Hh\|_{g_\varphi}^2
            d\mu_\varphi
    \geq&
        \int_{\Phi(\mathbb{B}^2_1(0,0))}
            \|Hh\|_{g_\varphi}^2
            d\mu_\varphi\\
    =&
        \int_{\mathbb{B}^2_1(0,0)}
            \left\|(Hh)\left(\Phi\left(x,y\right)\right)\right\|_{g_\varphi}^2
            e^{\varphi(\theta(x,y))}
            d\mathcal{L}^2(x,y)\\
    \geq&
        \min_{\theta \in ]-\pi,\pi]}
            e^{\varphi(\theta)}
            \int_{\mathbb{B}^2_1(0,0) \cap \{\theta(x,y) \in ]\theta_0-\epsilon,\theta_0+\epsilon[\}}
            \|(Hh)(\Phi(x,y))\|_{g_\varphi}^2
            d\mathcal{L}^2(x,y)\\
    =&
        \min_{\theta \in ]-\pi,\pi]}
            e^{\varphi(\theta)}
        \int_{\theta_0-\epsilon}^{\theta_0+\epsilon}
        \int_0^1
            \|(Hh)\left(\Phi\left(r\cos(\theta),r\sin(\theta)
            \right)\right)\|_{g_\varphi}^2
            rdr
            d\theta\\
    >&
        \min_{\theta \in ]-\pi,\pi]}
            e^{\varphi(\theta)}
        \int_{\theta_0-\epsilon}^{\theta_0+\epsilon}
        \int_0^1
            \frac{c_h^2}{r^2}
            rdr
            d\theta
        =
        \infty.
    \end{align}
    This implies that $\|Hh\|_{g_\varphi} \notin L^2\left(\mathbb{S}^2\setminus \{z\}, \mu_\varphi\right)$. Since the measure is finite, for all $p \geq 2$ we have $\left\|Hh\right\|_{g_\varphi} \notin L^p\left(\mathbb{S}^2\setminus \{z\}, \mu_\varphi\right)$. Thus for all $p \geq 2$ we have $h \notin W^{2,p}\left(\mathbb{S}^2 \setminus \{z\}, g_\varphi\right)$, finishing the proof.
\end{proof}

\subsection{Smooth approximating metrics}

Choose a function $\phi \in C^\infty\left(\left[0,1\right]\right)$ such that, for all integers $m \geq 1$,

\begin{align}
    \label{eq:definitionOfPhiFunc}
        \phi(0) &= 0,\quad
        \phi(1) = 1,\quad
        0
        \leq
        \phi(x)
        \leq
        1,\quad
        \phi'(x)> 0 \quad \text{for all } x \in ]0,1[,\\
        \partial_x^{m}\phi(0)
    &=
        \partial_x^{m}\phi(1)
    =
        0.
\end{align}

Given parameters $\eta > 0$ and $0<\epsilon<\frac{\pi}{2}$, define the function $\varphi_{\eta,\epsilon} \in C^\infty([-\pi,\pi])$ by
\begin{equation}
\label{eq:defOfConformalFactor}
    \varphi_{\eta,\epsilon}(\theta)
    =
        \begin{cases}
          \eta  &  \text{ if } |\theta| < \frac{\pi}{2}-\epsilon

          \\
          \eta-\eta\phi\left(\frac{|\theta|-(\frac{\pi}{2}-\epsilon)}{\epsilon}\right) & \text{ if }
          \frac{\pi}{2}-\epsilon
          \leq |\theta| \leq \frac{\pi}{2}\\
            0 & \text{ if } |\theta|\geq \frac{\pi}{2}
        \end{cases}
\end{equation}

Let $\dist_{\eta,\epsilon}$ denote the distance on $\mathbb{B}^2_2(0,0) \setminus \{(0,0)\}$ induced by the metric $e^{\varphi_{\eta,\epsilon}(\theta)}g_\euc$ through the infimum of the length of curves. The completion of $\left(\mathbb{B}^2_2(0,0) \setminus \{(0,0)\}, \dist_{\eta,\epsilon}\right)$ is $\left(\overline{\mathbb{B}^2_2(0,0)}, \overline{\dist}_{\eta,\epsilon}\right)$.

As above, consider the metric $g_{\eta,\epsilon} := g_{\varphi_{\eta,\epsilon}}$ on $\mathbb{S}^2\setminus \{z\}$, defined by \eqref{eq:metricOnSphere}, with distance $\dist_{\eta,\epsilon,\mathbb{S}^2}:= \dist_{\varphi_{\eta,\epsilon}}$ and measure $\mu_{\eta,\epsilon} := \mu_{\varphi_{\eta,\epsilon}}$. Consider the completion $\left(\mathbb{S}^2, \overline{\dist}_{\eta,\epsilon, \mathbb{
S}^2}\right)$. For $p \in \mathbb{S}^2$, let $\tilde{B}^{\eta,\epsilon}_r(p) \subset \mathbb{S}^2$ denote the ball of radius $r$ centered at $p$ with respect to $\overline{\dist}_{\eta,\epsilon}$. Let $B^{\eta,\epsilon}_r(x,y) \subset \mathbb{B}^2_2(0,0)$, denote the ball of radius $r$ centered at $(x,y)$ with respect to $\overline{\dist}_{\eta,\epsilon}$. In particular, there exists $r > 0$ such that, for all $0<\epsilon<\frac{\pi}{2}$,
\begin{equation}
\label{eq:ChartBallStuff}
    \Phi^{-1}\left(
        \tilde{B}^{\eta,\epsilon}_r(z)
    \right)
    =
    B^{\eta,\epsilon}_r(0,0).
\end{equation}

We study this metric for fixed $\eta$, and show that, as $\epsilon \rightarrow 0$, the balls centered at the origin converge to those of the discontinuous metric \eqref{eq:discontinuousMetricDef}. This will allow us to prove an analogue of Lemma \ref{lem:integralOfXDiscontMetric} for $\overline{\dist}_{\eta,\epsilon}$. We therefore investigate the balls and geodesics of this metric, particularly their behavior as $\epsilon \rightarrow 0$.

\begin{lemma}
\label{lem:characterizationOfBallsZero}
    Consider the distance $\overline{\dist}_{\eta,\epsilon}$ defined on $\mathbb{B}^2_2(0)$. For all $0<r<\frac{1}{2}$ and $\alpha>0$ such that $\alpha r < \frac{1}{2}$, we have that the balls $B_r^{\eta,\epsilon}(0,0)$ for the metric $\overline{\dist}_{\eta,\epsilon}$ satisfy
    \begin{align}
    \label{eq:scalingBalls}
        B_{\alpha r}^{\eta,\epsilon}(0,0)
    &=
        \alpha B_{r}^{\eta,\epsilon}(0,0),\\
        B_r^{\eta,\epsilon}(0,0)\cap\{x<0\}
    &=
        \mathbb{B}^2_r(0,0)
        \cap \left\{x<0\right\}.
    \label{eq:leftSideOfBall}
    \end{align}
    Moreover, for any $(x,y) \in \mathbb{B}^2_1(0,0)$ and $0<r<\frac{1}{2}$, we have
    \begin{equation}
    \label{eq:inclusionOfBALLS}
        \mathbb{B}^2_{r/e^{\eta/2}}(x,y)
    \subset
        B_r^{\eta,\epsilon}(x,y)
    \subset
        \mathbb{B}^2_r(x,y)
    \end{equation}
\end{lemma}
\begin{proof}
    To prove \eqref{eq:scalingBalls}, suppose that $0<r<\frac{1}{2}$ and $(x,y) \in B^{\eta,\epsilon}_r(0,0)$. Choose almost-minimizing curves $\gamma_{\delta} : [0,1] \rightarrow \mathbb{B}^2_2(0,0)$, that start at $(x,y)$, end at $(0,0)$ and satisfy
    \begin{equation}
        l_\euc\left(
            \gamma_{\delta}([0,1])
        \right)
        \leq
        \int_0^1
        g_{\eta, \epsilon}
        \left(
            \gamma_{\delta}'(t),
            \gamma_{\delta}'(t)
        \right)^\frac{1}{2}
        dt
    \leq
        \overline{\dist}_{\eta,\epsilon}
        \left(\left(x,y\right),\left(0,0\right)\right) + \delta.
    \end{equation}
    This implies that for $0<\delta<\frac{1}{4}$, we have $\gamma_\delta([0,1]) \subset \mathbb{B}^2_{r+\delta}(0,0) \subset \mathbb{B}^2_2(0,0)$.

    Given $\alpha > 0$ such that $\alpha r < \frac{1}{2}$, choose $\delta>0$ such that $\alpha\delta < \frac{1}{4}$ and consider the curves $\gamma_{\alpha,\delta}(t) = \alpha \gamma_\delta(t)$. Notice the curves satisfy $\gamma_{\alpha,\delta}([0,1]) \subset \mathbb{B}^2_2(0,0)$, begin at $\alpha(x,y)$, and end at $(0,0)$. Furthermore
    \begin{align}
        \int_0^1
            \left\|
                \gamma_{\alpha,\delta}'(t)
            \right\|
            e^{\frac{\varphi_{\eta,\epsilon}(\theta(\gamma_{\alpha,\delta}(t)))}{2}}
        dt
    =&
        \alpha
        \int_0^1
            \left\|
                \gamma_{\delta}'(t)
            \right\|
            e^{\frac{\varphi_{\eta,\epsilon}(\theta(\gamma_\delta(t)))}{2}}
        dt\\
    \leq&
        \alpha
        \left(
        \overline{\dist}_{\eta,\epsilon}
        \left(\left(x,y\right),\left(0,0\right)\right)
        +\delta
        \right).
    \end{align}
    Taking the limit as $\delta$ goes to zero shows
    \begin{equation}
        \overline{\dist}_{\eta,\epsilon}
        \left(\alpha \left(x,y\right),\left(0,0\right)\right)
        \leq
        \alpha
        \overline{\dist}_{\eta,\epsilon}
        \left(
            \left(x,y\right),
            \left(0,0\right)
        \right)
    \end{equation}
    Conversely, we may apply the same result with the factor $1/\alpha$. Since $\frac{1}{\alpha}\left(\alpha r\right) < \frac{1}{2}$ and $\alpha r < \frac{1}{2}$ we also have
    \begin{equation}
        \overline{\dist}_{\eta,\epsilon}\left(
            (x,y)
            ,(0,0)
        \right)
        =
        \overline{\dist}_{\eta,\epsilon}
        \left(
        \frac{1}{\alpha}
        \left(
        \alpha
            (x,y)
        \right)
            ,(0,0)
        \right)
        \leq
        \frac{1}{\alpha}
        \overline{\dist}_{\eta,\epsilon}
        \left(
        \alpha
            (x,y)
            ,(0,0)
        \right),
    \end{equation}
    which proves \eqref{eq:scalingBalls}.

    If $(x,y) \in \mathbb{B}^2_1(0,0)$ satisfies $x<0$, the curve $\gamma(t) = t(x,y)$ gives
    \begin{align}
        \|(x,y)\|
    &\leq
        \overline{\dist}_{\eta, \epsilon}
        \left(
            (x,y), (0,0)
        \right)
    \leq
        \int_0^1
            g_{\eta,\epsilon}
            \left(
                \gamma'(t)
                ,
                \gamma'(t)
            \right)^\frac{1}{2}
            dt\\
    &=
        \int_0^1
            \left\|
                \gamma'(t)
            \right\|
            dt
    =
        \|(x,y)\|,
    \end{align}
    proving equality \eqref{eq:leftSideOfBall}. Inclusions \eqref{eq:inclusionOfBALLS} follow from
    \begin{equation}
        \|v\|^2
    \leq
        g_{\eta,\epsilon}(v,v)
    =
        e^{\varphi_{\eta,\epsilon}(\theta)}
        g_\euc(v,v)
    \leq
        e^\eta\|v\|^2.
    \end{equation}
\end{proof}

Define the sets
\begin{align}
\label{eq:defOfAepsilon}
    A_\epsilon
:=&
    \{r(\cos(\theta),\sin(\theta))\in\mathbb{R}^2:0\leq r, |\theta| \leq \frac{\pi}{2}-\epsilon\}\\
    C_\epsilon
:=&
    \{x\geq 0\} \setminus \text{int}(A_\epsilon).
    \label{eq:defOfCepsilon}
\end{align}

\begin{prop}
\label{prop:convergenceOfBalls}
    Fix $\eta > 0$. For all $(x,y) \in \mathbb{B}^2_1(0,0)$, we have
    \begin{equation}
        \overline{\dist}_{\eta,\epsilon}\left(\left(x,y\right),\left(0,0\right)\right)
        \rightarrow
        \overline{\dist}_\eta\left(\left(x,y\right),\left(0,0\right)\right)
    \end{equation}
    as $\epsilon \rightarrow 0$. Consequently, for every $0<r<\frac{1}{2}$,
    \begin{equation}
        \chi_{B^{\eta,\epsilon}_r(0,0)}(x,y)
    \rightarrow
        \chi_{B^{\eta}_r(0,0)}(x,y)
    \end{equation}
    for $\mathcal{L}^2$-a.e. $(x,y)$.
\end{prop}
\begin{proof}
    By \eqref{eq:discontinuousMetricDef}, \eqref{eq:definitionOfPhiFunc}, and \eqref{eq:defOfConformalFactor}, every curve in $\mathbb{B}^2_2(0,0)$ satisfies
    \begin{align}
        \int_0^1
            g_{\eta,\epsilon}
            \left(
                \gamma'(t),
                \gamma'(t)
            \right)^\frac{1}{2}
            dt
        \leq
            \int_0^1
            g_{\eta}
            \left(
                \gamma'(t),
                \gamma'(t)
            \right)^\frac{1}{2}
            dt,
    \end{align}
    and so in general we have
    \begin{equation}
    \label{eq:basicDistIneq}
    \overline{\dist}_{\eta,\epsilon}\left(\left(x,y\right),\left(0,0\right)\right)
    \leq
        \overline{\dist}_\eta\left(\left(x,y\right),\left(0,0\right)\right).
    \end{equation}

    For the proof we can assume without loss of generality that $(x,y)$ satisfy $y \geq 0$ and $x > 0$ since $\overline{\dist}_\eta\left((x,y),(0,0)\right) = \overline{\dist}_{\eta,\epsilon}\left(
        (x,y),(0,0)
    \right)$ for $x \leq 0$ and in general \begin{equation}
    \overline{\dist}_{\eta,\epsilon}\left(
        (x,y), (0,0)
    \right) = \overline{\dist}_{\eta,\epsilon}\left(
        (x,-y), (0,0)
    \right),\quad
    \overline{\dist}_{\eta}\left(
        (x,y), (0,0)
    \right) = \overline{\dist}_{\eta}\left(
        (x,-y), (0,0)
    \right).
    \end{equation}

    Fix $(x,y) \in \mathbb{B}^2_1(0,0)$ such that $x>0$ and $y \geq 0$, and consider for each $0<\epsilon <\min\left\{
        \frac{\pi}{2}-\theta_{\left(x,y\right)},
        \frac{\left\|
            (x,y)
        \right\|}{e^{\eta/2}}
    \right\}$
    a curve $\gamma_{\eta,\epsilon}:[0,1]\rightarrow \mathbb{B}^2_2(0,0)$ between $(x,y)$ and $(0,0)$ such that
    \begin{equation}
    \label{eq:ineqWithCurveEpsilon}
    \int_{0}^{1}
        g_{\eta,\epsilon}(\gamma_{\eta,\epsilon}'(t),\gamma_{\eta,\epsilon}'(t))^\frac{1}{2}dt
    <
        \dist_{\eta,\epsilon}
        \left((x,y),(0,0)\right)
        +
        \epsilon.
    \end{equation}
    Since the curve $t \in [0,1] \mapsto \left(\left(\gamma_{\eta,\epsilon}\right)_x(t),\left|\left(\gamma_{\eta,\epsilon}\right)_y(t)\right|\right)$ is absolutely continuous and has the same length as $t \in [0,1] \mapsto \gamma_{\eta,\epsilon}(t)$, we may assume that $\left(\gamma_{\eta,\epsilon}\right)_y(t) \geq 0$ for all $t \in [0,1]$.

     The choice of $\epsilon$ implies by \eqref{eq:inclusionOfBALLS} that
    \begin{equation}
    \label{eq:integralOfCurveBoundByDist}
        \int_{0}^{1}
        g_{\eta,\epsilon}(\gamma_{\eta,\epsilon}'(t),\gamma_{\eta,\epsilon}'(t))^\frac{1}{2}dt
    <
        2\overline{\dist}_{\eta,\epsilon}
        \left(
            (x,y),(0,0)
        \right).
    \end{equation}

    Let $0 \leq t_{1,\epsilon} \leq 1$ be the first time such that $\gamma_{\eta,\epsilon}(t_{1,\epsilon}) \in C_\epsilon$ defined by \eqref{eq:defOfCepsilon}. Since $C_\epsilon$ is closed, and by the choice of $\epsilon$ the initial point lies outside it, then $t_{1,\epsilon}>0$. Notice that by \eqref{eq:inclusionOfBALLS} and \eqref{eq:integralOfCurveBoundByDist}, we have for all $s \in [0,1]$
    \begin{align}
         (\gamma_{\eta,\epsilon})_r(s)
    &=
        \|\gamma_{\eta,\epsilon}(s)\|
    \leq
        \overline{\dist}_{\eta,\epsilon}(\gamma_{\eta,\epsilon}(s), (0,0))
    \leq
        \int_{s}^{1}
        g_{\eta,\epsilon}
        \left(
            \gamma_{\eta,\epsilon}'(t),
            \gamma_{\eta,\epsilon}'(t)
        \right)^\frac{1}{2}
        dt\\
    &\leq
        2\overline{\dist}_{\eta,\epsilon}
        \left(
            \left(x,y\right),
            (0,0)
        \right)
    \leq
        2e^{\eta/2}\left\|(x,y)\right\|.
    \label{eq:ineqForRadiusOfCurve}
    \end{align}

    If $t_{1,\epsilon} = 1$, then $\gamma_{\eta,\epsilon}\left([0,1]\right) \subset A_\epsilon$ and so
    \begin{align}
        \overline{\dist}_{\eta}((x,y),(0,0))
    &\leq
        \int_0^1
        g_{\eta}
        \left(
            \gamma_{\eta,\epsilon}'(t)
            ,
            \gamma_{\eta,\epsilon}'(t)
        \right)^\frac{1}{2}
        dt
    =
        \int_0^1
        g_{\eta,\epsilon}
        \left(
            \gamma_{\eta,\epsilon}'(t)
            ,
            \gamma_{\eta,\epsilon}'(t)
        \right)^\frac{1}{2}
        dt\\
    &\leq
        \overline{\dist}_{\eta,\epsilon}((x,y),(0,0))
        + \epsilon.
        \label{eq:ineqCase1}
    \end{align}

    On the other hand, if $t_{1,\epsilon} < 1$, choose $t_{2,\epsilon}>0$ such that $t_{1,\epsilon} <t_{2,\epsilon}<1$. Let $\sigma_{2,\epsilon}:[t_{1,\epsilon}, t_{2,\epsilon}] \rightarrow \mathbb{B}^2_2(0,0)$ be the straight line between $\gamma_{\eta,\epsilon}(t_{1,\epsilon})$ and $(0,\|\gamma_{\eta,\epsilon}(t_{1,\epsilon})\|)$ and let $\sigma_{3,\epsilon}:[t_{2,\epsilon}, 1] \rightarrow \mathbb{B}^2_2(0,0)$ be the straight line between $(0,\|\gamma_{\eta,\epsilon}(t_{1,\epsilon})\|)$ and $(0,0)$. Define the curve $\sigma_{\epsilon}:[0,1] \rightarrow \mathbb{B}^2_2(0,0)$ by
    \begin{equation}
        \sigma_{\epsilon}(t)
    :=
        \begin{cases}
          \gamma_{\eta,\epsilon}(t)  & t \in [0,{t_{1,\epsilon}}] \\
          \sigma_{2,\epsilon}\left(t
          \right)& t \in ]t_{1,\epsilon},t_{2,\epsilon}]\\
          \sigma_{3,\epsilon}
          (t) &
          t \in ]{t}_{2,\epsilon}, 1]
        \end{cases}
    \end{equation}
    The points $\sigma_\epsilon(t_{1,\epsilon})$ and $\sigma_\epsilon(t_{2,\epsilon})$ have the same radius and respective angles $\frac{\pi}{2}-\epsilon$ and $\frac{\pi}{2}$. Since $\sigma_{2,\epsilon}$ is the straight line between points that lie on the same circle and differ by angle $\epsilon$, \eqref{eq:ineqForRadiusOfCurve} gives
    \begin{equation}
        l_{\euc}\left(\sigma_{2,\epsilon}\left([t_{1,\epsilon},t_{2,\epsilon}]\right)
        \right)
    \leq
        \|\gamma_{\eta,\epsilon}(t_{1,\epsilon})\|
        \epsilon
    \leq
        2e^{\eta/2}
        \|(x,y)\|\epsilon.
    \end{equation}
    This implies that
    \begin{equation}
        \int_{t_{1,\epsilon}}^{t_{2,\epsilon}}
        g_{\eta}(\sigma_{2,\epsilon}'(t),\sigma_{2,\epsilon}'(t))^\frac{1}{2}dt
    =
        e^{\eta/2}
        l_{\euc}\left(
            \sigma_{2,\epsilon}\left(
            [t_{1,\epsilon},t_{2,\epsilon}]\right)
        \right)
    \leq
        2e^{\eta}\|(x,y)\|\epsilon.
    \label{eq:ineq1GeoStuff}
    \end{equation}
    On the other hand, since $\gamma_{\eta,\epsilon}(1) = (0,0)$ and $\varphi_{\eta,\epsilon} \geq 0$, we also have
    \begin{align}
    \int_{t_{1,\epsilon}}^1
        \|\gamma_{\eta,\epsilon}'(t)\|
        e^{\frac{\varphi_{\eta, \epsilon}(\theta(\gamma_{\eta,\epsilon}(t)))}{2}}
        dt
    &\geq
        l_\euc\left(\gamma_{\eta,\epsilon}\left(\left[t_{1,\epsilon}, 1\right]\right)\right)\\
    &\geq
        \|\gamma_{\eta,\epsilon}(t_{1,\epsilon})\|\\
    &=
        l_\euc\left(\sigma_{3,\epsilon}\left([t_{2,\epsilon}, 1]\right)\right)\\
    &=
        \int_{t_{2,\epsilon}}^1
            g_\eta
            \left(
                \sigma_{3,\epsilon}'(t),
                \sigma_{3,\epsilon}'(t)
            \right)^\frac{1}{2}
            dt.
        \label{eq:ineq2GeoStuff}
    \end{align}
    Thus, combining \eqref{eq:ineqWithCurveEpsilon}, \eqref{eq:ineq1GeoStuff}, and \eqref{eq:ineq2GeoStuff}, we obtain
    \begin{align}
        \epsilon+
        \overline{\dist}_{\eta,\epsilon}\left(\left(x,y\right),\left(0,0\right)\right)
    \geq&\int_{0}^1
        \|\gamma_{\eta,\epsilon}'(t)\|
        e^{\frac{\varphi_{\eta, \epsilon}(\theta(\gamma_{\eta,\epsilon}(t)))}{2}}
        dt\\
    \geq&
        \int_0^{t_{1,\epsilon}}
        \|\gamma_{\eta,\epsilon}'(t)\|
        e^{\frac{\varphi_{\eta, \epsilon}(\theta(\gamma_{\eta,\epsilon}(t)))}{2}}
        dt
        +
        \int_{t_{1,\epsilon}}^1
        \|\gamma_{\eta,\epsilon}'(t)\|
        e^{\frac{\varphi_{\eta, \epsilon}(\theta(\gamma_{\eta,\epsilon}(t)))}{2}}
        dt\\
    =&\int_0^{t_{1,\epsilon}}
        \|\sigma_{\epsilon}'(t)\|
        e^{\eta/2}
        dt
        +
        \int_{t_{1,\epsilon}}^1
        \|\gamma_{\eta,\epsilon}'(t)\|
        e^{\frac{\varphi_{\eta, \epsilon}(\theta(\gamma_{\eta,\epsilon}(t)))}{2}}
        dt\\
    \geq&
        \int_0^{t_{1,\epsilon}}
        g_\eta(\sigma_{\epsilon}'(t),\sigma_{\epsilon}'(t))^\frac{1}{2}
        dt
    +
        \int_{t_{1,\epsilon}}^{t_{2,\epsilon}}
        g_{\eta}(\sigma_{2,\epsilon}'(t),\sigma_{2,\epsilon}'(t))^\frac{1}{2}
        dt
    -
        2e^{\eta}\|(x,y)\|\epsilon\\
    &\quad\quad
        +\int_{t_{2,\epsilon}}^1
            g_\eta
            \left(
                \sigma_{3,\epsilon}'(t),
                \sigma_{3,\epsilon}'(t)
            \right)^\frac{1}{2}
            dt\\
    =&
        \int_{0}^1
            g_\eta
            \left(
                \sigma_{\epsilon}'(t),
                \sigma_{\epsilon}'(t)
            \right)^\frac{1}{2}
            dt
    -
        2e^{\eta}\|(x,y)\|\epsilon\\
    \geq&
        \overline{\dist}_\eta\left(\left(x,y\right),\left(0,0\right)\right)
    -
        2e^\eta\|(x,y)\|\epsilon.
        \label{eq:ineqCase2}
    \end{align}
    Thus, as $\epsilon \rightarrow 0$, \eqref{eq:basicDistIneq}, \eqref{eq:ineqCase1}, and \eqref{eq:ineqCase2} imply
    \begin{equation}
        \overline{\dist}_{\eta,\epsilon}\left(\left(x,y\right),\left(0,0\right)\right)
    \rightarrow
    \overline{\dist}_{\eta}\left(\left(x,y\right),\left(0,0\right)\right).
    \end{equation}
    For the convergence of the characteristic functions, observe that for $\overline{\dist}_{\eta}\left((x,y),(0,0)\right) \neq r$ we have
    $\chi_{B^{\eta,\epsilon}_r(0,0)}(x,y) \rightarrow \chi_{B^{\eta}_r(0,0)}(x,y)$. By Lemma \ref{lem:ballsOfDiscontinuousMetric}, the set $\{(x,y)\in \mathbb{R}^2 : \overline{\dist}_{\eta}\left((x,y),(0,0)\right) = r\}$ has zero measure for $0<r<\frac{1}{2}$, which finishes the proof.
\end{proof}

\subsection{Non-convergence of AMV}

\begin{lemma}
\label{lem:negativeIntXContinuousMetric}
    There exist $\eta>0$, $\epsilon>0$, $a>0$, and $\delta > 0$ such that, for all $0<r<a$,
    \begin{equation}
        \int_{\tilde{B}^{\eta,\epsilon}_r(z)}
            h
        d\mu_{\eta,\epsilon}
    \leq
        -\delta r^3.
    \end{equation}
\end{lemma}
\begin{proof}

    By the metric bounds \eqref{eq:metricBoundsSphere}, and $0\leq \varphi_{\eta,\epsilon} \leq \eta$, there exists $0<a<\frac{1}{2}$  such that $\tilde{B}^{\eta,\epsilon}_a(z) \subset \Phi\left(\mathbb{B}^2_{\frac{1}{2}}(0,0)\right)$, and by \eqref{eq:metricOnSphere},  $B_a^{\eta,\epsilon}(0,0) = \Phi^{-1}\left(\tilde{B}^{\eta,\epsilon}_a(z)\right)$ for all $0<\epsilon< \frac{\pi}{2}$.

    Since $f(x,y):=h \circ \Phi(x,y) = x$ is linear, Lemma \ref{lem:characterizationOfBallsZero} and a change of variables give
    \begin{align}
        \int_{\tilde{B}^{\eta,\epsilon}_r(z)}
            h
        d\mu_{\eta,\epsilon}
    =&
        \int_{B^{\eta,\epsilon}_r(0,0)}
            h\circ \Phi(x,y)
            e^{\varphi_{\eta,\epsilon}(\theta(x,y))}
        d\mathcal{L}^2(x,y)\\
    =&
         \int_{B^{\eta,\epsilon}_r(0,0)}
            x
            e^{\varphi_{\eta,\epsilon}(\theta(x,y))}
        d\mathcal{L}^2(x,y)\\
    =&
        \frac{r^3}{a^3}
        \int_{B^{\eta,\epsilon}_a(0,0)}
            x
            e^{\varphi_{\eta,\epsilon}(\theta(x,y))}
        d\mathcal{L}^2(x,y)\\
    =&
        \frac{r^3}{a^3}
        \int_{B^{\eta,\epsilon}_a(0,0)}
            f(x,y)
            e^{\varphi_{\eta,\epsilon}(\theta(x,y))}
        d\mathcal{L}^2(x,y)\\
    =&
        \frac{r^3}{a^3}
        \int
            \chi_{B^{\eta,\epsilon}_a(0,0)}(x,y)
            f(x,y)
            e^{\varphi_{\eta,\epsilon}(\theta(x,y))}
            d\mathcal{L}^2(x,y).
    \label{eq:equalityOfFwithH}
    \end{align}
    By Proposition \ref{prop:convergenceOfBalls} and the definition of $\varphi_{\eta,\epsilon}$, we conclude
    \begin{equation}
        \chi_{B^{\eta,\epsilon}_a(0,0)}(x,y)
            f(x,y)
            e^{\varphi_{\eta,\epsilon}(\theta(x,y))}
        \rightarrow
            \chi_{B^{\eta}_a(0,0)}(x,y)
            f(x,y)
            G_\eta(x,y),
    \end{equation}
    as $\epsilon \rightarrow 0$ for $\mathcal{L}^2$-a.e. $(x,y)$. Since $B_{a}^{\eta,\epsilon}(0,0) \subset \mathbb{B}^2_1(0,0)$, we have
    \begin{equation}
    \left|
        \chi_{B^{\eta,\epsilon}_a(0,0)}(x,y)
            f(x,y)
            e^{\varphi_{\eta,\epsilon}(\theta(x,y))}
    \right|
        \leq
        e^{\eta}
        \chi_{\mathbb{B}^2_1(0,0)}(x,y).
    \end{equation}
    Therefore, by the dominated convergence theorem,
    \begin{align}
        \int
            \chi_{B^{\eta,\epsilon}_a(0,0)}(x,y)
            f(x,y)
            e^{\varphi_{\eta,\epsilon}(\theta(x,y))}
            d\mathcal{L}^2(x,y)
    \rightarrow&
        \int
        \chi_{B^{\eta}_a(0,0)}(x,y)
            f(x,y)
            G_\eta(x,y)
            d\mathcal{L}^2(x,y)\\
        =&
        \int_{B_{a}^\eta(0,0)}
        f\,d\mu_\eta.
        \label{eq:integralConvergence}
    \end{align}
    Let $\tilde{\delta}, \eta > 0$ be the constants from Lemma \ref{lem:integralOfXDiscontMetric}. By the previous convergence, there exists $\epsilon >0$ small enough that
    \begin{equation}
    \label{eq:approximationOfIntegralByMetric}
        \left|
            \int_{B^{\eta,\epsilon}_a(0,0)}
            f(x,y)
            e^{\varphi_{\eta,\epsilon}(\theta(x,y))}
        d\mathcal{L}^2(x,y)
        -
            \int_{B^{\eta}_a(0,0)}
            f(x,y)
            G_\eta(x,y)
        d\mathcal{L}^2(x,y)
        \right|
        \leq
        \frac{a^3\tilde{\delta}}{2}.
    \end{equation}
    Using Lemma \ref{lem:integralOfXDiscontMetric} and \eqref{eq:approximationOfIntegralByMetric}, we conclude that
    \begin{equation}
        \int_{\tilde{B}^{\eta,\epsilon}_r(z)}
            h
        d\mu_{\eta,\epsilon}
    \leq
        -\frac{\tilde{\delta} r^3}{2},
    \end{equation}
    and by taking $\delta = \frac{\tilde{\delta}}{2}$, this finishes the proof.
\end{proof}

\begin{prop}
\label{prop:AMVis1OverR}
    Fix $\eta, \epsilon > 0$ constants from Lemma \ref{lem:negativeIntXContinuousMetric}. Consider the metric measure space $(\mathbb{S}^2, \overline{\dist}_{\eta,\epsilon, \mathbb{S}^2}, \mu_{\eta,\epsilon})$. If $\Delta_r^{\eta,\epsilon}$ is the AMV Laplace operator of this metric measure space, then there exist $\sigma >0, r_0>0$ and $\tilde{\delta} > 0$ such that for all $0<r<r_0$ and $(x,y) \in \mathbb{B}^2_{\sigma r}(0,0)$
    \begin{equation}
        \left|
            \Delta_r^{\eta,\epsilon}
            h
            (\Phi(x,y))
        \right|
    \geq
        \frac{\tilde{\delta}}{r}.
    \end{equation}
\end{prop}

\begin{proof}
By Lemmas \ref{lem:characterizationOfBallsZero} and \ref{lem:negativeIntXContinuousMetric}, there exists $0<r_0<\frac{1}{2}$ such that for all $0<\sigma < \min\left\{\frac{1}{2},\frac{1}{2e^{\eta}}\right\}$, $0<r<r_0$ and  $(x,y) \in \mathbb{B}^2_{\sigma r}(0,0)$, we have $\tilde{B}_r^{\eta,\epsilon}\left(\Phi(x,y)\right) \subset \Phi\left(\mathbb{B}^2_1(0)\right)$, and
    \begin{align}
        \mu_{\eta,\epsilon}(\tilde{B}^{\eta,\epsilon}_r(\Phi(x,y)))
    &=\mu_{\eta,\epsilon }\left(\Phi\left(B^{\eta,\epsilon}_r\left(x,y\right)\right)\right)\\
    &\leq
        \mu_{\eta,\epsilon}\left(\Phi\left(\mathbb{B}^2_{r}(x,y)\right)\right)\\
    &\leq
        r^2
        e^\eta
        \mathcal{L}^2\left(
            \mathbb{B}^2_1(0,0)
        \right).
    \end{align}
    Shrink $r_0>0$ so that $r_0 < \min\left\{\frac{a}{1+e^{\eta/2}}, \frac{1}{2\left(1+e^{\eta/2}\right)} \right\} \leq \min\left\{\frac{a}{1+e^{\eta/2} \sigma }, \frac{1}{2\left(1+e^{\eta/2}\sigma \right)} \right\}$ where $a$ is the constant from Lemma \ref{lem:negativeIntXContinuousMetric}.

    Since $(x,y) \in \mathbb{B}^2_{\sigma r}(0,0)$, $\overline{\dist}_{\eta,\epsilon}\left((x,y),(0,0)\right) \leq e^{\eta/2}\|(x,y)\| \leq e^{\eta/2}\sigma r$. Thus, the triangle inequality gives
    \begin{equation}
        B^{\eta,\epsilon}_{r-e^{\eta/2}\sigma r}(0,0)
    \subset
        B^{\eta,\epsilon}_r(x,y)
    \subset
        B^{\eta,\epsilon}_{r+e^{\eta/2}\sigma r}(0,0).
    \end{equation}
    By Lemma \ref{lem:characterizationOfBallsZero} we have
    \begin{align}
        \mathcal{L}^2
        \left(
            B^{\eta,\epsilon}_{r+e^{\eta/2}\sigma r}(0,0)
        \setminus
            B^{\eta,\epsilon}_{r-e^{\eta/2}\sigma r}(0,0)
        \right)
    &=
        \mathcal{L}^2
        \left(
            B^{\eta,\epsilon}_{r+e^{\eta/2}\sigma r}(0,0)
        \setminus
            \frac{r-e^{\eta/2}\sigma r}{r+e^{\eta/2}\sigma r}B^{\eta,\epsilon}_{r+e^{\eta/2}\sigma r}(0,0)
        \right)\\
    &=
        \left(
            1
        -
        \left(
            \frac{r-e^{\eta/2}\sigma r}{r+e^{\eta/2}\sigma r}
        \right)^2
        \right)
        \mathcal{L}^2
        \left(B^{\eta,\epsilon}_{r+e^{\eta/2}\sigma r}(0,0)
        \right)\\
    &\leq
        \frac{4e^{\eta/2}\sigma}{\left(1+e^{\eta/2}\sigma\right)^2}
        r^2\left(1+e^{\eta/2}\sigma\right)^2
        \mathcal{L}^2\left(
            \mathbb{B}^2_1(0,0)
        \right)\\
    &=4e^{\eta/2}\sigma r^2\mathcal{L}^2\left(
            \mathbb{B}^2_1(0,0)
        \right).
    \end{align}
    Also using inclusion \eqref{eq:inclusionOfBALLS}, this implies
    \begin{align}
        &\int_{\left(B^{\eta,\epsilon}_{r+e^{\eta/2}\sigma r}(0,0)
        \setminus
            B^{\eta,\epsilon}_{r-e^{\eta/2}\sigma r}(0,0) \right)
            \cap \{\tilde{x}\geq 0\}}
            \tilde{x}
            e^{\varphi_{\eta,\epsilon}(\theta(\tilde{x},\tilde{y}))}
            d\mathcal{L}^2(\tilde{x},\tilde{y})\\
        \leq&
            r
            e^\eta
            \left(
                1+e^{\eta/2}\sigma
            \right)\mathcal{L}^2
            \left(
            B^{\eta,\epsilon}_{r+e^{\eta/2}\sigma r}(0,0)
        \setminus
            B^{\eta,\epsilon}_{r-e^{\eta/2}\sigma r}(0,0)
            \right)\\
        \leq&
            4r^3
            \sigma e^{3\eta/2}
            \left(
                1+e^{\eta/2}\sigma
            \right)
            \mathcal{L}^2\left(
                \mathbb{B}^2_1(0,0)
            \right).
    \end{align}
    Thus we have
    \begin{align}
        &\left|
        \Delta_r^{\eta,\epsilon}
        h\left(
            \Phi\left(
                x,y
            \right)
        \right)
        \right|\\
    =&
    \left|
        \frac{1}{r^2\mu_{\eta,\epsilon}(\tilde{B}_r^{\eta,\epsilon}(\Phi(x,y)))}
        \int_{\tilde{B}_r^{\eta,\epsilon}(\Phi(x,y))}
        \left(
            h(p)
            -
            h(\Phi(x,y))
        \right)
            d\mu_{\eta,\epsilon}(p)
    \right|\\
    =&
    \left|
        \frac{1}{r^2\mu_{\eta,\epsilon}(\tilde{B}_r^{\eta,\epsilon}(\Phi(x,y)))}
        \int_{B_r^{\eta,\epsilon}(x,y)}
        \left(
            h(\Phi(\tilde{x},\tilde{y}))
            -
            h(\Phi(x,y))
        \right)
        e^{\varphi_{\eta,\epsilon}(\theta(\tilde{x},\tilde{y}))}
            d\mathcal{L}^2(\tilde{x},\tilde{y})
    \right|\\
    \geq&
    \left|
        \frac{1}{r^2\mu_{\eta,\epsilon}(\tilde{B}_r^{\eta,\epsilon}(\Phi(x,y)))}
        \int_{B_r^{\eta,\epsilon}(x,y)}
            h(\Phi(\tilde{x},\tilde{y}))
        e^{\varphi_{\eta,\epsilon}(\theta(\tilde{x},\tilde{y}))}
            d\mathcal{L}^2(\tilde{x},\tilde{y})
    \right|
    -
    \left|\frac{ h(\Phi(x,y))}{r^2} \right|\\
    \geq&
        \left|
        \frac{1}{r^4\mathcal{L}^2(\mathbb{B}^2_1(0,0)) e^\eta}
        \int_{B_r^{\eta,\epsilon}(x,y)}
            h(\Phi(\tilde{x},\tilde{y}))
        e^{\varphi_{\eta,\epsilon}(\theta(\tilde{x},\tilde{y}))}
            d\mathcal{L}^2(\tilde{x},\tilde{y})
    \right|
    -
        \frac{\sigma}{r}.
    \label{eq:amvCalcMiddle}
    \end{align}
    Now decompose the integral as
    \begin{align}
        &\int_{B_r^{\eta,\epsilon}(x,y)}
            h(\Phi(\tilde{x},\tilde{y}))
        e^{\varphi_{\eta,\epsilon}\left(\theta(\tilde{x},\tilde{y})\right)}
            d\mathcal{L}^2(\tilde{x},\tilde{y})\\
    =&
        \int_{B_r^{\eta,\epsilon}(x,y) \cap \{\tilde{x}<0\}}
            \tilde{x}
        e^{\varphi_{\eta,\epsilon}(\theta(\tilde{x},\tilde{y}))}
            d\mathcal{L}^2(\tilde{x},\tilde{y})
        +
        \int_{B_r^{\eta,\epsilon}(x,y) \cap \{\tilde{x}\geq0\}}
            \tilde{x}
        e^{\varphi_{\eta,\epsilon}(\theta(\tilde{x},\tilde{y}))}
            d\mathcal{L}^2(\tilde{x},\tilde{y})\\
    \leq&
        \int_{B^{\eta,\epsilon}_{r-e^{\eta/2}\sigma r}(0,0) \cap \{\tilde{x}<0\}}
            \tilde{x}
        e^{\varphi_{\eta,\epsilon}(\theta(\tilde{x},\tilde{y}))}
            d\mathcal{L}^2(\tilde{x},\tilde{y})
        +
        \int_{B^{\eta,\epsilon}_{r+e^{\eta/2}\sigma r}(0,0) \cap \{\tilde{x}\geq0\}}
            \tilde{x}
        e^{\varphi_{\eta,\epsilon}(\theta(\tilde{x},\tilde{y}))}
            d\mathcal{L}^2(\tilde{x},\tilde{y})\\
    =&
        \int_{B^{\eta,\epsilon}_{r-e^{\eta/2}\sigma r}(0,0)}
            \tilde{x}
        e^{\varphi_{\eta,\epsilon}(\theta(\tilde{x},\tilde{y}))}
            d\mathcal{L}^2(\tilde{x},\tilde{y})
        +
        \int_{
        \left(B^{\eta,\epsilon}_{r+e^{\eta/2}\sigma r}(0,0) \setminus B^{\eta,\epsilon}_{r-e^{\eta/2}\sigma r}(0,0)
        \right)\cap \{\tilde{x}\geq0\}}
            \tilde{x}
        e^{\varphi_{\eta,\epsilon}(\theta(\tilde{x},\tilde{y}))}
            d\mathcal{L}^2(\tilde{x},\tilde{y})\\
    \leq&
        \int_{\tilde{B}^{\eta,\epsilon}_{r-e^{\eta/2}\sigma r}(z)}
            h
            d\mu_{\eta, \epsilon}
        +
            4r^3
            \sigma e^{3\eta/2}
            \left(
                1+e^{\eta/2}\sigma
            \right)
            \mathcal{L}^2\left(
                \mathbb{B}^2_1(0,0)
            \right).
    \end{align}

    Let $\delta > 0$ be as in Lemma \ref{lem:negativeIntXContinuousMetric}. Choose $\sigma$ small enough such that
    \begin{equation}
        4\sigma e^{3\eta/2}\left(
                1+e^{\eta/2}\sigma
            \right)\mathcal{L}^2\left(\mathbb{B}^2_1(0,0)\right) < \frac{\delta}{32},
        \quad 
        e^{\eta/2}\sigma < \frac{1}{2}.
    \end{equation}
    Then
    \begin{align}
        &\int_{B_r^{\eta,\epsilon}(x,y)}
            h(\Phi(\tilde{x},\tilde{y}))
        e^{\varphi_{\eta,\epsilon}(\theta(\tilde{x},\tilde{y}))}
            d\mathcal{L}^2(\tilde{x},\tilde{y})\\
    \leq&
        -\frac{\delta}{2}
        \left(
            r-e^{\eta/2} \sigma r
        \right)^3
        +
            4r^3\sigma e^{3\eta/2}
            \left(
                1+e^{\eta/2}\sigma
            \right)\mathcal{L}^2\left(\mathbb{B}^2_1(0,0)\right)\\
    \leq&
        -\frac{\delta r^3}{16}
        +
        \frac{\delta r^3}{32}
    \leq
        -\frac{\delta r^3}{32}.
        \label{eq:ineqOfIntegralOfh}
    \end{align}
    Substitute \eqref{eq:ineqOfIntegralOfh} into \eqref{eq:amvCalcMiddle}. If necessary, decrease $\sigma$ further so that $\sigma< \frac{\delta}{64 \mathcal{L}^2\left(\mathbb{B}^2_1(0,0)\right)e^\eta}$. We obtain
    \begin{align}
        \left|
        \Delta_r^{\eta,\epsilon}
        h\left(
            \Phi\left(
                x,y
            \right)
        \right)
        \right|
    \geq
        \left|
            \frac{1}{r^4 \mathcal{L}^2\left(\mathbb{B}^2_1(0,0)\right)e^\eta}
            \frac{\delta r^3}{32}
        \right|
        -
        \frac{\sigma}{r}
    \geq
        \frac{\delta}{ 64\mathcal{L}^2\left(\mathbb{B}^2_1(0,0)\right)e^\eta r},
    \end{align}
    finishing the proof by taking $\tilde{\delta} = \frac{\delta}{64\mathcal{L}^2\left(\mathbb{B}^2_1(0,0)\right)e^\eta}$.
\end{proof}

We now prove that on the metric measure space $(\mathbb{S}^2, \overline{\dist}_{\eta,\epsilon, \mathbb{S}^2}, \mu_{\eta,\epsilon})$ the AMV Laplace operator of $h$ does not converge to $\Delta_{g_{\varphi_{\eta,\epsilon}}} h$ in $L^p\left(\mathbb{S}^2, \mu_{\eta,\epsilon}\right)$, even though $\Delta_{g_{\varphi_{\eta,\epsilon}}} h \in L^\infty(\mathbb{S}^2, \mu_{\eta,\epsilon})$ and $h$ is Lipschitz. However, we note that $\|H h\|_{g_{\eta,\epsilon}} \notin L^2(\mathbb{S}^2, \mu_{\eta,\epsilon})$.

\begin{proof}[Proposition \ref{prop:nonConvergenceResult}]
    Set $z = \Phi(0,0)$, and fix $\eta, \epsilon > 0$ from Lemma \ref{lem:negativeIntXContinuousMetric} and Proposition \ref{prop:AMVis1OverR}. Define $g_{\psi} := g_{\varphi_{\eta,\epsilon}}$ as defined by \eqref{eq:metricOnSphere}, \eqref{eq:metricOnSphere2} and \eqref{eq:defOfConformalFactor}. By Proposition \ref{prop:AMVis1OverR} there exists $\sigma, \tilde{\delta}, r_0 > 0$ such that for all $0<r<r_0$ and $(x,y) \in \mathbb{B}^2_{\sigma r}(0,0)$ we have
    \begin{equation}
        \left|
            \Delta_r^{\eta,\epsilon}
            h
            \left(
                \Phi(x,y)
            \right)
        \right|
    \geq
        \frac{\tilde{\delta}}{r}.
    \end{equation}

    Let $r_n \rightarrow 0$ be an arbitrary sequence and set $A_n := \Phi\left(\mathbb{B}^2_{\sigma r_n}(0,0)\right)$. Then $\mu_{\eta, \epsilon}(A_n) \rightarrow 0$. Thus, if $\Delta_{r_n}^{\eta,\epsilon}h$ converges in $L^p(\mathbb{S}^2,\mu_{\eta,\epsilon})$ to some function $g$, we must have
    \begin{equation}
        \left(\int_{A_n}
            \left|
            \Delta_{r_n}^{\eta,\epsilon}
            h
        \right|^p
        d\mu_{\eta,\epsilon}
        \right)^\frac{1}{p}
    \leq
    \left(
        \int_{A_n}
            \left|
            \Delta_{r_n}^{\eta,\epsilon}
            h
            -
            g
        \right|^p
        d\mu_{\eta,\epsilon}
    \right)^\frac{1}{p}
        +
    \left(
        \int_{A_n}
        |g|^p
        d\mu_{\eta,\epsilon}
    \right)^\frac{1}{p}
    \rightarrow
        0,
    \end{equation}
    by absolute continuity of the integral of $|g|^p$. However, if $p \geq 2$, then for sufficiently large $N \in \mathbb{N}$ and every $n \geq N$
    \begin{align}
        \int_{A_n}
            \left|
            \Delta_{r_n}^{\eta,\epsilon}
            h
        \right|^p
        d\mu_{\eta,\epsilon}
    =&
        \int_{\mathbb{B}^2_{\sigma r_n}(0,0)}
            \left|
            \Delta_{r_n}^{\eta,\epsilon}
            h\left(\Phi(x,y)\right)
        \right|^p
        e^{\varphi_{\eta,\epsilon}\left(\theta_{(x,y)}\right)}
        d\mathcal{L}^2(x,y)\\
    \geq&
        \int_{\mathbb{B}^2_{\sigma r_n}(0,0)}
            \frac{\tilde{\delta}^p}{r_n^p}
        d\mathcal{L}^2(x,y)\\
    \geq&
        \tilde{\delta}^p
        \sigma^2
        \frac{r_n^2}{r_n^p}
        \mathcal{L}^2\left(
            \mathbb{B}^2_1(0,0)
        \right)
        \not\rightarrow 0.
    \end{align}
    Thus for $p \geq 2$, $\Delta_{r_n}^{\eta,\epsilon}
            h$ does not converge in $L^p(\mathbb{S}^2, \mu_{\eta,\epsilon})$. The same bound shows that, if $p>2$, then $\left\|\Delta_{r_n}^{\eta,\epsilon} h\right\|_{L^p\left(\mathbb{S}^2,\mu_{\eta,\epsilon}\right)} \rightarrow \infty$.

            By \eqref{eq:inclusionOfBALLS} we have
            \begin{equation}
                \frac{\overline{\dist}_{\eta, \epsilon}\left(
                    (x,y),(0,0)
                \right)}{e^{\eta/2}}
            \leq
                \left\|
                    (x,y)
                \right\|.
            \end{equation}
            There exists a neighborhood $z \in V_z \subset \Phi\left(\mathbb{B}^2_1(0,0)\right)$ such that, for $w \in V_z$, $\overline{\dist}_{\eta,\epsilon, \mathbb{S}^2}\left(w,z\right) = \overline{\dist}_{\eta,\epsilon}\left(
                \Phi^{-1}\left(
                    w
                \right)
                ,(0,0)
            \right)$.
            Therefore, by \eqref{eq:boundOnCurvature}, for $w \in V_z$
            \begin{align}
                K_{\eta,\epsilon}(w)
            \leq&
                \frac{C}{\left\|\Phi^{-1}(w)\right\|^2}
            \leq
                \frac{Ce^{\eta}}{\overline{\dist}_{\eta,\epsilon}\left(
                    \Phi^{-1}(w),
                    (0,0)
                \right)^2}\\
            =&
                \frac{Ce^{\eta}}{\overline{\dist}_{\eta,\epsilon, \mathbb{S}^2}\left(
                    w,
                    z
                \right)^2}.
            \end{align}
            Since $g_{\eta,\epsilon}$ is smooth on $\mathbb{S}^2 \setminus \{z\}$ and $\mathbb{S}^2 \setminus V_z$ is compact, after enlarging $C$ if necessary, for all $w \in \mathbb{S}^2 \setminus V_z$
            \begin{equation}
                K_{\eta,\epsilon}(w)
            \leq
                \frac{Ce^\eta}{\text{diam}_{g_{\eta,\epsilon}}\left(\mathbb{S}^2\right)^2}
            \leq
            \frac{Ce^\eta}{\tilde{\dist}_{\eta,\epsilon}\left(w,z\right)^2},
            \end{equation}
            proving \eqref{eq:curvatureBoundForExample}.

            Notice that $h\circ \Phi(x,y) = x$. Setting $h(z)=0$ makes $h$ smooth on $\mathbb{S}^2$. Since $g_{\euc} \leq g_{\eta,\epsilon}$, the function is locally Lipschitz in a neighborhood $U_z$ of $z$, with $\sup_{w \in U_z}\left\|\nabla h\right\|_{g_{\eta,\epsilon}} \leq 1$. On the compact set $\mathbb{S}^2\setminus U_z$, smoothness gives $\sup_{w \in \mathbb{S}^2\setminus U_z}\left\|\nabla h\right\|_{g_{\eta,\epsilon}}< L$ for some $L>0$. Hence $\sup_{w \in \mathbb{S}^2 \setminus \{z\}}\left\|\nabla h\right\|_{g_{\eta,\epsilon}} \leq \max\left\{1,L\right\}$. Integrating $\nabla h$ along curves shows that $h$ is Lipschitz on $\mathbb{S}^2 \setminus \{z\}$, and continuity at $z$ extends the Lipschitz estimate to $\mathbb{S}^2$.

            By Lemma \ref{lem:hIsHarmonicAndNotSobolev}, $\Delta_{g_{\eta,\epsilon}}h \in L^\infty\left(\mathbb{S}^2 \setminus \{z\},\mu_{\eta,\epsilon}\right)$. Since $\varphi_{\eta,\epsilon}' \neq 0$, we also have $h \notin W^{2,2}\left(\mathbb{S}^2 \setminus \{z\}, g_{\eta,\epsilon}\right)$. This finishes the proof.
\end{proof}

\section{Declaration of generative AI and AI-assisted technologies in the writing process}

The author is responsible for the architecture of the paper and proofs of the main results. During
the preparation of this manuscript, the authors used ChatGPT (Sol 5.6) to assist with editorial
restructuring, proof presentation and review, bibliographical research and \LaTeX consistency.  The author reviewed and edited all suggestions and takes full responsibility for
the content of the article.

\printbibliography

\end{document}